\input epsf

%auto-ignore (forbid arXiv's auto-postscript generator TeXing this file alone)
% Protection

\ifx\macrosloaded\relax    \else \let\macrosloaded\relax \fi

% Make sure reasonable size defaults (for A4 paper) are defined. These can be
% overruled selectively by explicit definitions after including this file, or
% globally by setting them before loading and defining \formloaded

\ifx\formloaded\undefined 

\hsize=16cm
\vsize=24.5 cm \advance \vsize -24pt % space for page numbers

\fi

\catcode`@=11

% TeXnical tricks

\def\letnext{\let\next= }
\def\defnext{\def\next}   % Not \outer, in case \next is.

\newif\ifamstex {\defnext{AmS-TeX} \ifx\fmtname\next \global\amstextrue \fi}

\def\reset#1{\global#1=0 }
\def\useouter #1 {\begingroup
                  \expandafter\let\csname#1\endcsname=\relax % not outer
                  \afterassignment\endgroup\global}
\def\?#1\then#2\else#3\fi{#1\defnext{#2}\else\defnext{#3}\fi \next}
\newif\ifnot
\newif\ifnomessage \nomessagetrue % used to suppress expansions in \message
\let\domessage=\message \def\message#1{{\nomessagefalse\domessage{#1}}}

% Fonts

% Define some fonts that are not present in plain.tex, and a general switch
% \sans that changes the main text fonts to sans-serif

\font\boldit=cmbxti10 % for subsection headings
\font\ssq=cmssq8 \font\ssqi=cmssqi8 % for quotations
\font\tenss=cmss10
\font\tenssit=cmssi10
\font\tenssbf=cmssbx10

\def\sans
{\def\rm{\fam0\tenss}\def\it{\fam\itfam\tenssit}\def\bf{\fam\bffam\tenssbf}\rm
}

% And a general mechanism to switch to nine point type
% This changes math as well, using sizes 9--7--5 pt

\font\ninerm=cmr9 \font\nineit=cmti9 \font\ninebf=cmbx9 \font\ninesl=cmsl9
\font\ninett=cmtt9 \font\ninei=cmmi9 \font\ninesy=cmsy9
\font\niness=cmss9 \font\ninessit=cmssi9 \font\ninessbf=cmssbx10 at 9 pt
\skewchar\ninei='177 \skewchar\ninesy='60

\def\ninept
{\def\rm{\fam0\ninerm}\def\it{\fam\itfam\nineit}\def\bf{\fam\bffam\ninebf}%
 \def\sl{\fam\slfam\ninesl}\def\tt{\fam\ttfam\ninett}\rm
 \def\sans{\def\rm{\fam0\niness}\def\it{\fam\itfam\ninessit}%
           \def\sl{\fam\slfam\ninessit}\def\bf{\fam\bffam\ninessbf}\rm}%
 \textfont0\ninerm \textfont1\ninei \textfont2=\ninesy
 \textfont\itfam\nineit \textfont\slfam\ninesl
 \textfont\bffam\ninebf \textfont\ttfam\ninett
 \normalbaselineskip=11pt \normalbaselines
 \setbox\strutbox=\hbox{\vrule height 8 pt depth 3 pt width 0 pt}%
}

% File handling

\newif\iftoc 
\def\writetoc{\write\tocfile}\def\tocfile{-1 }
\useouter newwrite
\def\opentoc{\newwrite\tocfile\openout\tocfile=\jobname.toc }
\edef\tocheader{\string\NotocSection\space Table of contents.
  \string\noindent\string\medskip}
\def\maketoc{\opentoc\writetoc{\tocheader}\toctrue}

\def\str#1{\string#1 }
% The following prohibits expansion between $..$,
% but there shouldn't be any groups there.
{\catcode`@=3 % weird math shifts, cannot occur in math mode
 \gdef\passord#1$#2\end{#1\ifx@#2@\else\getmath#2\end\fi}
 \gdef\getmath#1${$\noex#1\endmath$\passord}
 \gdef\noex#1\endmath{\ifx@#1@\else\donoex#1\endmath\fi}
 \gdef\donoex#1{\noexpand#1\noex}
 \gdef\nomathexp#1{\passord#1$\end}
}

\iffalse$\fi % to get emacs' font-lock adrift again

% Calling \putmark{Some title text} issues a \mark, in which the current
% values of \secno, \subsecno, and \proclno are recorded in the form of
% assignments to registers \count1-3; in that way when \firstmark is expanded
% during \output, these registers get their proper value, resulting in
% the display of meaningful values in the log after the page number. (Since
% the macros that call \putmark often reset inferior counters like \proclno to
% 0 before doing so, the log will usually contain less than three numbers
% after the page number.) At the same time the \mark records the relevant part
% of those counters in text form (from \lastlabel, assumed to be set before
% calling \putmark) together with the supplied title text in the token
% register \marktoks, where it is available for use in running heads. The
% macro \setmark combines this action with writing a line containing the
% expansion of \lastlabel and a restictively expanded (using \nomathexp)
% version of the title text to the table-of-contents file.

\newtoks\marktoks

\def\putmark#1{\mark
  {\count1=\the\secno \count2=\the\subsecno \count3=\the\proclno
   \marktoks={\ignorespaces\setchap\space#1}}}% set counters, running head
\def\setchap{\ifnum\secno<0 Appendix \lastlabel.\else\lastlabel\fi}
\def\puttoc#1{\toks0={\nomathexp{#1} \str\onpage\folio.}\edef\next
                      {{\str\tocitem \lastlabel=\the\toks0}}\writetoc\next}
\def\setmark#1{\putmark{#1}\puttoc{#1}}

\output=
{\everyhbox{}% don't interfere into output routine
 \firstmark \headline={\tenit \ifodd\pageno\hfill\fi \the\marktoks \hfil }
 \let~\relax % ~ shouldn't be expanded in index file, and doesn't occur in:
 \plainoutput
}

\newwrite\labfile
\def\tolab{\write\labfile}
\def\writelab{\openout\labfile=\jobname.lab
              \global\let\writelab=\tolab \tolab}
\def\chklab#1\lastlabel#2.{\def\next{#2}\ifx\next\empty}% if #1/=\lastlabel
\def\deflab#1{\chklab#1\lastlabel.% when #1 is some unique \cs, it will expand
                \immediate\writelab{\string\def\string#1{#1}}% OK any time, but
              \fi} % using \immediate is prudent (would allow #1{\lastlabel} )
\def\savecount#1{\writelab{#1=\number#1}}
\toks0={\savecount\secno\savecount\subsecno\savecount\proclno
        \savecount\eqnumber
        \savecount\pageno\writelab{\string\advancepageno}}
\toks2=\expandafter{\bye}
\edef\bye{\the\toks0 \the\toks2}

\newwrite\indexfile
\def\toindex{\write\indexfile}
\def\writeindex{\openout\indexfile=\jobname.inx
                \global\let\writeindex=\toindex \toindex}
\def\inx#1{\writeindex{\nomathexp{#1} @\folio.}}
\def\index#1{\inx{#1}#1}
\def\closeindex{\vfil\eject % Ensure final entries are written
                \ifx\writeindex\toindex \immediate\closeout\indexfile \fi}

% Layout style

% Redefine plain.tex's footnotes to be more independent of context settings,
% and in 9 pt type
\def\vfootnote#1%
{\insert\footins\bgroup \ninept
 \interlinepenalty=\interfootnotelinepenalty
 \splittopskip=\ht\strutbox % determines top baselines of broken footnotes
 \splitmaxdepth=\dp\strutbox % likewise for bottom lines
 \floatingpenalty=\@MM % please don't break footnotes, if possible
 \leftskip=1pc \rightskip=1pc \spaceskip=0pt \xspaceskip=\z@skip
 \parindent=\z@ \parskip=3pt minus 1pt \parfillskip=\z@ plus 1 fil
 \textindent{#1}\footstrut\futurelet\next\f@@t
}

\newcount\secno \newcount\subsecno \newcount\proclno \newcount\subsubsecno
\newcount\itemno
\reset\secno \reset\subsecno \reset\proclno \reset\itemno
\newdimen\itemindent \itemindent=\parindent % the default, but changeable

\global\subsubsecno=-1 % by default no subsubsections, set to 0 if using them

\newif\ifnewsec % communication with proclaim etc.

\let\Sec\S % Section mark
\useouter beginsection
 \outer\def\newsection #1.
  {\ifnewsec \endgroup \fi % previous (sub)section might be empty
   \global\advance\secno by 1
   \ifnum-1<\subsecno \reset\subsecno \fi % <0 means no subsections used
   \ifnum-1<\subsubsecno \reset\subsubsecno \fi % <0 means no subsubsections
   \reset\proclno
   \edef\lastlabel{\number\secno}%
   {\parskip0pt\beginsection\ifnomessage\Sec\fi\lastlabel. #1.\par
    \par}% immediately terminate the paragraph started by \beginsection
   \setmark{#1}% this mark follows the non-breakable space below the title
   \startsec % this replaces the indent-defeating code of \beginsection
  }
\def\subsection #1.
  {\secbigbreak
   \global\advance\subsecno by 1
   \ifnum-1<\subsubsecno \reset\subsubsecno \fi % <0 means no subsubsections
   \reset\proclno
   \edef\lastlabel{\number\secno.\number\subsecno}%
    \leftline{\bf\lastlabel.\enspace
              \defnext{#1}\ifx\next\empty\else \boldit#1.\fi}\setmark{#1}%
    \nobreak\smallskip\startsec}
\def\subsubsection #1.
  {\secmedbreak
   \global\advance\subsubsecno by 1 \reset\proclno
   \edef\lastlabel{\number\secno.\number\subsecno.\number\subsubsecno}%
    \leftline{\bf\lastlabel.\enspace
              \defnext{#1}\ifx\next\empty\else \it#1.\fi}\setmark{#1}%
    \nobreak\smallskip\startsec}
\useouter beginsection
 \outer\def\Appendix #1.
  {\ifnum\secno>0 \startapx\fi
   \global\advance\secno by -1
   \ifnum-1<\subsecno \reset\subsecno \fi % <0 means no subsections used
   \ifnum-1<\subsubsecno \reset\subsubsecno \fi % <0 means no subsubsections
   \reset\proclno
   \edef\lastlabel
    {\ifcase -\secno \or A\or B\or C\or D\or E\or F\or G\or H\or I\or J\fi}%
   {\parskip0pt\beginsection Appendix \lastlabel. #1.\par\par}%
   \setmark{#1}\startsec}
\def\startapx{\ifx\writetoc\totoc\totoc{\str\Appendices}\fi \secno=0 }
\useouter beginsection
 \outer\def\Supplement #1.
  {\let\lastlabel=\empty
   \global\ifnum\secno<0 \secno=1000 \else\advance\secno by 1 \fi
   \beginsection #1.\par\par\setmark{#1}\startsec}
\useouter beginsection
 \outer\def\NotocSection #1.
  {\let\lastlabel=\empty \reset\secno
   \beginsection #1.\par\par\putmark{#1}\startsec}

\def\sticksection % to be issued before \beginsection to prevent page break
{\vbox\bgroup\parskip=0pt \everypar{\egroup\noindent}}

\def\startsec
  {\begingroup \newsectrue\parskip0pt\parindent0pt \everypar{\endgroup}}
\def\secbigbreak
  {\ifnewsec \endgroup
     \ifdim\lastskip<\bigskipamount \removelastskip\bigskip \fi % no penalty
   \else \bigbreak
   \fi}
\def\secmedbreak
  {\ifnewsec \endgroup
     \ifdim\lastskip<\medskipamount \removelastskip\medskip \fi % no penalty
   \else \medbreak
   \fi}
\def\labelsec#1{\global\let#1=\lastlabel \deflab#1\nobreak}% prolong nobreak

% for \startsec to work properly, a paragraph should not begin with `{text}'
% here we fix those macros that do such a thing. In the meantime we deactivate
% the accents within messages, so that only the base character will appear

{\let\x=\expandafter
 % in the following macro #1 is the \cs to be redefined; \next is defined to
 % do the redefinition, which it achieves using the old expansion passed to it
 % as an argument; after defining \next, it is called with the expansion of #1
 \def~#1{\defnext##1{\gdef#1{\ifnomessage\leavevmode##1\fi}}\x\next\x{#1}}
 ~\copyright
 % the following is similar, but now \cs itself takes one argument. When
 % passing the 'expansion' of \cs, that argument is given as '#1', which will
 % capture the formal argument of \cs in the redefinition within \next
 \def~#1{\defnext##1{\gdef#1####1{\ifnomessage\leavevmode##1\else####1\fi}}%
         \x\next\x{#1{##1}}}
 ~\` ~\' ~\v ~\^^_ ~\u ~\^^S ~\^ ~\^^D ~\. ~\H ~\~ ~\" 
 \ifamstex ~\B \else ~\= ~\t \fi
}

\def\quotation
{\begingroup\everypar{}\parfillskip=0pt
 \leftskip=0pt plus \hsize \def\par{\ifhmode\/\endgraf\nobreak\fi}\obeylines
 \let\it=\ssq\let\rm=\ssqi\rm % switch roman and italic
 \baselineskip=10pt\parskip=0pt
 \def\author##1, ##2.{\vskip2pt\it ##1, \rm ##2\par\smallskip\endgroup}
}

\newcount\eqnumber \eqnumber=0
\def\neqn(#1#2){\global\advance\eqnumber by 1
                (\ifcat a\noexpand#1%
                        \errmessage{Incorrect label #1#2}\number\eqnumber
                 \else\xdef#1{\number\eqnumber}\deflab#1%
                      #1\fam\z@#2% #2 is possible tail, like "a,b"
                 \fi)}
\def\label{\ifinner\else \eqno \fi \neqn}
\def\nn{\label(\lastlabel)}
\newcount\acc
\def\add#1#2{\acc=#1\advance\acc#2\relax\number\acc}

\def\heading "#1"{\smallbreak\noindent{\sl #1}.\quad\ignorespaces}

\let\proclaim=\relax
\ifamstex
% kludge needed because \endtopmatter sets \proclaim to \innerproclaim@ (again)
\let\amsproclaim=\innerproclaim@
\def\innerproclaim@#1. #2#3\par
  {\secmedbreak \global\advance\proclno by 1 \reset\itemno
   \edef\lastlabel{\ifnum-1<\secno    \number\secno.\fi
                   \ifnum-1<\subsecno \number\subsecno.\fi
                   \ifnum-1<\subsubsecno \number\subsubsecno.\fi
                   \number\proclno}
   {\interlinepenalty=300 \parskip=0pt\def\par{\endgraf\penalty200 }% (item)
    \noindent % Emit \baselineskip & \parskip glue BEFORE doing \deflab
    \nottrue % here this means: no label to define
    \if\let\noexpand#2\ifx\undefined#2\notfalse\fi\fi% Is #2 undefined \cs?
    \ifnot\amsproclaim{\lastlabel.\enspace#1}#2#3\innerendproclaim@
    \else\global\let#2\lastlabel \deflab#2%
      \amsproclaim{\lastlabel.\enspace#1}#3\innerendproclaim@
    \fi
   }}
\else
\outer\def\proclaim#1. #2#3\par
  {\secmedbreak \global\advance\proclno by 1 \reset\itemno
   \edef\lastlabel{\ifnum-1<\secno    \number\secno.\fi
                   \ifnum-1<\subsecno \number\subsecno.\fi
                   \ifnum-1<\subsubsecno \number\subsubsecno.\fi
                   \number\proclno}
   {\interlinepenalty=300 \parskip=0pt\def\par{\endgraf\penalty200 }% (item)
    \noindent % Emit \baselineskip & \parskip glue BEFORE doing \deflab
    \if\let\noexpand#2% Is #2 a control sequence?
     \ifx\undefined#2% don't redefine a control sequence
      \global\let#2\lastlabel \deflab#2\letnext\ignorespaces
     \else\defnext{#2}\fi
    \else\defnext{#2}\fi
    \bf\lastlabel.\enspace#1.\enspace \sl\next#3\endgraf
    }%
   \ifdim \lastskip<\medskipamount \removelastskip\penalty55\medskip \fi
   }
\fi

\def\dueto#1.{\unskip~[#1].}

\def\newpar{\ifhmode\par\fi{\parskip0pt\noindent}}% no extra vertical space

\def\item{\doindent\itemindent}

\def\doindent#1#2{\newpar
  \hangindent#1\hangafter0 \llap{#2\enspace}\ignorespaces}

\def\statitemnr#1{(\number#1)}
\def\statitem{\global\advance\itemno by 1 \item{\statitemnr\itemno}}
\def\eqitemnr#1{(\romannumeral#1)}
\def\eqitem{\global\advance\itemno by 1 \item{\eqitemnr\itemno}}
\def\defitemnr#1{\count@=#1\advance\count@96\char\count@}
\def\defitem{\global\advance\itemno by 1 \item{\defitemnr\itemno.}}

\def\proof{\heading "Proof"} 
\def\QED{\dimen0=.5em\advance\dimen0\wd\QEDbox\advance\dimen0.4pt
         \unskip\kern\dimen0\vadjust{\vskip-\baselineskip}\unhcopy\tightbox
         \copy\QEDbox{\pretolerance=1000\parfillskip=0.4pt\smallbreak}}
\def\QEDQED{\unskip\nobreak\enskip\kern\wd\QEDbox
            {\setbox\QEDbox=\hbox{\copy\QEDbox\enspace\copy\QEDbox}\QED}}

\def\emph#1{{\it #1}\futurelet\next\doitcor}
\def\newterm#1{{\sl\index{#1}}\futurelet\next\doitcor}
\def\foreign#1{{\sl #1}\futurelet\next\doitcor}
\def\doitcor{\if.\next\else\if,\next\else\/\fi\fi}
\def\needfil{\penalty-100\vfil\penalty1500\vfilneg}% raggedbottom, if needed
\def\bigdisplay{\ifmmode \errmessage{Display begun in math mode.}
                \else\ifhmode\unskip\vadjust{\needfil}%
                     \else\needfil\noindent\fi$$\fi}
\def\frame#1{\leavevmode\hbox{\vrule\vtop{\vbox{\hrule\kern2pt
     \hbox{\strut\kern2pt#1\kern2pt}}\kern2pt\hrule}\vrule}} % framed box
\def\clap#1{\hbox to 0pt {\hss#1\hss}}% centered lap
\long\def\centerpar#1% centered paragraph(s)
{\begingroup
  \advance\leftskip by 0pt plus 1 fil minus 3pt \rightskip=\leftskip
	\parskip=0pt \parindent=0pt \parfillskip=0pt \linepenalty=100
  \noindent\ignorespaces#1\par
 \endgroup
}

% Changing the \baselineskip is a bit subtle, since we want to make sure it
% only kicks in _after_ the current paragraph. We place a dummy line after the
% paragraph at the old \baselineskip, and than back up at the new one.

\def\setbaselineskip{\par\afterassignment\setbskip\skip0 =}
\def\setbskip{\null\baselineskip=\skip0 \vskip-\baselineskip}
\def\forceeven{\eject\ifodd\pageno \null\mark{\marktoks{}}\vfil\eject \fi}
\def\forceodd{\eject\ifodd\pageno\else\null\mark{\marktoks{}}\vfil\eject\fi}
\def\topglue{\eject\afterassignment\topglu\skip0 }
\def\topglu{\hrule height0pt\vskip\skip0\vskip-\topskip
     \dimen0=\baselineskip\advance\dimen0-\topskip \prevdepth=\dimen0 }

% output routine resets topskip glue

\def\qdef#1{\?\ifx #1\undefined \then \def #1%
              \else \toks0={#1}\afterassignment\chkdef\def\next
              \fi}
\def\qqdef#1{%
\?\ifx #1\undefined
  \then \message{Patching definition:\string#1.}\let\Qdef=\qdef\def#1%
  \else \toks0={#1}\afterassignment\chkdef\def\next
  \fi}
\let\Qdef=\qqdef
\def\chkdef{\expandafter\let\expandafter\curcs\the\toks0 % space is crucial!
  \ifx\next\curcs\else
  \message{Altered definition:\the\toks0= \meaning\next,
                              should have been \meaning\curcs.}\fi}
\ifamstex

\def\mathtex#{\quad\hbox}
\def\mathttex#{\qquad\hbox}
\else
\def\text#1{\quad\hbox{#1}\quad} \def\ttext#1{\qquad\hbox{#1}\qquad}
\def\tex#{\quad\hbox}            \def\ttex#{\qquad\hbox}
\fi

% Symbols

\newcount\t@
\def\mathdef #1=#2 #3% #1 \cs,  #2: class, #3 model char (#3VarFam: #4 family)
{\count@=\z@
 {\def\t{#2}\nottrue
  \loop \xdef\next{\ifcase \count@
    Ord\or Op\or Bin\or Rel\or Open\or Close\or Punct\or VarFam\else #2\fi}%
    \ifx \next\t \notfalse \fi
  \ifnot \global\advance\count@ by 1 \repeat
 }\multiply\count@"1000
 \t@=\expandafter
 \ifx\expandafter\delimiter#3\divide\t@"1000 \nottrue
     % e.g., \mathdef\orbit=Bin \backslash sets \t@="026E since
     % \backslash -> \delimiter"026E30F 
 \else
   \ifcat \alpha#3#3\nottrue % a \cs (presumably \mathchardef'ed); use it
   \else \mathcode`#3 \ifnum\t@<"7000 \nottrue \else \notfalse \fi
   \fi
 \fi
 \advance\count@\t@
 \ifnot \divide\t@"1000 \multiply\t@"1000 \advance\count@-\t@
   \mathchardef#1=\count@ % use family and member from #3, but replace class
 \else \divide\t@"100 \multiply\t@"100 \advance\count@-\t@
   \def\next{#1}\afterassignment\mathd@f \t@=
   % for class VarFam ignore its family, using following number (#4) instead
 \fi
}
\def\mathd@f
{\multiply\t@"100 \advance\count@\t@ \expandafter\mathchardef\next=\count@}
\def\mathorddef #1={\mathdef #1=Ord }

\ifamstex  \let\romath=\roman
\else  \def\romath#1{{\fam\z@#1}}
\fi

\mathorddef\N=N\bffam
\mathorddef\Z=Z\bffam
\mathorddef\Q=Q\bffam
\mathorddef\R=R\bffam
\mathorddef\C=C\bffam
\mathorddef\P=P\bffam

\mathorddef\1=1\bffam
\mathorddef\2=2\bffam
\mathorddef\ii=i\bffam % imaginary unit
\mathorddef\\=\lambda
\mathdef\:=Punct :
\let\Id=\1
\mathorddef\Sym=S\bffam % symmetric group

\def\opdef#1{\def#1{\mathop\romath{\escapechar\m@ne\string#1}\nolimits}}
\def\altopdef#1#2{\def#1{\mathop\romath{#2}\nolimits}}
\def\li(#1#2..#3){%
  #1_{#2},\penalty\binoppenalty\ldots,\penalty\binoppenalty#1_{#3}}
\def\lix(#1#2.#3.#4){%
  #1_{#2},\penalty\binoppenalty\ldots,\widehat{#1_{#3}},\penalty\binoppenalty
  \ldots,#1_{#4}}
\def\fold#1(#2#3..#4){#2_{#3}#1\cdots#1#2_{#4}}
\def\foldx#1(#2#3.#4.#5){#2_{#3}#1\cdots#1\widehat{#2_{#4}}#1\cdots#1#2_{#5}}
\def\foldwith#1#2#3(#4..#5)% e.g., \foldwith{\mu^{(#1)}}\subset\cdots(0..l-1)
{\defnext##1{#1}\next{{#4}}#2#3#2\next{{#5}}}
\def\({\bigl(}\def\){\bigr)}
\def\inv{^{-1}}
\def\<#1>{{\langle #1\rangle}}

\mathdef\orbit=Bin \backslash
 % used to be \mathdef\union=Rel \cup
                      \let\Union=\bigcup
\let\thru=\cap        
     
\let\tensor=\otimes

 % function composition
\let\eps=\varepsilon
\def\setof#1:#2\endset{{\{\,#1\mid#2\,\}}}
\def\Setof#1:#2\endset{\left\{\,#1\,\,\vrule\,\,#2\,\right\}}
   % correct wrong names

\let\ssubset=\subset \let\subset=\subseteq
 
\let\baraccent=\= \def\={\relax\?\ifmmode\then\approx\else\baraccent\fi}

\opdef\sg \opdef\Ker\opdef\Im \opdef\Aut \opdef\End \opdef\rk
\altopdef\kar{char}
\def\Card#1{\##1}

\def\multiset#1{{\{\!\!\{#1\}\!\!\}}}
\def\multisetof#1:#2\endset{\multiset{\,#1\mid#2\,}}
\def\Multisetof#1:#2\endset{\left\{\!\!\left\{\,#1\,\,\vrule
                           \,\,#2\,\right\}\!\!\right\}}
\def\multichoose(#1,#2){\mathchoice
  {\left(\!{#1\choose#2}\!\right)}%
  {\bigl(\!{#1\choose#2}\!\bigr)}%
  {\left(\!{#1\choose#2}\!\right)}%
  {\left(\!{#1\choose#2}\!\right)}%
}

% Standard abbrev.'s

%\def\mo.{morphism}
%\def\iso.{iso\mo.}
%\def\cor.{correspond}
%\def\de.{definition}
%\def\De.{Definition}
%\def\lcit{loc.\thinspace cit.}

% Boxes

\def\topsmash#1{{\def\finsm@sh{\ht\z@\z@ \box\z@}\smash{#1}}}
\def\botsmash#1{{\def\finsm@sh{\dp\z@\z@ \box\z@}\smash{#1}}}

\newbox\tightbox % for "tight" breaks
\setbox\tightbox=
\hbox{\penalty -200               \skip0 = 0pt plus 0.05 \hsize
      \hskip \skip0 \penalty    0 
      \hskip \skip0 \penalty  200 \skip0 = .1\hsize
      \hskip \skip0 \penalty  283 
      \hskip \skip0 \penalty  346
      \hskip \skip0 \penalty  400
      \hskip \skip0 \penalty  447 
      \hskip \skip0 \penalty  490
      \hskip \skip0 \penalty  529
      \hskip \skip0 \penalty  567 
      \hskip \skip0 \penalty  600
      \hskip \skip0 \penalty  632
      \hskip 2\hsize % make it an offer TeX can't refuse
      \vadjust{}\hfil} % warning if you change \hsize do it FIRST. 
\newbox\QEDbox
\setbox\QEDbox=\hbox
{\kern-0.2pt\vrule height 6pt \kern-0.4pt % backspace
 \vbox to 6pt
  {\kern -0.2pt
   \hbox to 6.4pt{\hrulefill\kern-0.2pt\vrule}\vfil
   \hrule width 6.2pt
   \kern -0.2pt}%
 \kern-0.4pt \vbox to 6pt{\leaders\vrule\vfil\kern-0.2pt}\kern-0.2pt
}

% Nine point type (for footnotes and references); math sizes 9--7--5 pt

\ifx\ninept\undefined % don't overwrite more sophisticated \ninept macros

\font\ninerm=cmr9 \font\nineit=cmti9 \font\ninebf=cmbx9 \font\ninesl=cmsl9
\font\ninett=cmtt9 \font\ninei=cmmi9 \font\ninesy=cmsy9
\skewchar\ninei='177 \skewchar\ninesy='60

\def\ninept
{\def\rm{\fam0\ninerm}\def\it{\fam\itfam\nineit}\def\bf{\fam\bffam\ninebf}%
 \def\sl{\fam\slfam\ninesl}\def\tt{\fam\ttfam\ninett}\rm
 \textfont0\ninerm \textfont1\ninei \textfont2=\ninesy
 \textfont\itfam\nineit \textfont\slfam\ninesl
 \textfont\bffam\ninebf \textfont\ttfam\ninett
 \normalbaselineskip=11pt \normalbaselines
 \setbox\strutbox=\hbox{\vrule height 8 pt depth 3 pt width 0 pt}%
}

\fi

\def\caps#1{{\ninept #1}}

\ifamstex % then make sure the quotes required by AmSTeX are provided
\let\amsfootnote=\footnote
\def\footnote#1{\amsfootnote"#1"}
\fi

% References

\catcode`_=11 % control sequence names with this letter are private

\newif\ifamsrefs % whether AMS style refernece syntax is used

\ifamstex
\let\amsref=\ref
\let\amsendref=\endref
\let\amspubl=\publ
\let\amsinbook=\inbook
\let\amskey=\key
\let\amsed=\ed 
\let\amsvol=\vol
\let\amspages=\pages
\fi

\def\label_#1#2{\expandafter\xdef \csname#1.\endcsname{#2}}
\def\ref#1{\ref_#1,,,.}
\def\ref_#1,#2,,#3.{\expandafter\ref_lab\csname#1.\endcsname{#3#2}}
\def\ref_lab#1#2{\ifx #1\relax \errmessage{Reference to undefined label}%
                 \else [#1#2]\fi
                }

\newdimen\refindent
\refindent=36pt

\def\init_refs % For initialising reference labels at start of article
{\begingroup\setbox0=\vbox\bgroup \hbadness 10000 \hfuzz=\maxdimen
 \def\chapter ##1\par{}% avoid problems with chapter number and tocfile
 \def\endrefs{\egroup\egroup\endgroup }% end \beginrefs & \init_refs
 \ifamsrefs
 \else
  \def\endref{\iffalse{\fi}\ifx\the_key\undefined\else
     \label_{\expandafter\strip_space\the_key _}%
	    {\expandafter\strip_space\the_mark _}\fi
     \endgroup}%
 \fi
}
\def\strip_space #1 _{#1}

\def\readrefs{\init_refs \input \jobname.ref }

\useouter amsref \def\ams_start\key #1 \ident#2%
{\label_{#2}{#1}\amsref\amskey{#1}%
}
\outer\def\beginrefs{\noindent\par % fire up \everypar
 \ifamsrefs
  \bgroup \let\ref=\ams_start
  \let\publ=\amspubl \let\inbook=\amsinbook \let\key=\amskey
  \let\ed=\amsed \let\vol=\amsvol \let\pages=\amspages
 \else
  \bgroup \ninept
  \frenchspacing\tolerance=900
  \let\ref=\_ref \let\author=\_auth \let\title=\_tit
  \let\vol=\_vol \let\year=\_year
 \fi
}

\outer\def\endrefs{\egroup\smallbreak}

\def\ref_item#1{\iffalse{\fi}% must have an unbalanced explicit brace here
  \ifx#1\undefined \else
    \errmessage{Multiple definition of \string#1 within reference}\fi
  \edef#1{\iffalse}\fi % must be \edef to activate keywords (like \title)
}

\def\_ref{\begingroup
  \def\par{\endref\par}\def\ref{\endref\ref}\edef\the_mark{\iffalse}\fi}
\def\_auth{\ref_item\the_author}
\def\_tit{\ref_item\the_title}
\def\key{\ref_item\the_key}

\def\publ{\ref_item\the_publ}

\def\journal{\ref_item\the_journal}
\def\_vol{\ref_item\the_vol}
\def\_year{\ref_item\the_year}
\def\pages{\ref_item\the_pages\range_}

\def\inbook{\ref_item\the_inbook}
\def\ed{\ref_item\the_ed}
\def\report{\ref_item\the_report}

\def\range_{\noexpand\range_}   % these cancel the \edef in \ref_item
\def\amount_{\noexpand\amount_}
\def\do_amount{{\aftergroup\do_pp\aftergroup}}
\def\do_pp{~pp}
\def\expl_pp{\def\range_{pp.~}\let\amount_=\do_amount} % if `pp.' is desired
\def\impl_pp{\let\range_\empty\let\amount_=\do_amount} % if `pp.' not desired

\def\q_plur#1{\expandafter\q_s#1,_}\def\q_s#1,#2_{\ifx _#2_\else s\fi}

\def\_pages{\expandafter\pages_\the_pages\unskip---_}
\def\pages_#1-#2#3-#4_{#1\ifx _#4_\else\if -#2--#3\else --#2#3\fi\fi}

\ifamstex

\def\_opt#1#2{\ifx#1\undefined\else#2\fi}
\def\opt_pages{\_opt\the_pages{\amspages\expandafter\pages_\the_pages---_}}

\def\endref{\ifamsrefs \amsendref \else \ams_convert \fi}

\useouter amsref \def\ams_convert
{\iffalse{\fi}\let\par\endgraf
 \amsref \amskey\the_mark \by\the_author
 \ifx\the_inbook\undefined
   \ifx\the_journal\undefined
     \ifx\the_publ\undefined
       \ifx\the_report\undefined 
	 \errhelp{A reference must at least have a \journal \inbook \publ^^J
		  or \report field defined.}%
	 \errmessage{Unrecognised reference}%
       \else \make_report \fi
     \else \make_book \fi
   \else \make_article \fi
 \else\make_proceedings \fi
 \ifx\the_endnote\undefined\else \finalinfo\the_endnote \fi
 \amsendref
 \endgraf\endgroup
}

\def\make_report
{\paper\the_title \paperinfo\the_report\unskip \opt_pages
 \_opt\the_year{\yr\the_year}%
}
\def\make_book
{\book\the_title \opt_pages
 \_opt\the_series
 {\bookinfo\the_series
  \_opt\the_ed{\csname amsed\q_plur\the_ed\endcsname\the_ed}%
  \_opt\the_vol{\amsvol\the_vol\_opt\the_idno{\unskip, \the_idno}}%
 }%
 \amspubl\the_publ \yr\the_year \_opt\the_ISBN{\bookinfo\the_ISBN}%
}
\def\make_article
{\paper\the_title \jour\the_journal
 \_opt\the_vol{\amsvol\the_vol\_opt\the_idno{\unskip, \the_idno}}%
 \_opt\the_year{\yr\the_year}\opt_pages
}
\def\make_proceedings
{\paper\the_title \opt_pages \amsinbook\the_inbook
 \_opt\the_ed{\csname amsed\q_plur\the_ed\endcsname\the_ed}%
 \_opt\the_series{\bookinfo\the_series
   \_opt\the_vol{\amsvol\the_vol\_opt\the_idno{\unskip, \the_idno}}}%
 \_opt\the_publ{\amspubl\the_publ}\yr\the_year
}

\else % \amstexfalse

\def\_opt#1#2{\ifx#1\undefined\else{\def\__{#1\unskip}#2}\fi}
\def\opt_pages#1{\ifx\the_pages\undefined\else{\let\__\_pages#1}\fi}

\def\endref
{\iffalse{\fi}\let\par\endgraf\smallbreak\hangindent\refindent
 \setbox0=\hbox{[\the_mark\unskip]\enspace}\noindent
 \ifdim\wd0<\refindent \hbox to \refindent{\unhbox0\hfil}\else\unhbox0 \fi
 \format_author{\the_author\unskip},
 \ifx\the_inbook\undefined
   \ifx\the_journal\undefined
     \ifx\the_publ\undefined
       \ifx\the_report\undefined 
	 \errhelp{A reference must at least have a \journal \inbook \publ^^J
        	  or \report field defined.}%
         \errmessage{Unrecognised reference}%
       \else \format_report \fi
     \else \format_book \fi
   \else \format_article \fi
 \else\format_proceedings \fi
 \unskip \_opt\the_endnote{. \__}.
 \endgraf\endgroup
}

\def\format_author#1{{\rm#1}}

\def\format_report
{{\sl \the_title\unskip}, \the_report\unskip\expl_pp\opt_pages{, \__}%
 \_opt\the_year{, (\__)}}
\def\format_book
{{\sl\the_title\unskip},\expl_pp\opt_pages{ \__,}%
 \_opt\the_series{ \__\_opt\the_ed{ (\__, ed\q_plur\the_ed.)}%
 \_opt\the_vol{ \__}\_opt\the_idno{, \__},} \the_publ\unskip,
 \the_year\unskip\_opt\the_ISBN{, ISBN \__}}
\def\format_article
{``\the_title\unskip'', {\sl\the_journal\unskip\_opt\the_vol{\/ \bf\__}}%
 \_opt\the_idno{, \__}\_opt\the_year{, (\__)}\impl_pp\opt_pages{, \__}}
\def\format_proceedings
{``\the_title\unskip'',\expl_pp\opt_pages{ \__} in {\sl\the_inbook\unskip}%
 \_opt\the_ed{ (\__, ed\q_plur\the_ed.)},
 \_opt\the_series{\__\_opt\the_vol{ \__}\_opt\the_idno{, \__}, }%
 \_opt\the_publ{\__, }\the_year}

\fi % \ifamstex

\catcode`_=8 % Now underscores are for subscripts again
\catcode`@=12 % and `@' cannot be used in control sequences

%auto-ignore (forbid arXiv's auto-postscript generator TeXing this file alone)
\ifx\Tableauxloaded\relax   \else \let\Tableauxloaded\relax \fi

% Het tekenen van Young tableaux

\newcount\cols
{\catcode`,=\active\catcode`|=\active
 \gdef\Young(#1){\hbox{$\vcenter
 {\mathcode`,="8000\mathcode`|="8000
  \def,{\global\advance\cols by 1 &}%
  \def|{\cr
        \multispan{\the\cols}\hrulefill\cr
        &\global\cols=2 }%
  \offinterlineskip\everycr{}\tabskip=0pt
  \dimen0=\ht\strutbox \advance\dimen0 by \dp\strutbox
  \halign
   {\vrule height \ht\strutbox depth \dp\strutbox##
    &&\hbox to \dimen0{\hss$##$\hss}\vrule\cr
    \noalign{\hrule}&\global\cols=2 #1\crcr
    \multispan{\the\cols}\hrulefill\cr%
   }
 }$}}
 \gdef\Skew(#1:#2){\hbox{$\vcenter
 {\mathcode`,="8000\mathcode`|="8000
  \dimen0=\ht\strutbox \advance\dimen0 by \dp\strutbox
  \def\boxbeg{\vbox
    \bgroup\hrule\kern-0.4pt\hbox to\dimen0\bgroup\strut\vrule\hss$}%
  \def\boxend{$\hss\egroup\hrule\egroup}%
  \def,{\boxend\boxbeg}%
  \def|##1:{\boxend\vrule\egroup\nointerlineskip\kern-0.4pt
    \moveright##1\dimen0\hbox\bgroup\boxbeg}%
  \def\\##1\\##2:{\boxend\vrule\egroup\nointerlineskip\kern-0.4pt
    \kern ##1\dimen0\moveright##2\dimen0\hbox\bgroup\boxbeg}%
  \moveright#1\dimen0\hbox\bgroup\boxbeg#2\boxend\vrule\egroup
 }$}}
}

% We need script and scriptscript sizes of text italic in small squares

\font\sevenit=cmti7 \font\fiveit=cmti7 at 5 pt
\scriptfont\itfam=\sevenit \scriptscriptfont\itfam=\fiveit

% The following downgrades members of roman (0), math italic (1), symbol (2),
% text italic and boldface families by one step, and adapts \strutbox
% correspondingly. This will result in smaller boxes and smaller formulae.
% These parameters are used when wrapping up a math formula, so they cannot be
% made local to a subformula (unless the latter is manually \hbox-ed).

\def\smallsquares
{\textfont0=\scriptfont0 \scriptfont0=\scriptscriptfont0
 \textfont1=\scriptfont1 \scriptfont1=\scriptscriptfont1
 \textfont2=\scriptfont2 \scriptfont2=\scriptscriptfont2
 \textfont\itfam=\scriptfont\itfam \scriptfont\itfam=\scriptscriptfont\itfam
 \textfont\bffam=\scriptfont\bffam \scriptfont\bffam=\scriptscriptfont\bffam
 \setbox0=\hbox{\raise0.4pt\hbox{$($}}
 \setbox\strutbox=\hbox{\vrule width 0pt height\ht0 depth\dp0 }% smaller
}

\def\bigsquares
{\setbox\strutbox=\hbox
 {\vrule width 0pt height1.3\ht\strutbox depth1.3\dp\strutbox}% larger
}

\newbox\hpair \newbox\vpair
\def\installboxes
{\dimen0=\ht\strutbox \advance\dimen0 by \dp\strutbox
 \setbox0=\vbox{\hrule width \dimen0}%
 \setbox1=\rlap{\strut\vrule}%
 \setbox\hpair=\rlap{\raise\ht\strutbox\copy0 \lower\dp\strutbox\copy0}%
 \setbox\vpair=\rlap{\raise\dimen0\copy1 \kern\dimen0 \copy1}%
}

\newcount\ribbonlength
\newcount\joff \newcount\jmax \newcount\imid \newcount\jmid

% Note that in the user-level commands \ribbon, \edgeseq and \labeledges the
% first coordinate counts downwards, but in the internal calculations of
% \ribs, \edges, \edgelabs, the vertical counter \count1 counts the upwards
% displacement from the point of departure (since this is the way of ribbons).

\def\ribbon #1,#2:#3;#4 % (#1,#2)->bottom-left; #3=content; #4=shape
{\setbox2=\hbox
 {\count0=0 \count1=0
  \copy1 \rlap{\lower\dp\strutbox\copy0}\if|#4|\else\ribs#4 \fi
  \raise\count1\dimen0\rlap
    {\kern\count0\dimen0 \raise\ht\strutbox\copy0 \copy1\raise0.4pt\copy1}%
  \global\joff=\count0
  \divide\dimen0 by 2 \raise\imid\dimen0\rlap{\kern\jmid\dimen0
    \hbox to 2\dimen0{\hfil$#3$\hfil}}%
 }%
 \lower#1\dimen0\rlap{\kern#2\dimen0\box2}%
 \advance\joff by #2\relax \ifnum \joff>\jmax \global\jmax=\joff \fi
}
\def\ribs#1#2
{\raise\count1\dimen0\rlap
   {\kern\count0\dimen0 \copy\ifnum #1=0 \hpair \else \vpair \fi}%
 \advance\ribbonlength by -2
 % Now increment \count0 or \count1, depending on #1. Also, if
 % \ribbonlength=1, set (\jmid,\imid):=2*(\count0,\count1) after
 % incrementing, and if \ribbonlength=0 set (\jmid,\imid) to the sum of
 % (\count0,\count1) before and after incrementing. Since `2\count0' is not
 % a valid <integer>, the easiest way is to set (\jmid,\imid) explicitly both
 % for \ribbonlength=0 and \ribbonlength=1, and to then advance the
 % coordinates if \ribbonlength=0 OR \ribbonlength=1
 \ifnum \ribbonlength=0 \jmid=\count0 \imid=\count1 \fi
 \advance\count#1 by 1
 \ifnum \ribbonlength=1 \jmid=\count0 \imid=\count1 \fi
 \ifnum \ribbonlength>-1 \ifnum \ribbonlength<2
   \advance \jmid by \count0 \advance \imid by \count1
 \fi\fi
 \if |#2|\else \ribs#2 \fi % loop while there are bits in the specified string
}

\def\rtab#1
{\installboxes \ribbonlength=#1 \global\jmax=0
 \setbox0=\hbox{\aftergroup\endrtab \aftergroup}}
\def\endrtab
{\advance\dimen0 by \jmax\dimen0 % width is maximal horizontal coordinate + 1
 \wd0=\dimen0 \vcenter{\box0}}

\def\edgeseq#1,#2;#3 % (#1,#2)->bottom-left grid point; #3=shape
{\setbox2=\hbox{\count0=0 \count1=0 \edges#3
  \global\joff=\count0
  }%
 \rlap{\kern#2\dimen0
       \dimen1=#1\dimen0 \advance\dimen1 by -\dimen0 \lower\dimen1\box2
      }%
 \advance\joff by #2\advance\joff by -1 % compensate for +1 above
 \ifnum \joff>\jmax \global\jmax=\joff \fi
}
\def\edges#1#2
{\rlap
 {\kern\count0\dimen0
  \ifnum #1=0
     \dimen1=\count1\dimen0 \advance\dimen1 by -\dp\strutbox
     \raise\dimen1\copy0
  \else \raise\count1\dimen0\copy1
  \fi}%
 % Now increment \count0 or \count1, depending on #1.
 \advance\count#1 by 1
 \if |#2|\else \edges#2 \fi % loop while there are bits in the specified string
}

% In \labeledges, argument #3 must be a string with the same number of items
% as the bit string #4. Use braces around non-atomic items.

\def\labeledges#1,#2:#3;#4 % (#1,#2)->bottom-left; #3=content string; #4=shape
{\setbox2=\hbox{\count0=0 \count1=0 \divide\dimen0 by 2 \edgelabs#3;#4 }%
 \rlap{\kern#2\dimen0
       \dimen1=-\ht\strutbox % move up by \ht\strutbox
       \advance\dimen1 by \fontdimen22\textfont2 % but lower by axis height
       \advance\dimen1 by #1\dimen0 \lower\dimen1\box2}%
}
\def\edgelabs#1#2;#3#4
{\advance\count#3 by 1 % first increment of \count0 or \count1
 \rlap {\kern\count0\dimen0 \raise\count1\dimen0 \clap{$#1$}}%
 \advance\count#3 by 1 % second increment of \count0 or \count1
 \if |#4|\else \edgelabs#2;#4 \fi % loop while there are bits
}

\def\arrow#1,#2:#3
{\smash{\lower#1\dimen0\rlap{$\kern#2.5\dimen0
 \if r#3 \,\rightarrow \else
 \if l#3 \llap{$\leftarrow\,$} \else
 \if u#3 \setbox0=\hbox to 0pt{\hss$\uparrow$\hss}\ht0=0pt 
         \raise.5\dimen0 \box0 \else
 \if d#3 \setbox0=\hbox to 0pt{\hss$\downarrow$\hss}\dp0=-2pt
         \lower.5\dimen0 \box0 \else
 \if i#3 \setbox0=\hbox to 0pt{\hss$\nwarrow\,$}\ht0=0pt
         \raise.5\dimen0 \box0 \else
 \if o#3 \setbox0=\hbox to 0pt{$\,\searrow$\hss}\dp0=-2pt
         \lower.5\dimen0 \box0 \else
 \fi\fi\fi\fi\fi\fi
 $}}%
}

%auto-ignore (forbid arXiv's auto-postscript generator TeXing this file alone)

\def\ifabsent#1\intoks#2%
{\defnext##1#1##2##3\end {\ifx\ifabsent##2}%
 \expandafter\next\the#2#1\ifabsent\end
}
% #2 is token list; if token #1 absent in \the#2, then \iftrue

\let\intoks=\if % act as \if w.r.t. skipping conditionals

\def\ifremoved#1\intoks#2{%
    \defnext##1#1##2#1##3\end{\ifx#1##3\global#2{##1##2}}%
    \expandafter\next\the#2#1#1\end}
% if token #1 present in \the#2, then then remove it, yield true, else false

\def\addtotoks#1#2{\ifabsent#2\intoks#1%
    \edef\next{\global#1{\the#1\noexpand#2}}\next\fi}
% add token to token list (at end) unless already present

\newwrite\frwfile
\def\tofrw{\immediate\write\frwfile}
\def\writefrw{\immediate\openout\frwfile=\jobname.frw
              \global\let\writefrw=\tofrw \tofrw}

\newtoks\oldforward % previous forward references, not yet defined anew
\newtoks\newforward % labels used in this run before they are defined

\begingroup
  \openin15=\jobname.frw
  \ifeof15 \closein15
  \else
   \closein15
   \def\forward#1=#2{\gdef#1{#2}\addtotoks\oldforward#1}
   \input\jobname.frw
  \fi
\endgroup

\newif\ifforward
\def\deadref{{\bf?}}

% forward references are written \forward\csname (or \forward{\csname})
% If \csname is undefined, we print \deadref, and issue a warning.
% If \csname was defined in \jobname.frw, we include it in \newforward.
% If neither is the case, it is a backward reference, and we do nothing.

\def\forward#1% #1 is a control sequence
{\ifx#1\undefined \global\let#1\deadref \forwardtrue
   \immediate\write16{Warning: unresolved forward reference to \string#1.}%
 \else \ifabsent#1\intoks\oldforward \forwardfalse \else\forwardtrue \fi
 \fi
 \ifforward \addtotoks\newforward#1\fi
 #1%
}

\def\definefrw#1% define previously referenced symbol to \lastlabel
{\ifremoved#1\intoks\newforward % label was forward referenced
  \ifx#1\lastlabel\else % if it was neither the current label
    \ifx#1\deadref\else % nor \deadref, reference was wrong: warn!
     \immediate\write16{Warning:
       forward reference to \string#1 is incorrect, run again!
       Changed from #1 to \lastlabel.}%
  \fi\fi
  \writefrw{\string\forward\string#1={\lastlabel}}% write to file
  \global\let#1\lastlabel \deflab#1% declare 
  \letnext\ignorespaces % do nothing particular afterward
 \else % #1 was not forward referenced, maybe evaporated forward ref
  \ifremoved#1\intoks\oldforward % forget unused forward definition
   \global\let#1\lastlabel \deflab#1\letnext\ignorespaces % as above
  \else\defnext{#1}% #1 is not a reference at all, execute it normally
  \fi
 \fi
}
% afterwards, \next is defined as \ignorespaces if #1 was used as a symbol
% if not, we have set \defnext{#1}

\def\labelsec#1%
{\ifx\undefined#1% if label was undefined, just define it
 \global\let#1\lastlabel \deflab#1%
 \else\definefrw{#1}%
 \fi
 \nobreak % prolong nobreak
}

\outer\def\proclaim#1. #2#3\par
  {\secmedbreak \global\advance\proclno by 1 \reset\itemno
   \edef\lastlabel{\ifnum-1<\secno    \number\secno.\fi
                   \ifnum-1<\subsecno \number\subsecno.\fi
                   \ifnum-1<\subsubsecno \number\subsubsecno.\fi
                   \number\proclno}
   {\interlinepenalty=300 \parskip=0pt\def\par{\endgraf\penalty200 }% (item)
    \noindent % Emit \baselineskip & \parskip glue BEFORE doing \deflab
    \if\let\noexpand#2% Is #2 a control sequence?
     \ifx\undefined#2% if so, and undefined, define it
      \global\let#2\lastlabel \deflab#2\letnext\ignorespaces
     \else \definefrw{#2}% if already defined, handle forward definition
     \fi
    \else\defnext{#2}% not a control sequence at all, incorporate into body
    \fi
    \bf\lastlabel.\enspace#1.\enspace \sl\next#3\endgraf
    }%
   \ifdim \lastskip<\medskipamount \removelastskip\penalty55\medskip \fi
   }

\def\neqn(#1#2)%
{\global\advance\eqnumber by 1 \edef\lastlabel{\number\eqnumber}%
  (\ifcat a\noexpand#1\errmessage{Incorrect label #1#2}\lastlabel
   \else
    \ifx\undefined#1\global\let#1\lastlabel \deflab#1%
    \else\definefrw{#1}%
    \fi
    \lastlabel\fam0 #2% #2 is possible tail, like "a,b"
   \fi
  )}

% -*- mode: tex; tex-main-file: "root.tex"; -*-

% definitions proper to Stanley Festschrift paper

\font\fivesl=cmsl5
\font\sevensl=cmsl7
\scriptfont\slfam=\sevensl \scriptscriptfont\slfam=\fivesl

\mathorddef\comp=C2 % {\cal C}, compositions
\mathorddef\Part=P2 % {\cal P}, partitions
\def\bcomp{\comp^{\scriptscriptstyle\set2}}
\mathorddef\Mat=M2 % {\cal M}, matrices of finite support
\def\bMat{\Mat^{\scriptscriptstyle\set2}}

\def\tr{^{\sl t}}
\def\repex^#1{^{(#1)}}

% next line uses (Open,Close)=(4,5) and \bffam=6 from plain.tex
\def\Kr#1{\mathchar"465B\,#1\,\mathchar"565D\,} % Iverson symbol

\def\set#1{{[#1]}}
\def\sing#1{\{#1\}} % singleton
\def\edgprm#1#2{\eps(#1,#2)}
\altopdef\sym{Sym}

\def\Sn{{\Sym_n}}

\def\[#1,#2]{{#1[\if!#2!\,\else#2\fi]}}
\def\d#1#2{{\delta_{#1}\if|#2|\else+#2\fi}}
\def\dn{\d n}

\def\scalar<#1,#2>{{\bigl<\,#1\bigm|#2\,\bigr>}}
\def\Scalar<#1,#2>{\left<\,#1\,\,\vrule\,\,#2\,\right>}

\def\Mtt(#1,#2,#3,#4){\({#1\atop#3}~{#2\atop#4}\)} % 2x2 matrix

\let\leh=\leftharpoondown
\let\geh=\rightharpoonup
\let\lev=\leftharpoonup
\let\gev=\rightharpoondown
\def\nleh{\,\not\nobreak\!\leh}
\def\nlev{\,\not\nobreak\!\lev}
\def\rlt#1{\prec_{{\rm r}(#1)}}

\opdef\Tabl
\def\mat#1#2{\Mat_{#1,#2}}
\def\bmat#1#2{\bMat_{#1,#2}}
\def\IE#1{\Tabl(#1)}
\def\BE#1{\Tabl^{\scriptscriptstyle\set2}(#1)}
\opdef\LR

\opdef\row \opdef\col
\def\rsum#1{\row(#1)}
\def\csum#1{\col(#1)}

\opdef\wt
\def\weight#1{\wt(#1)}
\opdef\cwt \opdef\rwt
\def\colweight#1#2{\cwt_{#1}(#2)}
\def\rowweight#1#2{\rwt_{#1}(#2)}
\opdef\SST
\def\SSTw#1#2{\SST(#1,#2)}
\def\SSTx#1{\SST(#1)}

\altopdef\height{ht}

\def\Xn{X_\set n} 
\def\spol#1{\Lambda_\set{#1}} \def\spoln{\spol n}
\def\apol#1{A_\set{#1}}       \def\apoln{\apol n}

\def\commas#1#2 {#1\if|#2|\else,\commas#2 \fi}
\def\s#1 {s_{(\commas#1 )}}
\def\m#1 {m_{(\commas#1 )}}
\def\h#1 {h_{(\commas#1 )}}
\def\e#1 {e_{(\commas#1 )}}
\font \titf=cmbx12

\refindent=3pc

\readrefs

\maketoc

\hyphenation{semi-standard}

\newif\ifblownup  % false here, but will be true in definitive version
\def\maybebreak{} % may inset a break in definitive version

\mark{\marktoks={}}% no running head on this page
\topglue \smallskipamount
\centerline{\titf Schur functions and alternating sums}
\medskip

\begingroup \ninept
\parindent=0pt \leftskip=0pt plus 1fil \rightskip=\leftskip \parfillskip=0pt
\parskip 0pt \obeylines
Marc A. A. van Leeuwen
Universit{\'e} de Poitiers, Math{\'e}matiques
BP 30179, 86962 Futuroscope Chasseneuil, France
email: {\tt Marc.van-Leeuwen@math.univ-poitiers.fr}
\endgroup
\medskip\centerline{\boldit
  Dedicated to Richard Stanley on the occasion of his $60\rm^{th}$ birthday}
\medskip

\medskip
{\narrower\ninept \centerline{\bf Abstract} \nobreak
\noindent We discuss several well known results about Schur functions that can
 be proved using cancellations in alternating summations; notably we shall
 discuss the Pieri and Murnaghan-Nakayama rules, the Jacobi-Trudi identity and
 its dual (Von N{\"a}gelsbach-Kostka) identity, their proofs using the
 correspondence with lattice paths of Gessel and Viennot, and finally the
 Littlewood-Richardson rule. Our our goal is to show that the mentioned
 statements are closely related, and can be proved using variations of the
 same basic technique. We also want to emphasise the central part that is
 played by matrices over~$\{0,1\}$ and over~$\N$; we show that the
 Littlewood-Richardson rule as generalised by Zelevinsky has elegant
 formulations using either type of matrix, and that in both cases it can be
 obtained by two successive reductions from a large signed enumeration of such
 matrices, where the sign depends only on the row and column sums of the
 matrix.
\par
}

\subsecno=-1 % no subsections
\secno=-1    % start with section 0

% -*- mode: tex; tex-main-file: "root.tex"; -*-

\newsection Introduction.

Many of the more interesting combinatorial results and correspondences in the
basic theory of symmetric functions involve Schur functions, or more or less
equivalently the notions of semistandard tableaux or horizontal strips. Yet
the introduction of these notions, in any of the many possible ways, is not
very natural when considering only symmetric functions, cf.~\ref{Stanley EC2,
7.10}. One way the importance of Schur functions can be motivated is by
representation theory: interpreting symmetric functions in the representation
theory either of the symmetric groups or the general linear groups, the Schur
functions correspond to the irreducible representations. However, there is
another way of motivating it: if one broadens the scope slightly from
symmetric functions to alternating polynomials, then Schur functions do arise
quite naturally as quotients of alternants. It is this point of view, which
could also be reached from representation theory if use is made only of Weyl's
character formula, that we shall take in this paper; from this perspective it
is not so surprising that proofs of basic identities involving Schur functions
should involve alternating summations and cancellations.

The main point we would like to make in this paper is that the use of the
definition of Schur functions as quotients of alternants can be limited to the
deduction of a single simple formula (lemma~\forward\Schurmul), which
describes the multiplication by an arbitrary symmetric function in the basis
of Schur functions; after this, alternating polynomials need not be considered
any more. In general the formula produces an alternating sum of Schur
functions; in various particular cases, one can obtain classical results from
it (the Pieri and Murnaghan-Nakayama rules, Jacobi and Von
N{\"a}gelsbach-Kostka identities, and the Littlewood-Richardson rule) by
judiciously applying combinatorially defined cancellations. Our presentation
is nearly self-contained, but we do use an enumerative identity that follows
from the \caps{RSK}-correspondence; we shall omit its well known elementary
combinatorial proof, which is not directly related to the theme of this paper.

Our paper is structured as follows. In~\Sec\forward\prelimsec\ we give the
basic definitions concerning symmetric functions, alternating polynomials and
(skew) Schur functions. In~\Sec\forward\Pierisec\ we introduce our basic
lemma, and its most elementary applications giving the Pieri and
Murnaghan-Nakayama rules. In~\Sec\forward\Cauchysec\ we first establish the
duality of the bases of complete and minimal symmetric functions (this is
where the \caps{RSK}-correspondence is used). This allows us to interpret
(skew) Schur functions as generating series of semistandard tableaux (which
elsewhere is often used as their definition), and to deduce the Cauchy, Jacobi
and Von N{\"a}gelsbach-Kostka identities. In~\Sec\forward\GVsec\ we discuss
cancellations defined for intersecting families of lattice paths, in the style
of Gessel and Viennot, and relate them to identities derived from the Pieri
rules. These considerations lead to natural encodings of families of lattice
paths, and of the semistandard tableaux that correspond to non-intersecting
families of paths, by matrices with entries in~$\N$ or in~$\{0,1\}$; these
encodings are also important in the sequel. In~\Sec\forward\LRsec\ we give a
final application of our basic lemma to derive the Littlewood-Richardson rule.
Formulating (Zelevinsky's generalisation of) that rule in terms of binary or
integral matrices reveals in both cases an unexpected symmetry. We also
exhibit an equally symmetrical doubly alternating expressions for the same
numbers, in which no tableaux appear at all. We close by raising a question
inspired by these expressions, which will be taken up in a sequel to this
paper.

% Local IspellDict: default
% -*- mode: tex; tex-main-file: "root.tex"; -*-

\subsecno=0
\newsection Preliminaries and definitions.
\labelsec\prelimsec

Studying symmetric functions involves the use of various combinatorial
objects; we start with some general considerations concerning those. We shall
make much use of sequences (vectors) and matrices, of which the entries will
almost always be natural numbers. In some cases the entries are restricted to
be either $0$~or~$1$, in which case we shall refer to the objects as
``binary''. While all objects we shall encounter can be specified using finite
information, we shall consider vectors and matrices as associations of entries
to indices, without restricting those indices to a finite set (just like for
polynomials one usually does not give an a priori bound for the degrees of
their monomials). Thus vectors and matrices are ``finitely supported'', in
that the entries are zero outside a finite range of indices; finite vectors
and matrices are identified with infinite ones obtained by extension with null
entries. This convention notably allows addition of vectors or matrices
without concern about their sizes.

When displaying matrices we shall as usual let the first index increase
downwards and the second to the right, and the same convention will be used
whenever subsets of $\N\times\N$ are displayed, such as Young diagrams (in the
sequel to this paper we shall in fact encounter Young diagrams in the role of
subsets of indices in matrices). Some objects, notably tableaux, are defined
as sequences of vectors; in this case the indices for the sequence are written
as parenthesised superscripts to avoid confusion with the subscripts
indexing individual vectors.

We always start indexing at~$0$, in particular this applies to sequences, rows
and columns of matrices and tableaux, and entries of tableaux. Hence in the
situation where a sequence of objects is determined by the intervals between
members of another sequence (such as horizontal strips in a semistandard
tableau, which are given by successive members of a sequence of shapes), the
index used for an interval is the same as that of the \emph{first} of the
members bounding it. Our standard $n$-element set is
$\set{n}=\setof{i\in\N}:i<n\endset$. For the set theoretic difference
$S\setminus{T}$ we shall write $S-T$ when it is known that $T\subset{S}$.

We shall frequently use the ``Iverson symbol'': for any Boolean expression
{\it condition} one puts
$$
  \Kr{{\it condition}}=\cases{1 & if {\it condition} is satisfied, \cr
                              0 & otherwise. \cr}
$$
This notation, proposed in \ref{concrete, p.~24}, and taken from the
programming language \caps{APL} by K.~Iverson, generalises the Kronecker delta
symbol: instead of $\delta_{i,j}$ one can write $\Kr{i=j}$. Among other uses,
this notation allows us to avoid putting complicated conditions below
summations to restrict their range: it suffices to multiply their summands by
one or more instances of $\Kr{{\it condition}}$. By convention, in a product
containing such a factor, the factors to its right are evaluated only if the
condition holds; if it fails, the product is considered to be~$0$ even if some
remaining factor should be undefined.

\subsection Compositions and partitions.

The most basic combinatorial objects we shall use are finitely supported
sequences of natural numbers $\alpha=(\alpha_i)_{i\in\N}$. The entries
$\alpha_i$ are called the parts of~$\alpha$, and the main statistic on such
sequences is the sum of the parts, written $|\alpha|=\sum_{i\in\N}\alpha_i$.
The systematic name for such sequences~$\alpha$ with $|\alpha|=d$ would be
infinite weak compositions of~$d$, but we shall simply call them just
\newterm{compositions} of~$d$. The set of compositions of~$d$ will be
denoted by $\comp_d$ (this set is infinite when $d>0$), and
$\comp=\Union_{d\in\N}\comp_d$ denotes the set of all compositions. In order
to denote specific compositions, we shall specify an initial sequence of their
parts, which are implicitly extended by zeroes. When the parts of a
composition are restricted to lie in $\{0,1\}=\set2$, it will be called a
\newterm{binary composition}; we define
$\bcomp_d=\setof\alpha\in\comp_d:\forall{i\in\N}\:\alpha_i\in\set2\endset$ and
$\bcomp=\Union_{d\in\N}\bcomp_d$. Binary compositions of~$d$ correspond to
$d$-element subsets of~$\N$, while arbitrary compositions of~$d$ correspond to
multisets of size~$d$ on~$\N$. Among other uses, compositions parametrise
monomials; if $X_\N=\setof{X_i}:i\in\N\endset$ is a countable set of commuting
indeterminates, then the monomial~$\prod_{i\in\N}X_i^{\alpha_i}$ will be
denoted by~$X^\alpha$.

We shall consider permutations of indeterminates, and correspondingly of the
parts of compositions. The group that acts is the group~$\Sym_\infty$ of
permutations of~$\N$ that fix all but finitely many numbers. The permutation
$\sigma\in\Sym_\infty$ acts by simultaneously substituting
$X_i:=X_{\sigma(i)}$ for all indeterminates, and therefore operates on
compositions by permuting their parts:
$\sigma(\alpha)=(\alpha_{\sigma\inv(i)})_{i\in\N}$. Obviously
$|\sigma(\alpha)|=|\alpha|$, and the orbit of~$\alpha$ contains a unique
composition whose parts are weakly decreasing, which will be denoted
by~$\alpha^+$; for instance for $\alpha=(0,5,2,0,0,1,7,0,2)$ one has
$\alpha^+=(7,5,2,2,1)$. For $d\in\N$ we define the finite set
$\Part_d=\setof\\\in\comp_d:\forall{i\in\N}\:\\_i\geq\\_{i+1}\endset$, whose
elements are called \newterm{partitions} of~$d$; then $\alpha^+\in\Part_d$ for
any~$\alpha\in\comp_d$. We also put $\Part=\Union_{d\in\N}\Part_d$. All binary
compositions of~$d$ form a single orbit under permutations of their parts, so
there is just a single binary partition of~$d$: it is the partition
$\Kr{i\in\set{d}}_{i\in\N}$ whose $d$ initial parts are~$1$ and the rest~$0$,
and we shall denote it by~$1\repex^d$.

We shall usually denote compositions by Greek letters $\alpha,\beta,\ldots$,
but for partitions we use Greek letters further on in the alphabet: $\\$,
$\mu$, $\nu$, and sometimes~$\kappa$. Apart from listing its nonzero parts, a
partition $\\\in\Part$ can also be specified by drawing its \newterm{diagram}
$\set\\=\setof(i,j)\in\N^2:j\in\set{\\_i}\endset$. Elements of the diagram are
drawn (and usually referred to) as squares, so that for instance the diagram
of $\\=(7,5,2,2,1)$ would be drawn as
$$
  \smallsquares\Young(,,,,,,|,,,,|,|,|).
$$
The transpose partition of~$\\\in\Part$, which will be denoted by~$\\\tr$, is
the one whose parts give the lengths of the columns of~$\set\\$, so that
$\set{\\\tr}$ is the transpose diagram $\set\\\tr$; one has
$\\\tr_j=\Card\setof{i\in\N}:j\in\set{\\_i}\endset$.

We shall be considering several relations defined between partitions; we
collect their definitions here. The most fundamental relation is the partial
ordering~`$\subset$' defined by inclusion of diagrams: $\mu\subset\\$ means
that $\set\mu\subset\set\\$ or equivalently that $\mu_i\leq\\_i$ for
all~$i\in\N$. Note that if $\mu\subset\\$ then $\\-\mu$ and $\\\tr-\mu\tr$ are
compositions. The relation~`$\subset$' will be used mostly implicitly via the
notion of a \newterm{skew shape}~$\\/\mu$, which denotes the interval from
$\mu$ to~$\\$ in the poset~$(\Part,\subset)$; the corresponding skew diagram
is $\set{\\/\mu}=\set\\-\set\mu$, and we define $|\\/\mu|=|\\|-|\mu|$. Several
relations refining~`$\subset$' will be used; for the ones in the following
definition it will be convenient to define them on the set of all
compositions, although they will never hold unless both arguments are actually
partitions.

\proclaim Definition. \levhdef
The relations~`$\lev$' and~`$\leh$' on~$\comp$ are defined as follows. To have
either $\mu\lev\\$ or $\mu\leh\\$, it is necessary that $\\/\mu$ be a skew
shape (in other words $\\,\mu\in\Part$, and $\mu\subset\\$). If this is the
case, then $\mu\lev\\$ holds if and only if $\\-\mu\in\bcomp$, in which case
$\\/\mu$ is called a \newterm{vertical strip}; similarly $\mu\leh\\$ holds if
and only if $\\_{i+1}\leq\mu_i\leq\\_i$ for all~$i\in\N$, in which case
$\\/\mu$ is called a \newterm{horizontal strip}.

Note that the final condition for $\mu\leh\\$ already implies that $\\/\mu$ is
a skew shape; in addition it means that $\set{\\/\mu}$ has at most one square
in any column. Similarly, for a skew shape~$\\/\mu$, the condition $\mu\lev\\$
means that $\set{\\/\mu}$ has at most one square in any row. Therefore
$\mu\leh\\$ is equivalent to $\mu\tr\lev\\\tr$ when $\\,\mu\in\Part$. To
denote the opposite relations we shall rotate rather than reflect the symbol,
so $\\\geh\mu$ means the same as $\mu\leh\\$, while $\\\gev\mu$ means the same
as $\mu\lev\\$.

We illustrate concrete instances of these relations graphically by
superimposing the contours of the diagrams of the two partitions involved:
\bigdisplay
\eqalign
{ (7,5,2,2,1)\lev(8,6,3,3,1,1,1)
 &\:\qquad
  \hbox{$\smallsquares{\rtab1
           \edgeseq5,0  ;111110000000
           \edgeseq7,0;110101100010010
	   \edgeseq7,0;011100110001001 }
        $}
;\cr
\noalign{\smallskip}
  (7,5,2,2,1)\leh(11,6,4,2,1,1)
 &\:\qquad
  \hbox{$\smallsquares{\rtab1
           \edgeseq5,0 ;111110000000
           \edgeseq6,0;10101100010010000
	   \edgeseq6,0;01101001001000001 }
        $}
.\cr
}
$$

For the following definition we use the partitioning of~$\N\times\N$ into
\newterm{diagonals}~$D_d$, for~$d\in\Z$:
$$
  D_d=\setof(i,j)\in\N^2:j-i=d\endset
\ttex{for $d\in\Z$}.
\label(\diagonaldef)
$$
% Diagonals are called consecutive if the indices parametrising them are
% consecutive in~$\Z$.

\proclaim Definition.
For $k>0$, a relation `$\rlt{k}$' on~$\Part$ is defined as follows:
$\mu\rlt{k}\\$ means that $\\/\mu$ is a skew shape with $|\\/\mu|=k$, for
which the $k$ squares of~$\set{\\/\mu}$ lie on $k$~consecutive diagonals. In
this case we call the shape $\\/\mu$ a $k$-ribbon. The height
$\height(\\/\mu)$ of a $k$-ribbon $\\/\mu$ is the difference between the
initial (row) coordinates of the squares of~$\set{\\/\mu}$ on the first and
the last of those $k$~diagonals.

Again we give a graphic illustration in the same style as before, for
$\\/\mu=(7,6,6,3,3,2)/(7,5,2,2,1)$:
$$
  (7,5,2,2,1)\rlt{10}(7,6,6,3,3,2)
 \:\qquad
  \hbox{$\smallsquares{\rtab1 
           \edgeseq5,0 ;111110000000
           \edgeseq6,0;1010110001001
	   \edgeseq6,0;00101100011 }
        $}
.
$$
One has
$\set{\\/\mu}
 =\{(5,0),(5,1),(4,1),(4,2),(3,2),(2,2),(2,3),(2,4),(2,5),(1,5)\}$,
which diagram has its squares on the $10$ consecutive diagonals $D_d$ for
$-5\leq{d}<5$; moreover, we see that $\height(\\/\mu)=5-1=4$.

Finally we shall need the dominance partial ordering on each set~$\Part_d$
separately.

\proclaim Definition.
For any fixed~$d\in\N$ a relation~`$\leq$' on~$\Part_d$, called the dominance
ordering, is defined by $\mu\leq\\$ if and only if for every $k\in\N$ one has
$\sum_{i\in\set{k}}\mu_i\leq\sum_{i\in\set{k}}\\_i$.

\subsection Matrices and tableaux.

We shall use the two-dimensional counterparts of compositions: finitely
supported matrices with entries in~$\N$. Like for compositions the binary
case, where entries are restricted to $\set2=\{0,1\}$, will be of special
interest. The statistic given by the sum of all entries can be refined by
taking sums separately either of rows or of columns; in either case the result
is a composition.

\proclaim Definition.
Let $\Mat$ denote the set of matrices $(M_{i,j})_{i,j\in\N}$ with entries
$M_{i,j}$ in~$\N$, of which only finitely many are nonzero, and let $\bMat$
denote its subset of \newterm{binary matrices}, those of which all entries lie
in~$\{0,1\}$. Let $\row\:\Mat\to\comp$ be the map
$M\mapsto(\sum_{j\in\N}M_{i,j})_{i\in\N}$ that takes row sums, and
$\col\:\Mat\to\comp$ the map $M\mapsto(\sum_{i\in\N}M_{i,j})_{j\in\N}$ that
takes column sums; put
$\mat\alpha\beta=\setof{M\in\Mat}:\rsum{M}=\alpha,\csum{M}=\beta\endset$ and
$\bmat\alpha\beta=\mat\alpha\beta\thru\bMat$ for $\alpha,\beta\in\comp$.

These matrices can be used to record sequences of (binary) compositions with
finite support, either by rows or by columns. We shall denote row~$i$ of~$M$
by~$M_i$, and column~$j$ by~$M\tr_j$. We shall also need sequences of
partitions, but these will be subject to the condition that adjacent terms
differ by horizontal or vertical strips, and the condition of finite support
is replaced by the sequence becoming ultimately stationary. This gives rise to
the notion of semistandard tableau, and some variants of it.

\proclaim Definition. \SSTdef
Let $\\/\mu$ be a skew shape, and $\alpha\in\comp$. A semistandard tableau of
shape $\\/\mu$ and weight~$\alpha$ is a sequence of partitions
$(\\^{(i)})_{i\in\N}$ with $\\^{(i)}\leh\\^{(i+1)}$ and
$|\\^{(i+1)}/\\^{(i)}|=\alpha_i$ for all $i\in\N$, $\\^{(0)}=\mu$, and
$\\^{(N)}=\\$ for any $N$ that is so large that $\alpha_i=0$ for all
$i\geq{N}$. The weight of a tableau~$T$ is denoted by~$\weight T$, and the set
of all semistandard tableaux of shape~$\\/\mu$ by~$\SSTx{\\/\mu}$; we also put
$\SSTw{\\/\mu}\alpha=\setof{T\in\SSTx{\\/\mu}}:\weight{T}=\alpha\endset$. A
transpose semistandard tableau of shape $\\/\mu$ and weight~$\alpha$ is a
sequence of partitions defined similarly, with $\\^{(i)}\lev\\^{(i+1)}$
replacing $\\^{(i)}\leh\\^{(i+1)}$.

We shall reserve the qualification ``Young tableau'' to the case $\mu=(0)$, in
which case $\\/\mu$ will be abbreviated to~$\\$ in the notations just
introduced. There are maps from semistandard tableaux to transpose
semistandard tableau and vice versa, defined by transposing each partition in
the sequence; under these maps the shape of the tableau is transposed while
the weight is preserved. Another variation on the notion of semistandard
tableau is to replace the relations $\\^{(i)}\leh\\^{(i+1)}$ or
$\\^{(i)}\lev\\^{(i+1)}$ by their opposite relations $\\^{(i)}\geh\\^{(i+1)}$
respectively $\\^{(i)}\gev\\^{(i+1)}$. This gives the notions of reverse
(transpose) semistandard tableaux, which will occur in the sequel to this
paper; their shape $\\/\mu$ and weight~$\alpha$ are such that the sequence
starts at~$\\=\\^{(0)}$ and ultimately becomes $\mu$, while
$|\\^{(i)}/\\^{(i+1)}|=\alpha_i$ for all $i\in\N$.

The traditional way to display a semistandard tableau is to draw the diagram
of its shape filled with numbers, which identify for each square the
horizontal strip to which it belongs. We shall label with an entry~$i$ the
squares of $\set{\\^{(i+1)}/\\^{(i)}}$. The entries will then increase weakly
along rows, and increase strictly down columns, and for this reason
semistandard tableaux are also called column-strict tableaux (and transpose
semistandard tableaux are then called row-strict tableaux). Thus the
semistandard tableau $T=\((4,1)\leh (5,2)\leh (5,3,2)\leh (6,3,3,1)\leh
(6,4,3,2)\leh (7,5,4,3)\leh (9,5,5,3,1)\leh (9,8,5,5,3)\leh
(9,8,5,5,3)\leh\nobreak\cdots\)$, which is of shape $(9,8,5,5,3)/(4,1)$ and
weight~$(2,3,3,2,4,4,7)$, will be displayed as
\bigdisplay
T=
\Skew(4:0,2,4,5,5|1:0,1,3,4,6,6,6|0:1,1,2,4,5|0:2,3,4,6,6|0:5,6,6)
.\label(\Tdef)
$$
More important in our paper than this display will be two ways of representing
tableaux by matrices. Simply recording the partitions forming a tableau~$T$ in
the rows or columns of a matrix does not give a finitely supported matrix, but
we can obtain one by recording the differences between successive partitions.
We shall call the matrix so obtained an~\newterm{encoding} of~$T$, but one
should realise that decoding the matrix to reconstruct~$T$ requires knowledge
of at least one of the partitions forming the (skew) shape of~$T$. Various
ways are possible to record horizontal strips $\\^{(i+1)}/\\^{(i)}$: one may
either record the differences $\\^{(i+1)}-\\^{(i)}\in\comp$ or the differences
between the transpose shapes $(\\^{(i+1)})\tr-(\\^{(i)})\tr\in\bcomp$, and one
may record these compositions either in the rows or the columns of the matrix.
From the four possible combinations we choose the two for which one has a
correspondence either between the rows of the tableau and the rows of the
matrix, or between the columns of the tableau and the columns of the matrix.

\proclaim Definition. \encodingdef
Let $T=(\\^{(i)})_{i\in\N}$ be a semistandard tableau. The integral encoding
of~$T$ is the matrix $M\in\Mat$ defined by $M_{i,j}=(\\^{(j+1)}-\\^{(j)})_i$,
and the binary encoding of~$T$ is the matrix $M'\in\bMat$ defined by
$M'_{i,j}=((\\^{(i+1)})\tr-(\\^{(i)})\tr)_j$. The sets of integral and binary
encodings of semistandard tableaux $T\in\SSTx{\\/\mu}$ will be denoted by
$\IE{\\/\nu}$ and $\BE{\\/\mu}$, respectively.

For instance for the tableau~$T$ of~(\Tdef), one finds the integral and binary
encodings
\bigdisplay
M=
\pmatrix{1&0&1&0&1&2&0\cr
         1&1&0&1&1&0&3\cr
         0&2&1&0&1&1&0\cr
	 0&0&1&1&1&0&2\cr
	 0&0&0&0&0&1&2\cr}
\ifblownup \text{and} \else \ttext{and}\fi
M'=
\pmatrix{0&1&0&0&1&0&0&0&0\cr
         1&1&1&0&0&0&0&0&0\cr
         1&0&1&0&0&1&0&0&0\cr
	 0&1&0&1&0&0&0&0&0\cr
	 0&0&1&1&1&0&1&0&0\cr
	 1&0&0&0&1&0&0&1&1\cr
	 0&1&1&1&1&1&1&1&0\cr}
,\nn
$$
which finite matrices must be thought of as extended indefinitely by zeroes.
To reconstruct from either of these matrices the tableau~$T$ or the other
matrix, one must in addition know at least that $\mu=(4,1)$ or that
$\\=(9,8,5,5,3)$ for the shape~$\\/\mu$ of~$T$. Each entry $M_{i,j}$ counts
the number of entries~$j$ in row~$i$ of the displayed form of~$T$, while entry
$M'_{i,j}$ counts the number (at most one) of entries~$i$ in column~$j$.
Therefore the row $M_i$ records the weight of row~$i$ of the display of~$T$,
while the column $(M')\tr_j$ records the weight of its column~$j$. One has
$\rsum{M}=\\-\mu$, \ $\csum{M'}=\\\tr-\mu\tr$, and
$\csum{M}=\rsum{M'}=\weight{T}$.

\subsection Symmetric functions.

There are several equivalent ways to define the ring~$\Lambda$ of symmetric
functions. Following \ref{Stanley EC2}, we shall realise $\Lambda$ as a
subring of the ring $\Z[[X_\N]]$ of power series in infinitely many
indeterminates. The elements~$f$ of this subring are characterised by the fact
that the coefficients in~$f$ of monomials $X^\alpha,X^\beta$ are the same
whenever $\alpha^+=\beta^+$ (so $f$ is stable under the action
of~$\Sym_\infty$), and that the degree of monomials with nonzero coefficients
in~$f$ is bounded. Elements~$f\in\Lambda$ are called \newterm{symmetric
functions}. Since the indicated subring of~$\Z[[X_\N]]$ is just one
realisation of~$\Lambda$, we make a notational distinction between occurrences
of a symmetric function~$f$ that are independent of any realisation
of~$\Lambda$ (for instance in identities internal to~$\Lambda$), and
occurrences where the realisation inside~$\Z[[X_\N]]$ is essential (because
indeterminates $X_i$ occur explicitly in the same equation); in the latter
case we shall write $f[X_\N]$ instead of~$f$. If the nonzero coefficients of
$f[X_\N]$ only occur for monomials of degree~$d$, then $f$ is called
homogeneous of degree~$d$; due to the required degree bound, this makes
$\Lambda$ into a graded ring.

Another realisation of~$\Lambda$ is via its images in polynomial rings in
finite sets of indeterminates. This is for instance the point of view taken
in~\ref{Macdonald}; for us this realisation is important in order to be able
to consider alternating expressions, which is hard to do for infinitely many
indeterminates. For any $n\in\N$, let $\Xn=\setof{X_i}:i\in\set{n}\endset$ be
the set of the first~$n$ indeterminates. There is a ring morphism
$\Z[[X_\N]]\to\Z[[\Xn]]$ defined by setting $X_i:=0$ for all $i\geq{n}$, and
the image of the subring~$\Lambda$ under this morphism is the subring of the
symmetric polynomials in~$\Z[\Xn]$, those invariant under all permutations of
the indeterminates; we shall denote this image by~$\spoln$. For $f\in\Lambda$,
the image in $\spoln$ of $f[X_\N]\in\Z[[X_\N]]$ will be denoted by~$f[\Xn]$.
Thus each $f\in\Lambda$ gives rise to a family $(f[\Xn])_{n\in\N}$ of elements
$f[\Xn]\in\spoln$ of bounded degree, which family is coherent with respect to
the projections $\spol{n+1}\to\spoln$ defined by the substitution $X_n:=0$. We
shall write this final property as $f[X_\set{n+1}][X_n:=0]=f[\Xn]$ for
all~$n\in\N$. Conversely each family $(f_n)_{n\in\N}$ with $f_n\in\spoln$ for
all~$n\in\N$ that satisfies $f_{n+1}[X_n:=0]=f_n$ for all~$n$, and for which
$\deg{f_n}$ is bounded, forms the set of images of a unique element
$f\in\Lambda$. In other words, one can realise $\Lambda$ as the inverse limit
in the category of graded rings of the system $(\spoln)_{n\in\N}$ relative to
the given projections $\spol{n+1}\to\spoln$.

For any $\alpha\in\Part_d$, the sum $m_\alpha[X_\N]
=\sum_{\beta\in\comp_d}\Kr{\alpha^+ =\beta^+}X^\beta$ of all distinct
monomials in the permutation orbit of~$X^\alpha$ is a symmetric function.
Since no nonempty proper subset of its nonzero terms defines a symmetric
function, we shall call $m_\alpha$ a \newterm{minimal symmetric function} (we
avoid the more traditional term ``monomial'' symmetric function since the set
of all~$m_\alpha$ is not closed under multiplication). The set
$\setof{m_\\}:\\\in\Part_d\endset$ is a basis of the additive group of
homogeneous symmetric functions of degree~$d$.

The \newterm{elementary symmetric functions} $e_d$ for $d\in\N$ are instances
of minimal symmetric functions: they are defined as $e_d=m_{1\repex^d}$. One
can write more explicitly
$$
  e_d[X_\N]=\sum_{\alpha\in\bcomp_d}X^\alpha
 =\sum_{i_1,\ldots,i_d\in\N}\Kr{\fold<(i1..d)}\fold(Xi_1..i_d)
.\label(\edeq)
$$
The \newterm{complete (homogeneous) symmetric functions} $h_d$ for $d\in\N$ are
defined by $h_d=\sum_{\\\in\Part_d}m_\\$. Like the elementary symmetric
functions, they can be written  more explicitly
$$
  h_d[X_\N]=\sum_{\alpha\in\comp_d}X^\alpha
 =\sum_{i_1,\ldots,i_d\in\N}\Kr{\fold\leq(i1..d)}\fold(Xi_1..i_d)
.\label(\hdeq)
$$
The \newterm{power sum symmetric functions} $p_d$ for $d>0$ are defined by
$p_d=m_{(d)}$, so $p_d[X_\N]=\sum_{i\in\N}X_i^d$. These families of symmetric
functions have the following generating series, expressed in $\Z[[X_\N,T]]$.
$$ \eqalignno
{ \sum_{d\in\N}e_d[X_\N]T^d&=\prod_{i\in\N}(1+X_iT)
,& \label(\eseries) \cr
  \sum_{d\in\N}h_d[X_\N]T^d&=\prod_{i\in\N}\left(\sum_{k\in\N}(X_iT)^k\right)
                            =\prod_{i\in\N}{1\over1-X_iT}
,& \label(\hseries) \cr
  \sum_{k>0}p_k[X_\N]T^k   &=\sum_{i\in\N}\left(\sum_{k>0}(X_iT)^k\right)
			    =\sum_{i\in\N}{X_iT\over1-X_iT}
.& \nn \cr
}
$$
For any $\alpha\in\comp$ we define $e_\alpha=\prod_{i\in\N}e_{\alpha_i}$ and
$h_\alpha=\prod_{i\in\N}h_{\alpha_i}$; since $e_0=h_0=1$ the infinite products
converge, and it is clear by commutativity that $e_\alpha=e_{\alpha^+}$ and
$h_\alpha=h_{\alpha^+}$. The products $e_\alpha$ and~$h_\alpha$ can be
expanded into monomials combinatorially, in terms of binary respectively
integral matrices: by multiplying together copies of the first equality
in~(\edeq) respectively in~(\hdeq), one finds
$$
\eqalignno
{ e_\beta[X_\N]&=\sum_{\alpha\in\comp}\Card{\bmat\alpha\beta}X^\alpha
&\label(\ebetaeq) \cr
  h_\beta[X_\N]&=\sum_{\alpha\in\comp}\Card{\mat\alpha\beta}X^\alpha
 &\label(\hbetaeq) \cr
}
$$
We can obtain generating series in~$\Z[[X_\N,Y_\N]]$ in which all $e_\beta$ or
all $h_\beta$ appear, either from the preceding equations, or by substituting
$T:=Y_j$ into copies of (\eseries) or~(\hseries) for $j\in\N$ and multiplying
them, giving
$$
\eqalignno
{  \sum_{\beta\in\comp}e_\beta[X_\N]Y^\beta
  =\sum_{M\in\bMat}X^{\rsum{M}}Y^{\csum{M}}
 &=\prod_{i,j\in\N}(1+X_iY_j)
& \label(\ebetaseries) \cr
   \sum_{\beta\in\comp}h_\beta[X_\N]Y^\beta
  =\sum_{M\in\Mat}X^{\rsum{M}}Y^{\csum{M}}
 &=\prod_{i,j\in\N}{1\over1-X_iY_j}
.& \label(\hbetaseries) \cr
}
$$

\subsection Alternating polynomials and Schur functions.

Now fix $n\in\N$, and let $\apoln$ denote the additive subgroup of $\Z[\Xn]$
of alternating polynomials, i.e., of polynomials~$p$ such that for all
permutations~$\sigma\in\Sn$ the permutation of indeterminates given
by~$\sigma$ operates on~$p$ as multiplication be the sign~$\eps(\sigma)$.
Multiplying an alternating polynomial by a symmetric polynomial gives another
alternating polynomial, so if we view $\Z[\Xn]$ as a module over its
subring~$\spoln$, then it contains $\apoln$ as a submodule. Like for symmetric
polynomials, the condition of being an alternating polynomial can be expressed
by comparing coefficients of monomials in the same permutation orbit: a
polynomial $\sum_{\alpha\in\N^n}c_\alpha{X^\alpha}$ is alternating if and only
if for every $\alpha\in\N^n$ and $\sigma\in\Sn$ one has
$c_{\sigma\cdot\alpha}=\eps(\sigma)c_\alpha$. In particular this implies that
$c_\alpha=0$ whenever $\alpha$ is fixed by any odd permutation, which happens
as soon as $\alpha_i=\alpha_j$ for some pair~$i\neq{j}$. In the contrary case,
$\alpha$ is not fixed by any non-identity permutation, and the alternating
orbit sum $a_\alpha[\Xn]=\sum_{\sigma\in\Sn}\eps(\sigma)X^{\sigma\cdot\alpha}$
is an alternating polynomial that is minimal in the sense that its nonzero
coefficients are all~$\pm1$ and no nonempty proper subset of its nonzero terms
defines an alternating polynomial. The element $a_\alpha[\Xn]$ is called an
\newterm{alternant}, and can be written as a determinant
$$
a_\alpha[\Xn]=\det\(X_j^{\alpha_j}\)_{i,j\in\set{n}}
=\left|\matrix
{X_0^{\alpha_0}&X_0^{\alpha_1}&\cdots&X_0^{\alpha_{n-1}}\cr
 X_1^{\alpha_0}&X_1^{\alpha_1}&\cdots&X_1^{\alpha_{n-1}}\cr
 \vdots & \vdots & \ddots & \vdots \cr
 X_{n-1}^{\alpha_0}&X_{n-1}^{\alpha_1}&\cdots&X_{n-1}^{\alpha_{n-1}}\cr
}\right|
;\label(\adef)
$$
we define $a_\alpha[\Xn]$ by the same expression even when $\alpha$ is fixed
by some transposition, but in that case it is~$0$. The set of alternants
generates $\apoln$ as an additive group, but to obtain a $\Z$-basis one must
remove the null alternants, and for all other orbits of compositions choose
one of the two opposite alternants associated to it. Thus one finds the
$\Z$-basis $\setof{a_\alpha[\Xn]}:\alpha\in\N^n;\fold>(\alpha0..n-1)\endset$
of~$\apoln$. Our convention of interpreting finite vectors by extension with
zeroes as finitely supported ones, allows us to view $\N^n$ as a subset
of~$\comp$. Then putting $\dn{}=(n-1,n-2,\ldots,1,0)\in\N^n$ the above basis
of~$\apoln$ can be written as
$\setof{a_\dn\\[\Xn]}:\\\in\Part\thru\N^n\endset$.

Put $\Delta_n=a_\dn{}[\Xn]$; in other words, $\Delta_n\in\apoln$ is the
Vandermonde determinant, which evaluates to $\prod_{0\leq{i}<j<n}(X_i-X_j)$.
Alternating polynomials are all divisible by each factor $X_i-X_j$, and
therefore by~$\Delta_n$. So viewing $\apoln$ as an $\spoln$-module, it is
cyclic with generator~$\Delta_n$.
The map $\spoln\to\apoln$ of multiplication by~$\Delta_n$ is a $\Z$-linear
bijection, so one can apply its inverse to the basis of $\apoln$ consisting of
elements $a_\dn\\[\Xn]$. Thus defining
$$
  s_\\[\Xn]={a_\dn\\[\Xn]\over\Delta_n}\in\spoln
\nn
$$
for $\\\in\Part$ with $\\_n=0$, the set
$\setof{s_\\[\Xn]}:\\\in\Part;\\_n=0\endset$ forms a $\Z$-basis of~$\spoln$.
It is useful to define $s_\alpha[\Xn]$ for arbitrary $\alpha\in\N^n$ by the
same formula. Doing so does not introduce any new symmetric functions, since
one has $s_\alpha[\Xn]=0$ unless the $n$ components of $\dn\alpha$ are all
distinct, and in that case one has $s_\alpha[\Xn]={\eps}s_\\[\Xn]$, where
$\\\in\Part$ with $\\_n=0$ is determined by the condition
$(\dn\alpha)^+=\dn\\$, and $\eps$~is the sign of the (unique) permutation
$\sigma\in\Sn$ such that $\sigma(\dn\alpha)=\dn\\$. Finally we extend this
definition to any (infinite) composition $\alpha\in\comp$, by defining
$s_\alpha[\Xn]=0$ whenever $\alpha\notin\N^n$.

Here are some examples illustrating these definitions. One has
\bigdisplay \interdisplaylinepenalty=10000
\eqalignno
{s_{(3,1)}[X_\set2]
&={a_{(4,1)}[X_\set2]\over\Delta_2}
 ={X_0^4X_1-X_0X_1^4\over{X_0-X_1}}
 =X_0^3X_1+X_0^2X_1^2+X_0X_1^3
\cr&
 =m_{(3,1)}[X_\set2]+m_{(2,2)}[X_\set2]
\cr
 s_{(1,2)}[X_\set2]
&={a_{(2,2)}[X_\set2]\over\Delta_2}=0
\cr
 s_{(0,4)}[X_\set2]
&={a_{(1,4)}[X_\set2]\over\Delta_2}
 ={X_0X_1^4-X_0^4X_1\over{X_0-X_1}}=-s_{(3,1)}[X_\set2]
\cr\noalign{\vfil\penalty500\vfilneg}
 s_{(3,1)}[X_\set3]
&={a_{(5,2,0)}[X_\set3]\over\Delta_3}
 ={X_0^5X_1^2-X_0^2X_1^5-X_0^5X_2^2+X_1^5X_2^2+X_0^2X_2^5-X_1^2X_2^5
   \over(X_0-X_1)(X_0-X_1)(X_1-X_2)}
\cr
&=X_0^3X_1+X_0^3X_2
 +X_0^2X_1^2+2X_0^2X_1X_2+X_0^2X_2^2
\cr&\quad{}
 +X_0X_1^3+2X_0X_1^2X_2+2X_0X_1X_2^2+X_0X_2^3
 +X_1^3X_2+X_1^2X_2^2+X_1X_2^3
\cr
&=m_{(3,1)}[X_\set3]+m_{(2,2)}[X_\set3]+2m_{(2,1,1)}[X_\set3]
\cr
}
$$

\proclaim Proposition.
For all $n\in\N$ and all $\alpha\in\comp$ one has
$s_\alpha[X_\set{n+1}][X_n:=0]=s_\alpha[X_\set{n}]$.

For instance one sees in the example above that
$s_{(3,1)}[X_\set3][X_2:=0]=s_{(3,1)}[X_\set2]$: from the definition one has
$m_\alpha[X_\set{n+1}][X_n:=0]=m_\alpha[X_\set{n}]$ for all~$\alpha\in\comp$,
while $m_\alpha[X_\set{n}]=0$ unless $\alpha\in\N^n\ssubset\comp$, so in
particular $m_{(2,1,1)}[X_\set2]=0$.

\proof
The value $s_\alpha[X_\set{n+1}][X_n:=0]$ can be computed by applying the
substitution $X_n:=0$ separately to the numerator
$a_\d{n+1}\alpha[X_\set{n+1}]$ and the denominator $\Delta_{n+1}$ in the
definition of $s_\alpha[X_\set{n+1}]$, provided that the latter substitution
yields a nonzero value; this is the case since
$\Delta_{n+1}[X_n:=0]=\fold(X0..n-1)\Delta_n$. We may assume that
$\alpha\in\N^{n+1}$ holds, since otherwise both $s_\alpha[X_\set{n+1}]$
and~$s_\alpha[X_\set{n}]$ are zero by definition. Now put
$\beta=\d{n+1}\alpha=(n+\alpha_0,\ldots,1+\alpha_{n-1},\alpha_n)\in\N^{n+1}$,
so that the mentioned numerator is~$a_\beta[X_\set{n+1}]$. If one has
$\alpha\notin\N^{n}$, so that $s_\alpha[\Xn]=0$ by definition, then
$\alpha_n\neq0$, and all of the (first $n+1$) components of~$\beta$ are
nonzero; in this case the substitution $X_n:=0$ kills all terms of the
numerator, so that $s_\alpha[X_\set{n+1}][X_n:=0]=0$. On the other hand if
$\alpha\in\N^{n}$ then also $\beta\in\N^n$, and in this case
$a_\beta[X_\set{n+1}][X_n:=0] =a_\beta[X_\set{n}]
=\fold(X0..n-1)a_\dn\alpha[X_\set{n}]$; after simplification of the
substituted numerator and denominator by $\fold(X0..n-1)$, one obtains
$s_\alpha[X_\set{n+1}][X_n:=0]=s_\alpha[X_\set{n}]$ as desired.
\QED

Thus for fixed $\alpha\in\comp$, the families $(s_\alpha[\Xn])_{n\in\N}$ of
symmetric polynomials are coherent with respect to the projections
$\spol{n+1}\to\spoln$ defined by the substitution $X_n:=0$. This property
allows the following definition, which we already anticipated in our notation.

\proclaim Definition. \Schurfuncdef
For $\alpha\in\comp$, the symmetric function $s_\alpha$ is the unique element
of~$\Lambda$ whose image in~$\spoln$ under the substitutions $X_i:=0$ for all
$i\geq{n}$ is $s_\alpha[\Xn]
=\Kr{\alpha\in\N^n}{a_\dn\alpha[X_\set{n}]\over{a_\dn{}[X_\set{n}]}}$, for
all~$n\in\N$.

The set $\setof{s_\\}:\\\in\Part\endset$ forms a $\Z$-basis of~$\Lambda$,
whose elements are called \newterm{Schur functions}. They are the central
subject of this paper, and we shall now introduce several notations to
facilitate their study. Firstly we shall denote by
$\scalar<\cdotp,\,\cdotp\!>$ the scalar product on~$\Lambda$ for which the
basis of Schur functions is orthonormal. Thus one has
$$
  f=\sum_{\\\in\Part}\scalar<f,s_\\>s_\\ \ttex{for any~$f\in\Lambda$}
.\label(\scalarprodid)
$$
The operation of multiplication by a fixed Schur function~$s_\mu$ has an
adjoint operation $s_\mu^*$ for this scalar product, i.e., which satisfies
$\scalar<s_\mu^*(f),g>=\scalar<f,s_\mu{g}>$ for all~$f,g\in\Lambda$; in terms
of this the \newterm{skew Schur functions} are defined by
$$
  s_{\\/\mu}=s_\mu^*(s_\\)
.\label(\skewSchurdef)
$$
They typically arise when one expresses the multiplication by a fixed symmetric
function in the basis of Schur functions, as skew Schur functions are
characterised by
$$
  \scalar<s_\mu{f},s_\\>=\scalar<f,s_{\\/\mu}>
\ttex{for all $\\,\mu\in\Part$ and $f\in\Lambda$}
.\nn
$$
For $f=h_\alpha$ and $f=e_\alpha$ these scalar products are of particular
interest, and are called \newterm{Kostka numbers}.

\proclaim Definition.
For $\mu,\\\in\Part$ and $\alpha\in\comp$, we set
$K_{\\/\mu,\alpha} =\scalar<h_\alpha,s_{\\/\mu}>$ and
$K'_{\\/\mu,\alpha}=\scalar<e_\alpha,s_{\\/\mu}>$.

When $\mu=(0)$, we abbreviate $K_{\\/\mu,\alpha}$ to
$K_{\\,\alpha}$ and $K'_{\\/\mu,\alpha}$ to $K'_{\\,\alpha}$. Then as a special
case of (\scalarprodid) one has
$$
  h_\alpha=\sum_{\\\in\Part}K_{\\,\alpha}s_\\
\ttext{and}
  e_\alpha=\sum_{\\\in\Part}K'_{\\,\alpha}s_\\
.\label(\heKostkaeq)
$$
We shall later give combinatorial descriptions of the Kostka numbers
(corollary~\forward\Kostkacor), from which it will become clear that they are
non-negative (this also follows from representation theoretic considerations),
and are related by $K'_{\\/\mu,\alpha}=K_{\\\tr/\mu\tr,\alpha}$.

Computing scalar products $\scalar<f,s_\\>$ directly from the definition is
usually quite hard. But if $f=s_\alpha$ with $\alpha\in\comp$, then $f$ can be
considered as slightly generalised Schur function: either $f=0$, or
$f=\pm{s_\\}$ for some~$\\\in\Part$. Setting
$$
  \edgprm\alpha\\=\scalar<s_\alpha,s_\\>
\ttex{for $\alpha\in\comp$ and $\\\in\Part$}
,\label(\edgprmdef)
$$
one has $\edgprm\alpha\\\in\{-1,0,1\}$, and given $\alpha$ there is at most
one~$\\$ with $\edgprm\alpha\\\neq0$; this symbol will be used as a signed
variant of the Iverson symbol. Any expression in terms of such~$s_\alpha$ can
be converted to one in terms of Schur functions using
$$
%\postdisplaypenalty=2000
  s_\alpha=\sum_{\\\in\Part}\edgprm\alpha\\s_\\
,\label(\normaleq)
$$
which removes null terms, and replaces the remaining terms~$s_\alpha$ by the
appropriate $\pm{s_\\}$ with~$\\\in\Part$. This process, which is the main
source of alternating sums in this paper, will be called normalisation.

\subsection Diagram boundaries.

The partition~$\\$ with $\edgprm\alpha\\\neq0$, if any, is characterised by
the condition $\dn\\=(\dn\alpha)^+$, where $n$ is so large that
$\alpha\in\N^n$, and $\edgprm\alpha\\$ is the sign of the
permutation~$\sigma\in\Sn$ such that $\dn\\=\sigma(\dn\alpha)$. The given
condition gives rise to an equivalent one when $n$ is increased, and the
permutation involved does not change either, if each $\Sn$ is considered as a
subgroup of~$\Sym_\infty$. Nonetheless, it is convenient to have a description
that does not involve~$n$ at all. If we subtract $n$ from each of the
$n$~entries of~$\dn\alpha=(n-1-i+\alpha_i)_{i\in\set{n}}$, then
coefficient~$i$ becomes $\alpha_i-1-i$, and therefore independent of~$n>i$;
moreover this transformation is compatible with the action of~$\Sn$ by
permutation of the coefficients. Now increasing~$n$ allows to associate
with~$\alpha$ a unique infinite sequence of numbers; note however that it will
contain negative entries, and will no longer be finitely supported. We shall
denote this sequence by $\[\alpha,]=(\[\alpha,i])_{i\in\N}$, where
$$\[\alpha,i]=\alpha_i-1-i
 \ttex{for $\alpha\in\comp$ and $i\in\N$}
.\nn
$$
Thus if $\edgprm\alpha\\$ is nonzero, it is the sign of the unique permutation
$\sigma\in\Sym_\infty$, that transforms $\[\alpha,]$ into $\[\\,]$.

Instead of using sequences $\[\alpha,]$, one could understand $\edgprm\alpha\\$
using an alternative way for the group~$\Sym_\infty$ to act on finitely
supported sequences of integers: defining $\sigma\cdot\alpha=
\(\alpha_{\sigma\inv(i)}+i-\sigma\inv(i)\)_{i\in\N}$, one may interpret
$\edgprm\alpha\\$ as the sign of the unique permutation
$\sigma\in\Sym_\infty$, if any, such that $\sigma\cdot\alpha=\\$. However,
considering sequences $\[\alpha,]$ has the advantage that they can be
interpreted directly in a graphical manner, in particular when $\alpha$ is a
partition. Each part~$\\_i$ of a partition~$\\$ corresponds to a row of its
diagram~$\set\\$, and therefore to a unique vertical segment of the bottom
right boundary of~$\set\\$, namely the segment that delimits row~$i$
of~$\set\\$ to the right. In case $\\_i=0$, row~$i$ of~$\set\\$ is empty, but
we can still consider it as delimited to the right by the segment on the
vertical axis that crosses row~$i$. Now, whereas the sequence~$\\$ records the
horizontal coordinates of these vertical segments (taken from top to bottom),
the sequence~$\[\\,]$ records their ``diagonal coordinates'', where a segment
has diagonal coordinate~$d$ if it connects points on diagonals $d$ and~$d+1$.
For instance, for $\\=(7,5,2,2,1)$ the shape of the boundary can be drawn as
$$
  \smallsquares{\rtab1 \edgeseq8,0;111010110001001000000 }
,
$$
and the sequence of the diagonal coordinates of its vertical segments, taken
from top to bottom, is $(6,3,-1,-2,-4,-6,-7,\ldots)=\[\\,]$. For general
compositions~$\alpha$ one could define a ``diagram'' in the same way as for
partitions (although we did not), and while their boundaries are more ragged,
reading off the diagonal coordinates of the vertical segments from top to
bottom still gives the sequence~$\[\alpha,]$.

Let us make more precise some terms that were used above. Both squares and
points are formally elements of~$\N\times\N$, but we shall treat them as
different types of objects. The square~$(i,j)$ has four corners, which are the
points in the set $\{i,i+1\}\times\{j,j+1\}$, so the square has the same
coordinates as (and is formally speaking identified with) its top left corner;
this is in keeping with our convention to index by~$i$ the interval from $i$
to~$i+1$. Here we shall treat the diagonals~$D_d$ (which were defined
in~(\diagonaldef)) as sets of points, which should not be too confusing since
a square ``lies on the same diagonal'' as its top left corner. The boundary
associated to~$\\\in\Part$ consists of a set of points, one on each
diagonal~$D_d$, which are connected by edges between the point on $D_d$ and
the point on~$D_{d+1}$ for each~$d\in\Z$. Using the fact that each diagonal is
totally ordered under the coordinate-wise partial ordering of points, the
point on~$D_d$ can be characterised as the minimal element of the difference
set $D_d\setminus\set\\$, in other words the top-leftmost point on~$D_d$ that
is not the top-left corner of any square of~$\set\\$. The edges of this
boundary are also indexed by numbers~$d\in\Z$ (their diagonal coordinate), as
described above. Thus the vertical edge associated to the part~$\\_i$ of~$\\$
is the one between the points $(i,\\_i)$ and~$(i+1,\\_i)$, which has as
diagonal coordinate the index $\\_i-(i+1)=\[\\,i]$ of the diagonal containing
the latter point. In the example above the vertical segments are
\begingroup
\def\seg
#1#2{\count0=#1\advance\count0by1\hbox{$(#1,#2)$--$(\number\count0,#2)$, }}%
\seg07 \seg15 \seg22 \seg32 \seg41 \seg50 \seg60\dots
\endgroup

We did not use the horizontal segments of the boundary associated
to~$\\\in\Part$. However, since the sequence~$\[\\,]$ is strictly decreasing,
it is entirely determined by the \emph{set}
$\setof\[\\,i]:i\in\N\endset\ssubset\Z$, and the set of diagonal coordinates
of the horizontal segments of the boundary forms the complement of that set
in~$\Z$, which can be written as $\setof-1-\[\\\tr,j]:j\in\N\endset$. It is
often useful to think of~$\\$ as represented by the characteristic function of
the set $\setof\[\\,i]:i\in\N\endset$, i.e., the function $\Z\to\{0,1\}$ that
takes the value~$1$ at~$d$ if $d$ is the diagonal coordinate of a vertical
segment of the boundary, and the value~$0$ if $d$ is the diagonal coordinate
of a horizontal segment. In turn that function can be thought of as a doubly
infinite sequence of bits (the edge sequence of~$\\$), which describes the
form of the boundary, traversed from bottom left to top right. The relations
$\mu\lev\\$, $\mu\leh\\$ and $\mu\rlt{k}\\$ can be easily understood in terms
of transpositions of bits in edge sequences; notably $\mu\rlt{k}\\$ means that
the edge sequence of~$\\$ is obtained from that of~$\mu$ by transposing a
single bit~$1$ with a bit~$0$ that is $k$ places to its right. This is
graphically obvious if we draw the diagram of the $k$ by superimposing the
boundaries associated to the partitions involved:
$$
  (7,5,2,2,1)\rlt{10}(7,6,6,3,3,2)
 \:\qquad
  \hbox{$\smallsquares{\rtab1 
           \edgeseq9,0;1111010110001001000
	   \edgeseq6,0   ;00101100011 }
        $}
;
$$
the transposition involves the underlined bits in
$\cdots111\underline1010110001\underline001000\cdots$. Moreover
$\height(\\/\mu)$ is the sum of the bits that are being ``jumped over'' during
the transposition, which is $4$ in the example.

\proclaim Proposition. \ribboncharprop
For $k>0$ and $\mu,\\\in\Part$, the relation $\mu\rlt{k}\\$ holds if and only
if there exist indices $i_0,i_1\in\Z$ such that
$\setof\[\mu,i]:i\in\N-\{i_0\}\endset=\setof\[\\,i]:i\in\N-\{i_1\}\endset$ and
$\[\mu,i_0]+k=\[\\,i_1]$. In this case one has moreover
$\height(\\/\mu)=i_0-i_1$.

\proof
By definition, $\mu\rlt{k}\\$ means that $k$ consecutive points on the
boundary associated to~$\mu$ must be moved one place down their diagonal to
reach the boundary associated to~$\\$. This means that only the two segments
that link a point of this subset of~$k$ points and a point of its complement
change orientation, and their diagonal coordinates differ by~$k$. In terms of
the set $\setof\[\mu,i]:i\in\N\endset$, this means that one element
$\[\mu,i_0]$ is replaced by $\[\mu,i_0]+k=\[\\,i_1]$, and since
$i_0$~and~$i_1$ are the row numbers of the vertical boundary segments with
diagonal coordinates $\[\mu,i_0]$ and $\[\\,i_1]$, respectively, one has
$\height(\\/\mu)=i_0-i_1$.
\QED

% Local IspellDict: default
% -*- mode: tex; tex-main-file: "root.tex"; -*-

\subsecno=-1
\newsection The Pieri and Murnaghan-Nakayama rules.
\labelsec\Pierisec

In this section we shall compute the product of Schur polynomials by
respectively elementary, power-sum and complete symmetric functions,
expressing the result in the basis of Schur functions. These lead to nice
combinatorial formulae known as the Pieri and Murnaghan-Nakayama rules, in
which one encounters the notions of vertical and horizontal strips, and of
ribbons.

Before starting our computations, we consider the validity of summations
involving the values~$s_\alpha$ for all $\alpha\in\comp_n$ at once. Although
each of these is either zero or up to a sign equal to one of the finitely many
Schur functions $s_\\$ with $\\\in\Part_n$, this does not prove that only
finitely many $s_\alpha$ are nonzero, and such summations therefore well
defined. The following proposition shows that this is nevertheless the case.

\proclaim Proposition. \alphafiniteprop
If $s_\alpha\neq0$ for some $\alpha\in\comp_n$, then $\alpha\in\N^n$, where
$\N^n$ is considered as subset of\/~$\comp$.

\proof
If $s_\alpha\neq0$, then there exists $\\\in\Part_n\subset\N^n$ with
$\edgprm\alpha\\\neq0$. Now the sequence $\[\alpha,]$ is a permutation of
$\[\\,]$, and therefore in particular it contains all values $-1-i$ for
$i\geq{n}$ as entries. If $\alpha$ had any part $\alpha_i>0$ with $i\geq{n}$,
then for the largest such~$i$, the value $-1-i$ would not occur in the
sequence $\[\alpha,]$, which would give a contradiction. Therefore one must
have $\alpha\in\N^n$.
\QED

The starting point of our approach is a very simple formula describing the
multiplication between a Schur function and an arbitrary symmetric function.

\proclaim Lemma. \Schurmul
Let $f\in\Lambda$, and write $f[X_\N]$ as a sum of monomials
$f[X_\N]=\sum_{\alpha\in\comp}c_\alpha{X^\alpha}$, then
$$
  fs_\beta=\sum_{\alpha\in\comp}c_\alpha s_{\alpha+\beta}
 \ttex{for all~$\beta\in\comp$.}
$$

\proof
First, we must assure that the summation is effectively finite (or stated more
fancily: that it converges in~$\Lambda$ for the discrete topology), so that
the right hand side is well defined. For this it suffices to observe that in
order to have $c_\alpha\neq0$ one must have $|\alpha|\leq\deg{f}$, and that in
order to have $s_{\alpha+\beta}\neq0$, one must have
$\alpha+\beta\in\N^{|\alpha|+|\beta|}$ by proposition~\alphafiniteprop, and
therefore certainly $\alpha\in\N^{|\alpha|+|\beta|}$, which leaves only
finitely many possibilities for~$\alpha$. With this point settled, it will
suffice to prove that the image
$f[\Xn]s_\beta[\Xn]=\sum_{\alpha\in\N^n}{c_\alpha}s_{\alpha+\beta}[\Xn]$ of
the identity in~$\spoln$ holds for all~$n\in\N$. We may assume $\beta\in\N^n$;
it is then obvious that
$f[\Xn]X^\dn\beta=\sum_{\alpha\in\N^n}c_\alpha{X^\dn{\alpha+\beta}}$. To this
equation we apply every $\sigma\in\Sn$, acting by permutation of the
indeterminates, and take the alternating sum by~$\eps(\sigma)$ of the
resulting equations. Since $f[\Xn]$ is invariant under such permutations, we
obtain the equation $f[\Xn]a_\dn\beta[\Xn]
=\sum_{\alpha\in\N^n}{c_\alpha}a_\dn{\alpha+\beta}[\Xn]$ in~$\apoln$. Now
division by~$\Delta_n$ gives the desired identity.
\QED

This lemma can be paraphrased as follows: when multiplying symmetric functions
by Schur functions, one may behave as if $X^\alpha{s_\beta}$ were equal
to~$s_{\alpha+\beta}$ (although of course it is not). This lemma will be our
main tool for this paper. It allows formulae for products of symmetric
functions by Schur functions to be obtained very easily, but the resulting
expression on the right hand side is not normalised, even
when~$\beta\in\Part$. After normalisation one gets an alternating summation
with possibly many internal cancellations; the remaining task is to perform
such cancellations explicitly, and describe the surviving terms in a
combinatorial way. This often involves a certain amount of arbitrariness: in
general an alternating sum neither uniquely determines the set of terms that
will survive the cancellation, nor the way in which the other terms are
paired~up so as to cancel (in so far as they are not already null initially).
In our applications however, a set of candidates for survival that are already
normalised will usually present itself, and the arbitrariness is limited to
the way in which the remaining terms cancel.

Here is a concrete example of a computation using lemma~\Schurmul, namely that
of the decomposition of the product $m_{(3,1)}s_{(2,1)}$ into Schur functions.
We traverse the $\Sym_\infty$-orbit of $(3,1)$ in decreasing lexicographic
order ($(3,1)$, $(3,0,1)$, $(3,0,0,1)$, \dots, $(1,3)$, $(1,0,3)$, \dots,
$(0,3,1)$, \dots, $(0,1,3)$, \dots, $(0,0,3,1)$, \dots), but omit those
elements~$\alpha$ for which $\[{((2,1)+\alpha)},]$ has repeated entries. After
normalisation there are no cancellations, but in two cases a pair of terms is
combined into one:
$$
\def\breakline{\cr &\qquad}% this must be a macro to be used conditionally
\eqalign
{ \m31 \s21 &= \sum_{\alpha\in\Sym_\infty\cdot(3,1)}s_{(2,1)+\alpha}
\cr &=\s52 +\s511
     +\s313 +\s31003 +\s241 \breakline
     +\s2203 +\s22003 +\s2131 +\s211003 +\s210031
\cr &=\s52 +\s511 -\s322 +\s31111 -\s331
\ifblownup\breakline\fi
      -2\s2221 +\s22111 +2\s211111
.\cr}
$$

In the remainder of this section we shall apply the lemma three times, in each
case for a symmetric function~$f$ all of whose nonzero coefficients~$c_\alpha$
are equal to~$1$. The first and easiest application will be for $f=e_d$, where
the monomials~$X^\alpha$ with nonzero coefficients are those for binary
compositions~$\alpha$ of~$d$.

\proclaim Proposition. \Pierielemprop
For $d\in\N$ and $\mu\in\Part_n$ one has
$e_ds_\mu=\sum_{\\\in\Part_{n+d}}\Kr{\mu\lev\\}s_\\$.

\proof
Applying lemma~\Schurmul\ gives $e_ds_\mu
=\sum_{\alpha\in\bcomp_d}s_{\mu+\alpha}$.
Now $\mu+\alpha$ can only fail to be a partition if for some index~$i$ one
has $\mu_i=\mu_{i+1}$ while $(\alpha_i,\alpha_{i+1})=(0,1)$; in that case one
has $s_{\mu+\alpha}=0$ since $\[(\mu+\alpha),i]=\[(\mu+\alpha),i+1]$. The
remaining terms are already normalised, and the corresponding values of
$\mu+\alpha$ run through the set $\setof\\\in\Part_{n+d}:\mu\lev\\\endset$
by the definition (\levhdef) of the relation~`$\lev$'.
\QED

As a first consequence of this proposition one finds for $\mu=(0)$ that
$s_{1\repex^d}=e_d$, since $s_{(0)}=1$ and
$\setof\\\in\Part_d:0\lev\\\endset=\Part_d\thru\bcomp=\{1\repex^d\}$; this
fact can also be easily deduced directly from the definition of Schur
functions. Another instance of the proposition gives
$e_1s_\mu=\sum_{\\\in\Part_{n+1}}\Kr{\mu\subset\\}s_\\$ for~$\mu\in\Part_n$.

In this first application of lemma~\Schurmul, it sufficed to discard the null
terms in the resulting summation. Our second application, which will be for
$f=p_k$, is slightly more complicated: it still suffices to discard the null
terms (there are no cancellations of pairs of terms), but the remaining terms
are not all normalised, and in normalising them some signs do remain
(for~$k\geq2$). This will in fact be the only case in our paper where we
obtain a result that is not a positive sum of Schur functions. The monomials
occurring in~$p_k$ are just the powers $X_i^k$, which are the
monomials~$X^\alpha$ where $\alpha$ is a permutation of the composition~$(k)$,
in other words a composition of~$k$ with just one nonzero part $\alpha_i=k$
for some~$i\in\N$.

\proclaim Proposition. \MurNakprop
For $k>0$ and $\mu\in\Part_n$ one has $p_ks_\mu=
\sum_{\\\in\Part_{n+k}}\Kr{\mu\rlt{k}\\}(-1)^{\height(\\/\mu)}s_\\$.

\proof
Applying lemma~\Schurmul\ gives $p_ks_\mu =\sum_{i\in\N}s_{\mu+\alpha^{(i)}}$
where $\alpha^{(i)}=(\Kr{i=j}k)_{j\in\N}$ is the composition of~$k$ with a
unique nonzero part at position~$i$. If for any $j<i$ one has
$\[\mu,i]+k=\[\mu,j]$, then the corresponding term $s_{\mu+\alpha^{(i)}}$
vanishes, because the sequence $\[(\mu+\alpha^{(i)}),]$ has equal entries at
positions $i$ and~$j$. We discard such terms, and henceforth suppose that no
such~$j$ exists. Then the sequence $\[(\mu+\alpha^{(i)}),]$ can be transformed
into a strictly decreasing sequence by successively transposing the entry
$\[(\mu+\alpha^{(i)}),i]=\[\mu,i]+k$ with its preceding entry some number of
times (possibly zero), and the resulting sequence is of the form $\[\\,]$ for
some $\\\in\Part_{n+k}$. If $j\leq{i}$ is the final position of the term
$\[\mu,i]+k$, so that $\[\mu,i]+k=\[\\,j]$, then the sign of the permutation
that has been applied is $(-1)^{i-j}=\edgprm{\mu+\alpha^{(i)}}\\$. Now by
proposition~\ribboncharprop\ one has $\mu\rlt{k}\\$ and $\height(\\/\mu)=i-j$,
and the proposition follows.
\QED

One obtains for instance, for~$\mu=(0)$, an expression of $p_k$ as an
alternating sum of ``hook shape'' Schur functions:
$p_k=\sum_{i\in\set{k}}(-1)^is_{\nu^{(i)}}$ where $\nu^{(i)}=(k-i,1\repex^i)$
is the partition whose nonzero parts consist of one part~$k-i$ followed by $i$
parts~$1$; a concrete example is
$p_4=s_{(4)}-s_{(3,1)}+s_{(2,1,1)}-s_{(1,1,1,1)}$. 

Using~(\scalarprodid), one sees that this proposition is equivalent
to the statement that for $\\\in\Part_n$, \ $k>0$, and $f\in\Lambda$ of degree
$n-k$, one has
$$\scalar<p_kf,s_\\>
  =\sum_{\mu\in\Part_{n-k}}
  \Kr{\mu\rlt{k}\\}(-1)^{\height(\\/\mu)}\scalar<f,s_\mu>
.\nn
$$
Since it is known that a scalar product of the form
$\scalar<\fold(pk_1..k_l),s_\\>$ can be interpreted as the value of the
symmetric group character parametrised by~$\\$ on the class of permutations
with cycles of lengths $\li(k1..l)$, this gives a recurrence relation for
symmetric group characters. For instance, the value of the
$\Sym_{17}$-character parametrised by~$\\=(7,5,2,2,1)$ on a permutation of
cycle type $(4,3,3,2,2,1,1,1)$ is $\scalar<p_4\,p_3^2\,p_2^2\,p_1^3,s_\\>$,
which can be computed either as $-\scalar<p_3^2\,p_2^2\,p_1^3,s_{(4,4,2,2,1)}>
+\scalar<p_3^2\,p_2^2\,p_1^3,s_{(7,5,1)}>$ using the recurrence relation for
$k=4$, or as $\scalar<p_4\,p_3\,p_2^2\,p_1^3,s_{(7,2,2,2,1)}>
-\scalar<p_4\,p_3\,p_2^2\,p_1^3,s_{(7,5,2)}>$ using the recurrence relation
for $k=3$ (one could also use $k=2$ or~$k=1$, but that gives expressions with
$3$ and $4$ terms, respectively); if the character tables for $\Sym_{13}$ or
$\Sym_{14}$ have been previously determined, this can be evaluated to give the
value $-3+-2=-4-1=-5$. This recurrence relation is known as the
Murnaghan-Nakayama rule, and appears in \ref{Murnaghan} and \ref{Nakayama}.
That name is also given to the combinatorial expression for
$\scalar<\fold(pk_1..k_l),s_\\>$ in terms of chains of partitions related
by~`$\rlt{k}$' for varying~$k$, that can be obtained by expanding the
recurrence relation recursively; however, in that form, which can already be
found in \ref{Littlewood Richardson}, the rule is rather less practical for
computations, due to the large number of terms it often gives (which in
addition depends heavily on the order in which the factors $p_{k_i}$ are
extracted).

Our third application of lemma~\Schurmul\ will be for $f=h_d$; here the
coefficients $c_\alpha$ are nonzero for all compositions $\alpha\in\comp_d$.
Thus the summation obtained from the lemma gives many more terms than in the
previous applications, and we must perform some real cancellations. Yet,
surprisingly, the resulting linear combination of Schur functions will be
quite similar to the one obtained in proposition~\Pierielemprop, notably all
nonzero coefficients are equal to~$1$. In fact the only difference is the
replacement of vertical strips by horizontal strips. That proposition and the
following one are known as the Pieri rules (or Pieri formulae).

\proclaim Proposition. \Piericomplprop
For $d\in\N$ and $\mu\in\Part_n$ one has
$h_ds_\mu=\sum_{\\\in\Part_{n+d}}\Kr{\mu\leh\\}s_\\$.

\proof
Application of lemma~\Schurmul\ gives $h_ds_\mu
=\sum_{\alpha\in\comp_d}s_{\mu+\alpha}$. In this summation we must cancel all
the terms with $\mu\nleh\mu+\alpha$, in other words those terms with
$(\mu+\alpha)_{i+1}>\mu_i$ for some index~$i\in\N$. That inequality is
equivalent to $\[(\mu+\alpha),i+1]\geq\[\mu,i]$, which means that
$\[(\mu+\alpha),i+1]$ lies outside the ``safe interval''
$I_{i+1}=\setof{j\in\Z}:\[\mu,i+1]\leq{j}<\[\mu,i]\endset$. Setting in
addition to this $I_0=\setof{j\in\Z}:\[\mu,0]\leq{j}\endset$, the
intervals~$I_i$ form a partition of~$\Z$, so $\alpha$ determines a unique map
$f\:\N\to\N$ such that $\[(\mu+\alpha),i]\in{I_{f(i)}}$ for all~$i$. Moreover,
one clearly has $f(i)\leq{i}$ for all~$i$, so if $f$ is injective, it must be
the identity, which is equivalent to $\mu\leh\mu+\alpha$. We must then cancel
the terms for those~$\alpha$ for which $f$ is not injective. Indeed these
terms can be cancelled, since if $f(i_0)=f(i_1)$ for $i_0\neq{i_1}$, then the
transposition of entries $i_0$~and~$i_1$ of the sequence $\[(\mu+\alpha),]$
gives a sequence of the form $\[(\mu+\alpha'),]$, with $\alpha'\in\comp_d$ and
$s_{\mu+\alpha}+s_{\mu+\alpha'}=0$; either these are two distinct terms that
cancel each other, or if $\alpha=\alpha'$ they are twice the same term, which
is then null already. It remains to show that cancelling terms can be paired
up effectively, in other words that a concrete pair of indices $(i_0,i_1)$ can
be chosen for each $\alpha\in\comp_d$ with $\mu\nleh\mu+\alpha$, in such a way
that for the resulting $\alpha'$ the same pair $(i_0,i_1)$ will be chosen.
This is easy: since the function $f$ is the same for any $\alpha'$ that would
cancel with~$\alpha$, any choice that depends only on~$f$ will do; for
instance one can take the minimal possible pair~$(i_0,i_1)$ in lexicographic
order.
\QED

Similarly to what we saw for proposition~\Pierielemprop, a first consequence
of proposition~\Piericomplprop\ is that $s_{(d)}=h_d$, which is not quite as
easy to obtain directly from definition~\Schurfuncdef\ as $s_{1\repex^d}=e_d$
is. As an illustration we choose a rather modest example, since the number of
terms in the summation tends to grow rather rapidly: one has $ h_4\s4221
=\s9221 +\s8321 +\s8222 +\s82221 +\s7421 +\s7322 +\s73211 +\s72221 +\s6422
+\s64211 +\s63221 +\s54221 $. Although a dozen terms remain, the number of
nonzero terms produced by application of lemma~\Schurmul\ before normalisation
is considerably larger, namely~$178$. To illustrate the cancellation, consider
the cancelled term $\s42422 $ that arises for $\alpha=(0,0,2,1,2)$. Putting
$\mu=(4,2,2,1)$, and comparing $\[(\mu+\alpha),] =\[{(4,2,4,2,2)},]
=(3,0,1,-2,-3,-6,-7\ldots)$ with $\[\mu,]=(3,0,-1,-3,-5,-6,-7\ldots)$, we see
that the entries of the former sequence that lie outside their respective save
intervals are $\[(\mu+\alpha),2]=1\notin{I_2}=\sing{-1}$
and~$\[(\mu+\alpha),4]=-3\notin{I_4}=\sing{-5,-4}$. Since in fact
$1\in{I_1}=\sing{0,1,2}$ and $-3\in{I_3}=\sing{-3,-2}$, the values
$f(i)\neq{i}$ are $f(2)=1$ and $f(4)=3$. Therefore a term~$s_{\mu+\alpha'}$
cancelling $s_{\mu+\alpha}$ can be obtained either by transposing the entries
$0$ and~$1$ at indices $1$~and~$2$ of $\[(\mu+\alpha),]$ or by transposing the
entries $-2$ and~$-3$ at indices $3$~and~$4$. The result can be written as
$\[{(4,3,3,2,2)},]$ respectively as~$\[{(4,2,4,1,3)},]$, and indeed one has
$\s43322 =\s42413 =-\s42422 $, while the occurring compositions can be written
as $\mu+\alpha'$ for $\alpha'=(0,1,1,1,2),(0,0,2,0,3)$, respectively, so that
both $\s43322 $ and~$\s42413 $ occur in the summation. The rule of choosing
the lexicographically first pair of indices would in fact make the term
$\s42422 $ cancel against $\s43322 $.

As with proposition \MurNakprop, propositions \Pierielemprop\
and~\Piericomplprop\ lead to recurrence relations for certain scalar products.
In particular we obtain such relations for the Kostka numbers:
$$
\eqalignno
{ K'_{\\/\mu,\alpha}= \scalar<s_\mu e_\alpha,s_\\>
 & = \sum_{\mu'\in\Part_{|\mu|+\alpha_0}}\Kr{\mu\lev\mu'}
     K'_{\\/\mu',(\alpha_1,\alpha_2,\ldots)}
,&\nn \cr
  K_{\\/\mu,\alpha}= \scalar<s_\mu h_\alpha,s_\\>
 & = \sum_{\mu'\in\Part_{|\mu|+\alpha_0}}\Kr{\mu\leh\mu'}
     K_{\\/\mu',(\alpha_1,\alpha_2,\ldots)}
,&\nn \cr}
$$
for all $\\,\mu\in\Part$ and $\alpha\in\comp$. One can recursively apply these
recurrence relations; using the basic case
$K'_{\\/\\,\alpha}=K'_{\\/\\,\alpha}=\Kr{\alpha=(0)}$ and comparing with
definition~\SSTdef, one finds the
promised combinatorial description of Kostka numbers.

\proclaim Corollary. \Kostkacor
For $\\,\mu\in\Part$ and $\alpha\in\comp$, the Kostka number
$K_{\\/\mu,\alpha}$ equals the number~$\Card{\SSTw{\\/\mu}\alpha}$ of
semistandard tableaux of shape $\\/\mu$ and weight~$\alpha$, and
$K'_{\\/\mu,\alpha}$ is the number of transposed semistandard tableaux of
shape $\\/\mu$ and weight~$\alpha$. As a consequence the two are related by
$K'_{\\/\mu,\alpha}=K_{\\\tr/\mu\tr,\alpha}$.
\QED

This description shows that one has $\Card{\SSTw{\\/\mu}\alpha}
=\Card{\SSTw{\\/\mu}{\sigma(\alpha)}}$ for any $\sigma\in\Sym_\infty$, and it
turns equation~(\heKostkaeq) into a combinatorial expression of $h_\alpha$ and
$e_\alpha$ in terms of Schur functions. In fact, we may consider
$h_\alpha$ and $e_\alpha$ as generating polynomials in~$\Lambda$ of
semistandard (respectively transposed semistandard) Young tableaux of
weight~$\alpha$ by shape, in the basis of Schur functions:
$$
\eqalignno
{h_\alpha
&=\sum_{\\\in\Part_{|\alpha|}}\sum_{T\in\SSTx\\}
  \Kr{\weight{T}=\alpha}s_\\
,&\nn\cr
 e_\alpha
&=\sum_{\\\in\Part_{|\alpha|}}\sum_{T\in\SSTx{\\\tr}}
  \Kr{\weight{T}=\alpha}s_\\
.&\nn\cr
}
$$
Besides the combinatorial description of corollary~\Kostkacor, one also has
the following expressions for $K_{\\/\mu,\beta}$ and $K'_{\\/\mu,\beta}$ that
can be obtained directly from lemma~\Schurmul, using (\hbetaeq) and~(\ebetaeq):
$$
\eqalignno
{
  K_{\\/\mu,\beta}
 &=\sum_{\alpha\in\comp}\Card{\mat\alpha\beta}\scalar<s_{\mu+\alpha},s_\\>
  =\sum_{M\in\Mat}\Kr{\csum{M}=\beta}\edgprm{\mu+\rsum{M}}\\
,&\label(\Kostkaalteq)
\cr
  K'_{\\/\mu,\beta}
 &=\sum_{\alpha\in\comp}\Card{\bmat\alpha\beta}\scalar<s_{\mu+\alpha},s_\\>
  =\sum_{M\in\bMat}\Kr{\csum{M}=\beta}\edgprm{\mu+\rsum{M}}\\
.&\label(\Kostkatralteq)
\cr}
$$
These expressions produce the Kostka numbers by an alternating enumeration of
matrices rather than by a direct enumeration of tableaux, and are therefore
rather inefficient as a means of computing those numbers. However it will
prove useful to have, as a complement to corollary~\Kostkacor, these
descriptions that do not involve semistandard tableaux, or horizontal or
vertical strips.

% Local IspellDict: default
% -*- mode: tex; tex-main-file: "root.tex"; -*-

\newsection The Cauchy, Jacobi-Trudi and Von N{\"a}gelsbach-Kostka identities.
\labelsec\Cauchysec

We have expressed products of elementary, of power sum, and of complete
symmetric functions as $\Z$-linear combinations of Schur functions. However we
have not yet expressed Schur functions themselves in terms of anything else,
not even in terms of monomials, other than via the polynomial division used in
their definition. In this section we shall provide such expressions of (skew)
Schur functions, in terms of minimal symmetric functions (which implies an
expression in terms of monomials), in terms of products $h_\alpha$ of complete
symmetric functions, and in terms of products $e_\alpha$ of elementary
symmetric functions. The latter two expressions will eventually be obtained
via another application of lemma~\Schurmul, but the expression in terms of
minimal symmetric functions will be obtained in an indirect manner, by
establishing that the dual basis of $\setof{m_\\}:\\\in\Part\endset$ with
respect to the scalar product on~$\Lambda$ is in fact the set
$\setof{h_\\}:\\\in\Part\endset$ (with the obvious correspondence between
their elements, in other words $\scalar<m_\\,h_\mu>=\Kr{\\=\mu}$). Another way
of formulating that property is that scalar products $\scalar<f,h_\\>$ can be
interpreted as the coefficient of $m_\\$ in the expression of~$f$ as linear
combination of minimal symmetric functions, and hence that
$\scalar<f,h_\alpha>$ can be interpreted as the coefficient of $X^\alpha$ in
the expansion of $f[X_\N]$ into monomials. This will imply in particular that
$K_{\\/\mu,\alpha}$ can be interpreted as the coefficient of $X^\alpha$
in~$s_{\\/\mu}[X_\N]$.

Let us first show that $\setof{h_\\}:\\\in\Part\endset$ is a $\Z$-basis
of~$\Lambda$ in the first place. Since its elements are precisely the
monomials that can be formed from the set $\setof{h_i}:i>0\endset$, this
amounts to showing that those $h_i$ are algebraically independent generators
of~$\Lambda$, which can then be viewed as a polynomial ring
$\Z[h_1,h_2,\ldots]$. One has $h_\mu=\sum_{\\\in\Part_n}K_{\\,\mu}s_\\$
for $\mu\in\Part_n$, while $\setof{s_\\}:\\\in\Part\endset$ is known to be a
$\Z$-basis of~$\Lambda$, so we can obtain our goal by showing that the matrix
$(K_{\\,\mu})_{\\,\mu\in\Part_n}$ of Kostka numbers is invertible for
every~$n$. That will follow from the fact that this matrix is unitriangular
for the dominance ordering on~$\Part_n$.

\proclaim Proposition.
For $\\\in\Part_n$ and $\alpha\in\comp_n$, $K_{\\,\alpha}\neq0$ implies
$\alpha^+\leq\\$; moreover $K_{\\,\alpha}=1$ if $\alpha^+=\\$.

\proof
Since $K_{\\,\alpha}=K_{\\,\alpha^+}$ we may assume $\alpha=\alpha^+\in\Part$.
We use the description of $K_{\\,\alpha}$ in corollary~\Kostkacor. For any
Young tableau $(\\^{(i)})_{i\in\N}\in\SSTx\\$, the diagram $\set{\\^{(k)}}$ is
contained in the topmost $k$ rows of~$\set\\$, for all $k\in\N$. The existence
of $(\\^{(i)})_{i\in\N}\in\SSTw\\\alpha$ therefore implies
$\sum_{i\in\set{k}}\alpha_i=|\\^{(k)}|\leq\sum_{i\in\set{k}}\\_i$. If
$\alpha=\\$, we see that $\set{\\^{(k)}}$ must coincide with the topmost $k$
rows of~$\set\\$, whence $\Card\SSTw\\\\=1$.
\QED

A combinatorial proof of this proposition for general~$\alpha\in\comp$,
without using $K_{\\,\alpha}=K_{\\,\alpha^+}$, can also be given, but is
slightly more awkward to formulate (one basically argues that the number of
squares contributed by any column to the first $k$ rows of~$\set\\$ is at
least as large as the contribution, in any $T\in\SSTw\\\alpha$, of the $k$
largest parts of~$\alpha$ to that column). But we shall only use the case
$\alpha\in\Part$.

The proposition implies that the matrix $(K_{\\,\mu})_{\\,\mu\in\Part_n}$ is
upper unitriangular, and hence invertible, if $\Part_n$ is ordered in such a
way that $\\$ comes before~$\mu$ whenever $\\\geq\mu$, for instance by
decreasing lexicographic order. Since, by corollary~\Kostkacor,
$(K'_{\\,\mu})_{\\,\mu\in\Part_n}$ is obtained from that matrix by a
permutation of its rows (the permutation being given by the map
$\\\mapsto\\\tr\: \Part_n\to\Part_n$), we may conclude from
$e_\mu=\sum_{\\\in\Part_n}K'_{\\,\mu}s_\\$ that $\Lambda$ is also a polynomial
ring $\Z[e_1,e_2,e_3\ldots]$, a fact that is sometimes called the fundamental
theorem of symmetric functions, and that is usually proved without using Schur
functions.

To prove that the basis $\setof{h_\\}:\\\in\Part\endset$ is dual to
$\setof{m_\\}:\\\in\Part\endset$, one may establish, for $f$ running through
some basis of~$\Lambda$, that $f=\sum_{\\\in\Part}\scalar<f,h_\\>m_\\$ or
equivalently that $f[X_\N]=\sum_{\alpha\in\comp}\scalar<f,h_\alpha>X^\alpha$.
It is not practical to use the basis of Schur functions for this, since in
that case those expressions are precisely what we would like to conclude from
duality. However for $f=e_\beta$ or for $f=h_\beta$ we already know their
expressions in terms of monomials, which are given in equations (\ebetaeq) and
(\hbetaeq), while the scalar products that should coincide with the
coefficients in those equations are easily computed. Indeed using
equations~(\heKostkaeq) to develop these elements in the orthonormal basis of
Schur functions one finds
$\scalar<e_\beta,h_\alpha>=\sum_{\\\in\Part}K'_{\\,\beta}K_{\\,\alpha}$ and
$\scalar<h_\beta,h_\alpha>=\sum_{\\\in\Part}K_{\\,\beta}K_{\\,\alpha}$. Thus
we need to show the following.

\proclaim Proposition. \Cauchyprop
For every $\alpha,\beta\in\comp$ one has
$$ \eqalignno
{ \Card{\mat\alpha\beta}&=\sum_{\\\in\Part}K_{\\,\alpha}K_{\\,\beta}
& \nn \cr
\noalign{\hbox{and}}
  \Card{\bmat\alpha\beta}&=\sum_{\\\in\Part}K_{\\,\alpha}K'_{\\,\beta}
.& \nn \cr
}
$$

\proof
This is the one place where we import a result from elsewhere. The proposition
expresses the enumerative consequences of the bijections between integer
matrices and pairs of semistandard tableaux of equal shape defined
in~\ref{Knuth PJM}, the first of which is known as the
\caps{RSK}-correspondence.
\QED

\proclaim Corollary. \dualitycor
The $\Z$-bases $\setof{h_\\}:\\\in\Part\endset$ and
$\setof{m_\\}:\\\in\Part\endset$ of~$\Lambda$ are dual to each other with
respect to $\scalar<\cdotp,\,\cdotp\!>$: for every $f\in\Lambda$ one has
$f=\sum_{\\\in\Part}\scalar<h_\\,f>m_\\$, and consequently
$f[X_\N]=\sum_{\alpha\in\comp}\scalar<h_\alpha,f>X^\alpha$.

\proof
Since we have shown that $\setof{h_\\}:\\\in\Part\endset$ and
$\setof{e_\\}:\\\in\Part\endset$ are $\Z$-bases of~$\Lambda$, either of
the special cases $f=h_\beta$ or $f=e_\beta$ that are given by
proposition~\Cauchyprop\ suffice to prove the general case~$f\in\Lambda$.
\QED

Applying this corollary to $f=s_{\\/\mu}$ gives the expressions that we were
after, of (skew) Schur functions in terms of minimal symmetric function and in
terms of monomials:
$$
\eqalignno
{s_{\\/\mu}
&=\sum_{\nu\in\Part}K_{\\/\mu,\nu}m_\nu
,&\nn
\cr
 s_{\\/\mu}[X_\N]
&=\sum_{\alpha\in\comp}K_{\\/\mu,\alpha}X^\alpha
 =\sum_{T\in\SSTx{\\/\mu}}X^{\weight{T}}
.&\label(\Schurgenf)
\cr}
$$
For~$\mu=(0)$, the final expression describes $s_\\[X_\N]$ as the generating
series in~$\Z[[X_\N]]$ of semistandard Young tableaux of shape~$\\$ by weight.
This is often taken as the definition of Schur functions (for instance
in~\ref{Stanley EC2}), but in our setting it is a nontrivial statement that
provides an alternative combinatorial interpretation of semistandard Young
tableaux $T\in\SSTw\\\alpha$: rather than corresponding to constituent Schur
functions~$s_\\$ of $h_\alpha$, they now correspond to constituent
monomials~$X^\alpha$ of~$s_\\[X_\N]$. One may ask whether this result can be
proved in a bijective manner. The way in which we obtained it does not qualify
as a bijective proof: although the main ingredients (lemma~\Schurmul\ and
proposition~\Cauchyprop) are based on explicit correspondences, we did need
some (linear) algebraic reasoning. The most direct possible bijective proof
would be to give, for
$f=\sum_{\alpha\in\N^n}\sum_{T\in\SSTw\\\alpha}X^\alpha\in\Z[\Xn]$, a
correspondence between terms that establishes $f\Delta_n=a_\dn\\[\Xn]$. That
would require a fairly complicated correspondence, but there is a simpler
possibility that nevertheless is more or less equivalent: taking
$f=\sum_{T\in\SSTx\\}X^{\weight{T}}$ and $\beta=(0)$ in lemma~\Schurmul, its
right hand side $\sum_{T\in\SSTx\\}s_{\weight{T}}$ should reduce by
cancellations to~$s_\\$, and one may seek to describe such a cancellation
process explicitly. For instance for $\\=(2,1)$ the right hand side,
restricted to $\weight{T}\in\N^3$ thanks to proposition~\alphafiniteprop,
gives $\s21 +\s201 +\s12 +2\s111 +\s102 +\s021 +\s012 $; here the terms $\s201
$, $\s12 $, and~$\s012 $ are null, and the remaining terms except the initial
one can be cancelled using $0=\s111 +\s102 =\s111 +\s021 $. We shall see
below, when we consider products of Schur functions, that an explicit and
general description of such a cancellation can indeed be given.

Now that we have identified the generating series of semistandard Young
tableaux of a given shape by weight, we may turn the equations of
proposition~\Cauchyprop\ into generating series identities by multiplying them
by~$X^\alpha{Y^\beta}$ and summing over $\alpha,\beta\in\comp$; one obtains
$$\eqalignno
{ \prod_{i,j\in\N}{1\over1-X_iY_j} &=\sum_{\\\in\Part}s_\\[X_\N]s_\\[Y_\N]
& \label(\Cauchyid)
\cr
\noalign{\hbox{and}}
  \prod_{i,j\in\N}(1+X_iY_j) &=\sum_{\\\in\Part}s_\\[X_\N]s_{\\\tr}[Y_\N]
.& \label(\dualCauchyid)
\cr}
$$
The first of these is called the Cauchy identity (although it is difficult to
justify the attribution; for details see \ref{Stanley EC2, p.~397}),
and the second the dual Cauchy identity. From (\hbetaseries) and~(\Cauchyid)
one deduces (after combining terms in the left hand side of the former with a
common value $\beta^+=\\$) that
$$ \sum_{\\\in\Part}h_\\[X_\N]m_\\[Y_\N]
  =\sum_{\\\in\Part}s_\\[X_\N]s_\\[Y_\N]
\label(\hortheq)
$$
in $\Z[[X_\N,Y_\N]]$, which is an alternative formulation of the duality of
corollary~\dualitycor. Similarly (\ebetaseries) and~(\dualCauchyid) give
$$ \sum_{\\\in\Part}e_\\[X_\N]m_\\[Y_\N]
  =\sum_{\\\in\Part}s_{\\\tr}[X_\N]s_\\[Y_\N]
.\label(\eortheq)
$$
This equation does not express a duality like the previous one, but one can
relate the two equations using the $\Z$-linear involution
$\omega\:\Lambda\to\Lambda$ that maps $s_\\\mapsto{s_{\\\tr}}$ for all
$\\\in\Part$. By~(\heKostkaeq), the fact that
$K'_{\\,\alpha}=K_{\\\tr,\alpha}$ implies $\omega(h_\alpha)=e_\alpha$ for all
$\alpha\in\comp$, so $\omega$ coincides with the ring morphism
$\Lambda\to\Lambda$ that sends $h_i\mapsto{e_i}$ for all~$i>0$ (which exists
since $\Lambda=\Z[h_1,h_2,\ldots]$). Then (\eortheq) can be obtained
from~(\hortheq) by applying, to the subring of~$\Z[[X_\N,Y_\N]]$ that is the
the image of $\Lambda\tensor\Lambda$ via $f\tensor{g}\mapsto{f[X_\N]}g[Y_\N]$,
the ring automorphism corresponding to the automorphism $\omega\tensor\Id$
of~$\Lambda\tensor\Lambda$. (Actually that is not quite true: one needs to
take the closure of the subring in the topology of power series rings, and
extend the automorphism to that closure by continuity.)

We now proceed to express skew Schur functions as $\Z$-linear combinations of
elements $h_\alpha$. Since duality of bases is a symmetric notion, the
coefficient of~$h_\nu$ in the expression of any~$f\in\Lambda$ in the basis
$\setof{h_\\}:\\\in\Part\endset$ can be computed as $\scalar<m_\nu,f>$.
Applying this for $f=s_{\\/\mu}$, we see that its coefficient of~$h_\nu$ is
$\scalar<m_\nu,s_{\\/\mu}>=\scalar<s_\mu{m_\nu},s_\\>$. This number can be
determined by applying lemma~\Schurmul\ for~$m_\nu$ (a case that we did not
yet exploit before), which gives as result the expression
$\sum_{\alpha\in\comp}\Kr{\alpha^+=\nu}\edgprm{\mu+\alpha}\\$. This shows that
$$
  s_{\\/\mu} =\sum_{\nu\in\Part}\sum_{\alpha\in\comp}
	      \Kr{\alpha^+=\nu}\edgprm{\mu+\alpha}\\h_\nu
             =\sum_{\alpha\in\comp}\edgprm{\mu+\alpha}\\h_\alpha
,\label(\preJacobiTrudi)
$$
since $h_\nu=h_\alpha$ when $\alpha^+=\nu$. This achieves the expression of
skew Schur functions in terms of complete symmetric functions, but we shall
proceed to make the expression more transparent. The process will be somewhat
different from what we have seen before, as we shall not actually cancel any
pairs of terms, nor even group together terms for values of~$h_\alpha$ that
are equal.

The nonzero terms of the final summation of~(\preJacobiTrudi) can be found by
permuting the terms of~$\[\\,]$ in all possible ways such that each term of
the resulting sequence $\sigma(\[\\,])$ is at least as large as the
corresponding term of $\[\mu,]$ (in formula, $\[\\,\sigma\inv(i)]\geq\[\mu,i]$
for all~$i\in\N$), which condition we shall abbreviate as
$\sigma(\[\\,])\geq\[\mu,]$; then $\alpha=\sigma(\[\\,])-\[\mu,]$ gives a
contribution $\eps(\sigma)h_\alpha$. From the proof of
proposition~\alphafiniteprop\ we see that if $\\\in\N^l$, then one needs to
consider only permutations~$\sigma\in\Sym_l$. For instance for
$\\/\mu=(5,4,2)/(3,1)$, only permutations of the $3$ initial terms of
$\[\\,]=(4,2,-1,-4,\ldots)$ need to be considered, and only the permutations
giving $(4,2,-1)$, $(2,4,-1)$, $(4,-1,2)$, and~$(2,-1,4)$ satisfy
$\sigma(\[\\,])\geq\[\mu,]=(2,-1,-3,\ldots)$; for this case (\preJacobiTrudi)
therefore gives $s_{(5,4,2)/(3,1)}=\h232 -\h052 -\h205 +\h007 =\h322 -2\h52
+h_7$. We can now write (\preJacobiTrudi) as
$$
  s_{\\/\mu} =\sum_{\sigma\in\Sym_l}
              \Kr{\sigma(\[\\,])\geq\[\mu,]}\,\eps(\sigma)
              \,h_{\sigma(\[\\,])-\[\mu,]}
.\nn
$$
With the convention that $h_i=0$ when $i<0$ (which will take care of the
factor $\Kr{\sigma(\[\\,])\geq\[\mu,]}$), we can write this expression as the
determinant of a matrix whose entries are complete symmetric functions:
$$
  s_{\\/\mu} = \det\(h_{\[\\,i]-\[\mu,j]}\)_{i,j\in\set{l}}
 \ttex{when $\\\in\N^l$}
.\label(\JacobiTrudi)
$$
This identity is called the Jacobi-Trudi identity. For the example above, it
gives
$$
  s_{(5,4,2)/(3,1)}
  = \left|\matrix{h_2&h_5&h_7\cr h_0&h_3&h_5\cr 0&h_0&h_2\cr}\right|
  =\h322 -2\h52 +h_7
.
$$

To obtain an expression of skew Schur functions in terms of elementary
symmetric functions, we may apply the automorphism~$\omega$ to
equation~(\preJacobiTrudi), giving
$$
  s_{\\\tr/\mu\tr} =\sum_{\alpha\in\comp}\edgprm{\mu+\alpha}\\e_\alpha
,\label(\preVonNagelsbachKostka)
$$
from which one deduces as above, again with the convention that $e_i=0$ when
$i<0$, that
$$\catcode`\@=11 % \eqalignno with extra column
  \displ@y \tabskip\centering
  \halign to\displaywidth{\hfil$\@lign\displaystyle#$\tabskip\z@skip
    &$\@lign\displaystyle{}#$\hfil
     \tabskip2em
    &\@lign#\hfil
     \tabskip\centering
    &\llap{$\@lign#$}\tabskip\z@skip\crcr
  s_{\\\tr/\mu\tr} &= \det\(e_{\[\\,i]-\[\mu,j]}\)_{i,j\in\set{l}}
  & when $\\\in\N^l$, & \nn \cr
\noalign{\hbox{or equivalently}}
  s_{\\/\mu} &= \det\(e_{\[\\\tr,i]-\[\mu\tr,j]}\)_{i,j\in\set{r}}
  & when $\\_0\leq{r}$. & \nn \cr
}
$$
This final equation is sometimes referred to as the dual Jacobi-Trudi
identity, but we prefer to call it the Von N\"agelsbach-Kostka identity
(following Littlewood, see \ref{Littlewood, p.~87}, but restoring the spelling
that had been lost in translation). We note that Kostka proved the identity in
\ref{Kostka} by an ingenious manipulation of determinants representing
alternants, see \ref{Muir, p.~155}; see also~\ref{Stanley EC2, p.~397}.

% Local IspellDict: default
% -*- mode: tex; tex-main-file: "root.tex"; -*-

\newsection The Gessel-Viennot correspondences.
\labelsec \GVsec

Let us compare equation~(\Kostkaalteq), which expresses the Kostka number
$K_{\\/\mu,\beta}$ as an alternating sum over integral matrices with column
sums given by~$\beta$, with the description of~$K_{\\/\mu,\beta}$ as counting
semistandard tableaux of shape~$\\/\mu$ and weight~$\beta$. The integral
encodings~$M$ of such tableaux (given in definition~\encodingdef) are also
matrices with $\csum{M}=\beta$, and since moreover $\rsum{M}=\\-\mu$, these
matrices contribute $+1$ to the alternating sum in~(\Kostkaalteq). Therefore
these integral encodings give a subset of its terms that completely accounts
for the value of the alternating sum; apparently all remaining terms cancel
out. We can combine the observed identity for all~$\beta$ into a generating
series identity
$$
  \sum_{\beta\in\comp}K_{\\/\mu,\beta}X^\beta
 =\sum_{M\in\Mat}\eps(\mu+\rsum{M},\\)X^{\csum{M}}
 =\sum_{M\in\IE{\\/\mu}}X^{\csum{M}}
,\label(\GVbasiseq)
$$
where the second member is obtained from~(\Kostkaalteq), and the third member
from corollary~\Kostkacor\ by taking integral encodings. In this section we
shall discuss the reduction of the second member to the third member by
cancellations, and a similar question for an identity derived from
equation~(\Kostkatralteq).

We know from~(\Schurgenf) that the first and last members of~(\GVbasiseq)
equal $s_{\\/\mu}[X_\N]$, and the cancellation that we consider is often
associated with the Jacobi-Trudi identity. Indeed, if one takes
equation~(\hbetaeq), slightly rewritten as $h_\alpha[X_\N]
=\sum_{M\in\Mat}\Kr{\rsum{M}=\alpha}X^{\csum{M}}$, and substitutes it
into~(\preJacobiTrudi), then one obtains
$$
  s_{\\/\mu}[X_\N]=\sum_{M\in\Mat}\edgprm{\mu+\rsum{M}}\\X^{\csum{M}}
,\nn
$$
whose second member matches that of~(\GVbasiseq). We note however that we did
not need (\Schurgenf) or~(\preJacobiTrudi) to obtain~(\GVbasiseq).

In fact, the way we obtained~(\GVbasiseq) suggests a way to understand how to
cancel terms in it. Its second member was originally obtained from an
application of lemma~\Schurmul\ for multiplication by~$h_\beta$, while its
third member was obtained by repeated application of
proposition~\Piericomplprop\ for individual multiplications by
factors~$h_{\beta_j}$, each of which uses lemma~\Schurmul\ followed by a
pruning of the summation. It suffices to postpone such intermediate pruning.
Stated differently, for each term of the second member that needs to be
cancelled, one needs to find the factor~$h_{\beta_j}$ responsible for this
fact, and cancel the term correspondingly.

Let $M\in\Mat$ with $\csum{M}=\beta$. Each column $\smash{M\tr_j}$ represents
the exponent of a term in the factor~$h_{\beta_j}$ of~$h_\beta$, so that the
entire matrix represents the monomial~$X^{\rsum{M}}$ in the expansion
of~$h_\beta$. If we have $M\in\IE{\\/\mu}$, which means that $M$ is the
integral encoding of some semistandard tableau~$T$ of shape~$\\/\mu$, then the
sequence of partitions $(\\^{(l)})_{l\in\N}$ that defines~$T$ can be found by
starting with~$\mu$ and successively adding the compositions $M\tr_j$ for
increasing values of~$j$: $\\^{(l)}=\mu+\sum_{j<l}M\tr_j$. Therefore the
condition that $M\in\IE{\\/\mu}$ amounts to the requirement that
$\mu+\sum_{j<l}M\tr_j\leh\mu+\sum_{j\leq{l}}M\tr_j$ holds for all~$l\in\N$,
and that moreover one has $\mu+\rsum{M}=\\$.

Now suppose that $M\notin\IE{\\/\mu}$, so that $M$ defines a term of the
second member of~(\GVbasiseq) that has to be cancelled. Should only the final
part $\mu+\rsum{M}=\\$ of the condition $M\in\IE{\\/\mu}$ fail, then
$\mu+\rsum{M}$ is a partition distinct from~$\\$ (recall that
$\alpha\leh\alpha'$ implies that $\alpha,\alpha'\in\Part$), which means that
$\edgprm{\mu+\rsum{M}}\\=0$, and the term defined by~$M$ is null (this case
only arises because~(\GVbasiseq) just deals with the coefficient
$\scalar<h_\beta,s_{\\/\mu}>=\scalar<s_\mu{h_\beta},s_\\>$ rather than with
all of $s_\mu{h_\beta}$).

Putting that case aside, there must be an index~$l$ for which
$\mu+\sum_{j<l}M\tr_j\nleh\mu+\sum_{j\leq{l}}M\tr_j$, and we choose the
minimal such~$l$; this means that in the computation where $s_\mu$ is
successively multiplied by the factors~$h_{\beta_j}$, the term determined
by~$M$ is cancelled after the multiplication by~$h_{\beta_l}$. At that point
the term is $s_\alpha$ with $\alpha=\mu+\sum_{j\leq{l}}M\tr_j$, and it is
cancelled against a term $s_{\alpha'}$ defined similarly, but possibly with a
different composition instead of the final term~$M\tr_l$ of the summation (as
usual we do not exclude the possibility $\alpha'=\alpha$, if already
$s_\alpha=0$). From the proof of proposition~\Piericomplprop\ we recall that
$\alpha'$ is such that the sequence $\[\alpha',]$ is obtained
from~$\[\alpha,]$ by interchanging its entries at appropriately chosen
distinct indices $i_0$~and~$i_1$. Those indices were chosen in a manner that
ensures that the value that should replace~$M\tr_l$ in the expression
for~$\alpha'$ is a composition, i.e., it has no negative coefficients. That
composition is in fact obtained from~$M\tr_l$ by a transfer between the entries
$M_{i_0,l}$~and~$M_{i_1,l}$, by which we mean the replacement of those entries
by some pair of natural numbers with the same sum.

We have considered $M$ truncated to its columns $M\tr_j$ with $j\leq{l}$, and
found another such partial matrix such that their contributions to
$s_\mu\fold(h\beta_0..\beta_l)$ cancel. It is clear algebraically that the
total contribution to $s_\mu{h_\beta}$ from matrices that extend the first
partial matrix cancels the total contribution from matrices that extend the
other partial matrix. However, to describe the cancellation in~(\GVbasiseq)
combinatorially, we have to match up the term for~$M$ with a specific term
cancelling it, taking into account also the columns $M\tr_j$ with~$j>l$. The
obvious way to do this is to exchange, in addition to the transfer between
$M_{i_0,l}$~and~$M_{i_1,l}$, the entire remainders of rows $i_0$~and~$i_1$ to
the right of those entries; in this way it is guaranteed that the values of
$\[(\mu+\rsum{M}),i_0]$ and $\[(\mu+\rsum{M}),i_1]$ are interchanged between a
pair of matching matrices.

Here is a simple example to illustrate the construction of cancelling
matrices. Let $\\/\mu=(7,5,4)/(1)$ and
$$
  M=\pmatrix{1&2&3\cr2&1&2\cr3&1&0\cr}
,
$$
so that $\beta=\csum{M}=(6,4,5)$. Then one already has
$\mu+\sum_{j<l}M\tr_j\nleh\mu+\sum_{j\leq{l}}M\tr_j$ for $l=0$, the concrete
values being $\mu=(1)\nleh(2,2,3)=\mu+M\tr_0$. Setting $\alpha=\mu+M\tr_0$ one
has $\[\alpha,]=(1,0,0,-4,\ldots)$, and the indices $(i_0,i_1)$ can be chosen
in various ways. For instance for $(i_0,i_1)=(0,1)$ one finds the sequence
$\[\alpha',]=(0,1,0,-4,\ldots)$ for $\alpha'=(1,3,3)$, which can be realised
by replacing~$M\tr_0$ by $(0,3,3)$. But one can also take $(i_0,i_1)=(1,2)$,
which permutes two identical entries of~$\[\alpha,]$ and therefore gives
$\alpha''=\alpha$, and an unchanged column~$M\tr_0$ (and indeed one has $\s223
=0$), or $(i_0,i_1)=(0,2)$, which realises $\alpha'''=(1,2,4)$ with
$\[\alpha''',]=(0,0,1,-4,\ldots)$ by replacing~$M\tr_0$ by $(0,2,4)$.
Performing the indicated changes to~$M$, we find that the contribution to
$s_\mu{h_\beta}$ from~$M$ could be cancelled by either one of the
contributions from the matrices
$$
  M'=\pmatrix{0&1&2\cr3&2&3\cr3&1&0\cr}
,\qquad
  M''=\pmatrix{1&2&3\cr2&1&0\cr3&1&2\cr}
,\ttext{or}
  M'''=\pmatrix{0&1&0\cr2&1&2\cr4&2&3\cr}
.
$$
Indeed one has $\s754 =-\s484 =-\s736 =-\s259 $ (the subscripts were
respectively computed as $\mu+\rsum{M}$, \ $\mu+\rsum{M'}$, \
$\mu+\rsum{M''}$, and $\mu+\rsum{M'''}$), which shows that the contribution of
$M$ to the second member of~(\GVbasiseq) could cancel against that of~$M'$,
of~$M''$, or of~$M'''$; the rule of choosing the minimal possible $(i_0,i_1)$
would in fact select~$M'$. Apart from the operations involved in the
cancellation, the example illustrates the fact that when we cancel $s_\alpha$
against $s_{\alpha'}$ with $\alpha\neq\alpha'$, one might still have
$s_\alpha=0$; it also shows that cancelling a null term $s_\alpha$ in
$s_\mu\fold(h\beta_0..\beta_l)$ ``against itself'' might lead to the
cancellation of two distinct terms in~$s_\mu{h_\beta}$.

There is a beautiful graphical way of visualising this cancellation, and the
terms that survive it, due to Gessel and Viennot~\ref{Gessel Viennot}. Set
$\alpha^{(j)}=\mu+\sum_{j'<j}M\tr_{j'}$ for~$j\in\N$, so that $M$ identifies
the term $s_{\alpha^{(j)}}$ in the development of
$s_\mu\fold(h\beta_0..\beta_{j-1})$. The idea is to trace for each fixed
index~$i$ the evolution of the number~$\[\alpha^{(j)},i]$, as $j$ increases,
by means of a lattice path: a path connecting points of~$\Z^2$ using unit
steps in two orthogonal directions, in the current case down and to the right.
One has a separate path for each $i\in\N$; path~$i$ passes through all of the
points~$\(j,\[\alpha^{(j)},i]\)$ for $j\in\N$. These points are connected as
follows: after arriving at~$\(j,\[\alpha^{(j)},i]\)$, the path takes $M_{i,j}$
horizontal steps (so the second coordinate becomes $\[\alpha^{(j)},i]+M_{i,j}
=\[\alpha^{(j+1)},i]$), followed by a vertical step to arrive at
$\(j+1,\[\alpha^{(j+1)},i]\)$. If $i>0$, this portion of path~$i$ remains
disjoint from path~$i-1$ if and only if
$\[\alpha^{(j+1)},i]<\[\alpha^{(j)},i-1]$, in other words if the addition of
$M_{i,j}$ to~$\[\alpha^{(j)},i]$ leaves it inside what we called the safe
interval~$I_i$ in the proof of proposition~\Piericomplprop. Thus all paths
will remain disjoint when taking into account column~$M\tr_j$ if and only if
$\alpha^{(j)}\leh\alpha^{(j)}+M\tr_j=\alpha^{(j+1)}$; they will remain
disjoint forever if and only if $M\in\IE{\\'/\mu}$ for some~$\\'\in\Part$.

For any~$M$, path~$i$ starts at $(0,\[\mu,i])$ and eventually runs along the
vertical line $\(*,\[(\mu+\rsum{M}),i]\)$. The term for~$M$ in the second
member of~(\GVbasiseq) is nonzero if and only if these vertical lines, for
$i\in\N$, form a permutation of the lines $(*,\[\\,i])_{i\in\N}$, whose sign
then gives the coefficient of the term. If moreover $\mu+\rsum{M}$ is a
partition, which must necessarily be~$\\$ given the previous condition, then
this permutation will be the identity, so that path~$i$ runs from the
point~$(0,\[\mu,i])$ to the line $(*,\[\\,i])$. This is in particular (but not
exclusively) the case whenever the paths are disjoint, so that
$M\in\IE{\\/\mu}$.

As an example of a term that survives the cancellation, and of the family of
disjoint paths that corresponds to it, we recall from~\Sec\prelimsec\ the
integral encoding~$M\in\IE{\\/\mu}$ of the tableau~$T\in\SSTx{\\/\mu}$
of~(\Tdef), where $\\/\mu=(9,8,5,5,3)/(4,1)$. To facilitate recognition of the
compositions (in fact partitions)~$\alpha^{(j)}$, we have written the
coefficients of~$\mu$ along the left border of~$M$. The paths, taken from
right to left, correspond to the rows of the matrix from top to bottom, and to
the rows of~$T$ in the same order.
\bigdisplay
\advance\abovedisplayskip by -6pt
T:\quad
\Skew(4:0,2,4,5,5|1:0,1,3,4,6,6,6|0:1,1,2,4,5|0:2,3,4,6,6|0:5,6,6)
\qquad
\bordermatrix{\mu&\multispan8\hss$M$\hss\cr
         4&1&0&1&0&1&2&0&0\cr
         1&1&1&0&1&1&0&3&0\cr
         0&0&2&1&0&1&1&0&0\cr
	 0&0&0&1&1&1&0&2&0\cr
	 0&0&0&0&0&0&1&2&0\cr
	 0&0&0&0&0&0&0&0&0\cr}
\maybebreak\qquad{}
\epsfxsize=0.28\hsize
\vcenter{\hbox{\epsfbox{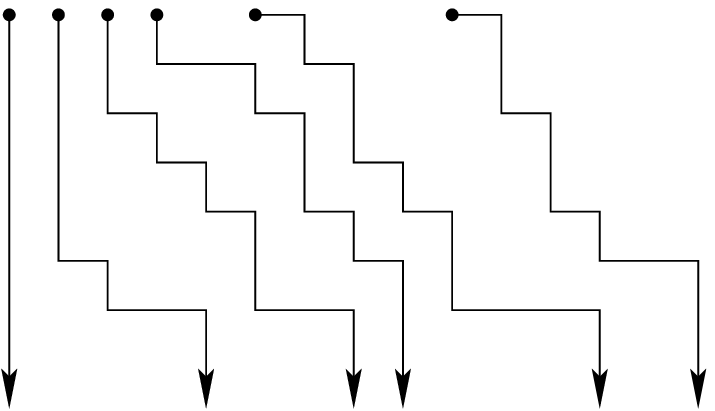}}}
$$
On the other hand, for the terms that must be cancelled, we have that whenever
two paths of the associated family touch, then a family of paths with opposite
contribution can be formed by exchanging these two paths after the point of
contact; if the term is nonzero, this exchange will change the sign of the
permutation defined by the family. This exchange corresponds to the type of
cancellation we considered above. Indeed if these are paths $i_0$~and~$i_1$,
and they have point~$(j,m)$ in common, then
$$\max(\[\alpha^{(j)},i_0],\[\alpha^{(j)},i_1])
  \leq m
  \leq \min(\[\alpha^{(j+1)},i_0],\[\alpha^{(j+1)},i_1])
,\label(\touchineq)
$$
and the exchange corresponds to transferring the difference
$\[\alpha^{(j+1)},i_0]-\[\alpha^{(j+1)},i_1]$ from $M_{i_0,j}$ to $M_{i_1,j}$
(which leaves both entries~$\geq0$ due to the inequality between the outer
members of~(\touchineq)), followed by the interchange for all $j'>j$ of the
entries $M_{i_0,j'}$~and~$M_{i_1,j'}$.

Looking back at proposition~\Piericomplprop\ from this path point of view, we
see that in its proof there is even more freedom in matching cancelling terms
than we used: the existence of a common point on paths $i_0$~and~$i_1$ is all
that is needed in order to obtain a cancellation by a transfer between parts
$i_0$~and~$i_1$, while we used the stricter requirement $f(i_0)=f(i_1)$. (To
make this clear, recall that the function~$f$ was defined relative to the
multiplication by a single factor~$h_j$ and under the assumption that
$\alpha^{(j)}$ is a partition, by the condition
$\[\alpha^{(j+1)},i]\in{I_{f(i)}}$, where
$I_0=\setof{k\in\Z}:\[\alpha^{(j)},0]\leq{k}\endset$ and
$I_i=\setof{k\in\Z}:\[\alpha^{(j)},i]\leq{k}<\[\alpha^{(j)},i-1]\endset$
for~$i>0$. Then for $i_1>i_0$ there is a common point with vertical
coordinate~$j$ on paths $i_0$ and~$i_1$ if and only if $i_0\geq{f(i_1)}$,
which is certainly the case if $f(i_1)=f(i_0)$.)

This extra freedom does not increase the number of terms that cancel, it just
gives more possibilities of choosing $(i_0,i_1)$ for those terms that do allow
cancellation. Of course, to complete a proof that the second member
of~(\GVbasiseq) reduces to the third member, one must give up this freedom and
make a concrete choice that defines a sign-reversing involution on the set of
cancelling terms. A fairly natural choice to do so is to choose the
top-leftmost point that lies on more than one path, in other words take $j$
minimal in the discussion above and then $m=\[\alpha^{(j)},i_0]$ where
$i_0=\max\setof{i}:\[\alpha^{(j+1)},i+1]\geq\[\alpha^{(j)},i]\endset$; then
$i_1=i_0+1$, and only paths $i_0$~and~$i_1$ pass through $(j,m)$, which makes
it easy to show that after exchange of paths this rule finds the same
point~$(j,m)$ and the same pair~$(i_0,i_1)$. Indeed the fact that we took some
pair~$(i_0,i_1)$ with $f(i_1)=f(i_0)$ in the proof of
proposition~\Piericomplprop\ seems a bit artificial in retrospect (regardless
of exactly which pair was chosen); in particular there may be no such pair
with $i_1=i_0+1$, so the common point used need not be the first one where
path~$i_1$ meets another path.

One can proceed analogously for equation~(\Kostkatralteq), comparing it with
the description of~$\smash{K'_{\\/\mu,\beta}}$ as counting transpose
semistandard tableaux of shape~$\\/\mu$ and weight~$\beta$. We have not
defined binary encodings of such tableaux~$T$ (or, for that matter, integral
encodings), but if we take as columns differences of successive partitions
of~$T$, we obtain a binary matrix~$M$ whose transpose is the binary encoding
of the transpose of~$T$. The situation is in fact quite similar to the one
considered above for~(\Kostkaalteq), as we have $\csum{M}=\beta$ and
$\rsum{M}=\\-\mu$, whence $M$ contributes $+1$ to the summation
in~(\Kostkatralteq), and all such matrices account for the entire value of the
summation. Forming a generating series by~$\beta$, we obtain
$$
  \sum_{\beta\in\comp}K'_{\\/\mu,\beta}X^\beta
 =\sum_{M\in\bMat}\eps(\mu+\rsum{M},\\)X^{\csum{M}}
 =\sum_{M\tr\in\BE{\\\tr/\mu\tr}}X^{\csum{M}}
.\label(\GVtrbasiseq)
$$
The first equality is related to the Von-N\"agelsbach-Kostka identity, since
the first member of the equation equals $s_{\\\tr/\mu\tr}[X_\N]$
by~(\Schurgenf), and the second member can be obtained by
substituting~(\ebetaeq) into~(\preVonNagelsbachKostka).

Before we consider the reduction of the second member of~(\GVtrbasiseq) to its
third member by cancellations, we transform the equation slightly by taking
transposes of tableaux and matrices, so that we will be dealing with ordinary
semistandard tableaux and their binary encodings. At the same time we replace
$\\/\mu$ by~$\\\tr/\mu\tr$ so that the members of the equation represent
$s_{\\/\mu}[X_\N]$ rather than $s_{\\\tr/\mu\tr}[X_\N]$; we get
$$
% \sum_{\alpha\in\comp}K_{\\/\mu,\alpha}X^\alpha=
  \sum_{M\in\bMat}\eps(\mu\tr+\csum{M},\\\tr)X^{\rsum{M}}
 =\sum_{M\in\BE{\\/\mu}}X^{\rsum{M}}
.\label(\GVtreq)
$$
Let $M\in\bMat$ with $\rsum{M}=\alpha$. Each row~$M_i$ represents a
monomial~$X^{M_i}$ in the factor~$e_{\alpha_i}$ of~$e_\alpha$, and the entire
matrix represents an occurrence of the monomial~$X^{\csum{M}}$ in the
expansion of~$e_\alpha$. Put $\beta^{(i)}=\mu\tr+\smash{\sum_{i'<i}M_{i'}}$
for~$i\in\N$, then $M\in\BE{\\/\mu}$ means that $\beta^{(i)}\in\Part$ for
all~$i$ (the remaining requirements for the relations
$\beta^{(i)}\lev\beta^{(i+1)}$ are automatically satisfied) and that
$\beta^{(i)}=\\\tr$ for $i$ sufficiently large (in other words
$\mu\tr+\csum{M}=\\\tr$). As $i$ increases, $M$ identifies a term
$s_{\beta^{(i)}}$ in the development of
$s_{\mu\tr}\fold(e\alpha_0..\alpha_{i-1})$, until possibly at some point
$\beta^{(i)}\notin\Part$, in which case $s_{\beta^{(i)}}=0$. If this never
happens, then $M$ also identifies a term $s_{\mu\tr+\csum{M}}$ in the
expression of the final value $s_{\mu\tr}e_\alpha$ in the basis of Schur
functions. In this case $M$ gives equal contributions to both members
of~(\GVtreq), which contributions are nonzero if and only if
$\mu\tr+\csum{M}=\\\tr$ or equivalently $M\in\BE{\\/\mu}$ (the contributions
are zero when $M\in\BE{\\'/\mu}$ for some partition $\\'\neq\\$; as before
this case needs no special consideration).

Now suppose to the contrary that one does have $s_{\beta^{(i)}}=0$ for
some~$i$, and fix~$i$ to the minimal such value. We see that the situation is
somewhat simpler than before, in that there is no need to modify row~$M_i$ to
find a matching term that accounts for the cancelling of this contribution to
$s_{\mu\tr}\fold(e\alpha_0..\alpha_{i-1})$. However the contribution
$s_{\mu\tr+\csum{M}}$ to~$s_{\mu\tr}e_\alpha$ might still be nonzero, and in
order to find a matrix with a cancelling contribution (whose rows up to and
including~$M_i$ will be the same as those of~$M$), one needs to choose a
column index~$j$ that ``causes'' the vanishing of $s_{\beta^{(i)}}$.
Concretely, the fact that $s_{\beta^{(i)}}=0$, while $\beta^{(i)}$ was
obtained by adding the binary composition~$M_i$ to the
partition~$\beta^{(i-1)}$, implies that the sequence~$\[\beta^{(i)},]$ has
equal adjacent entries at some position (perhaps at several different
positions). So choose some~$j$ for which
$\[\beta^{(i)},j]=\[\beta^{(i)},j+1]$, for instance the minimal one. Then to
find a matrix whose term cancels that of~$M$, we may interchange the parts of
columns $j$~and~$j+1$ beyond row~$M_i$; the cancellation is assured by the
fact that this modification causes the interchange of the entries at positions
$j$~and~$j+1$ of the sequence~$\[(\mu\tr+\csum{M}),]$.

A simple example will make this more concrete. We take $\\/\mu=(7,5,2)/(2)$
and
$$
  M=
  \pmatrix{0&0&1&1&0&0&0\cr1&0&1&0&1&0&1\cr0&1&0&0&1&0&0\cr1&1&0&0&1&1&0\cr}
,
$$
so that $\alpha=\rsum{M}=(2,4,2,4)$. One has $\beta^{(0)}=\mu\tr=(1,1)$ and
$\beta^{(1)}=(1,1,1,1)$ which is still a partition, but
$\beta^{(2)}=(2,1,2,1,1,0,1)$ is not, and indeed $s_{\beta^{(2)}}=0$. As cause
for the vanishing of~$s_{\beta^{(2)}}$ one can indicate either of the two
repetitions in the sequence $\[\beta^{(2)},]
=(+1,-1,-1,-3,-4,-6,-6,-8,-9,\ldots)$, and interchanging either columns
$1$~and~$2$ or columns $5$~and~$6$ beyond row~$1$ leads to two matrices whose
contribution could cancel against that of~$M$, namely
$$
  M'=
  \pmatrix{0&0&1&1&0&0&0\cr1&0&1&0&1&0&1\cr0&0&1&0&1&0&0\cr1&0&1&0&1&1&0\cr}
\ttext{and}
  M''=
  \pmatrix{0&0&1&1&0&0&0\cr1&0&1&0&1&0&1\cr0&1&0&0&1&0&0\cr1&1&0&0&1&0&1\cr}
.
$$
Indeed one has $\s3321311 =-\s3141311 =-\s3321302 $ (the subscripts were
respectively computed as $\mu\tr+\csum{M}$, \ $\mu\tr+\csum{M'}$, and
$\mu\tr+\csum{M''}$), which shows that the contribution of $M$ to the left
hand side of~(\GVtreq) could cancel against that of $M'$ or of~$M''$; since
$\eps((3,3,2,1,3,1,1),\\\tr)=-1$, the former contribution is in fact
$-X^{(2,4,2,4)}$, and the latter are both $+X^{(2,4,2,4)}$. The rule suggested
above of choosing~$j$ would pair up $M$ with~$M'$.

Again we can visualise the cancellation using lattice paths that trace the
evolution of the individual terms of the sequence~$\[\beta^{(i)},]$ as $i$
increases. In order that all possible paths (subject only to the condition of
eventually running in a prescribed direction) can be obtained, each
bit~$M_{i,j}$ should directly control the direction of segment~$i$ of
path~$j$, choosing between two possible directions (note the difference with
the correspondence above for integral matrices, where $M_{i,j}$ determined the
number of steps to the right on path~$i$ at horizontal level~$j$). For reasons
that will become clear below we choose the possible directions in the binary
case to be downwards (for bits~$0$) and leftwards (for bits~$1$). We want the
points on all paths that represent terms of the same
sequence~$\[\beta^{(i)},]$ to line up, and given the directions in which the
paths progress, this alignment must be along diagonals in~$\Z^2$. Thus we
arrive at the following rule: path~$j$ starts at the point
$(-\[\mu\tr,j],-\[\mu\tr,j])$ on the main diagonal, and passes through the
point $(i-\[\beta^{(i)},j],-\[\beta^{(i)},j])$, which lies $i$ diagonals below
the main diagonal, for each~$i\in\N$. If for $j\neq{j'}$ paths $j$~and~$j'$
have a point in common, then this happens for the same value of~$i$, and it
means that $\[\beta^{(i)},j]=\[\beta^{(i)},j']$. All paths will eventually run
straight downwards, with path~$j$ running along the vertical line
$(*,-\[(\mu\tr+\csum{M}),j])$. The matrix~$M$ will give a nonzero contribution
to the left hand side of~(\GVtreq) if and only if those lines for $j\in\N$ are
a permutation of the lines $(*,-\[\\\tr,j])$, and $M$ will then give a nonzero
contribution to the right hand side of~(\GVtreq) if and only if the paths
remain disjoint, in which case the permutation must be the identity, so that
$\mu\tr+\csum{M}=\\\tr$. We illustrate such a family of disjoint paths
corresponding to the same semistandard tableau~$T$ that was used before,
together with its binary encoding~$M$. Note that columns of~$T$ correspond to
those of~$M$, and to the paths, all in the same order from left to right.
\bigdisplay
T:\quad
\Skew(4:0,2,4,5,5|1:0,1,3,4,6,6,6|0:1,1,2,4,5|0:2,3,4,6,6|0:5,6,6)
\qquad
\bordermatrix{\mu\tr& 2&1&1&1&0&0&0&0&0&0 \cr
           &0&1&0&0&1&0&0&0&0&0\cr
           &1&1&1&0&0&0&0&0&0&0\cr
           &1&0&1&0&0&1&0&0&0&0\cr
	   &0&1&0&1&0&0&0&0&0&0\cr
\kern-.5em\smash{\raise 0.5\baselineskip \hbox{$M$}}
	   &0&0&1&1&1&0&1&0&0&0\cr
	   &1&0&0&0&1&0&0&1&1&0\cr
	   &0&1&1&1&1&1&1&1&0&0\cr
	   &0&0&0&0&0&0&0&0&0&0\cr}
\maybebreak\qquad{}
\epsfysize=10\baselineskip
\vcenter{\hbox{\epsfbox{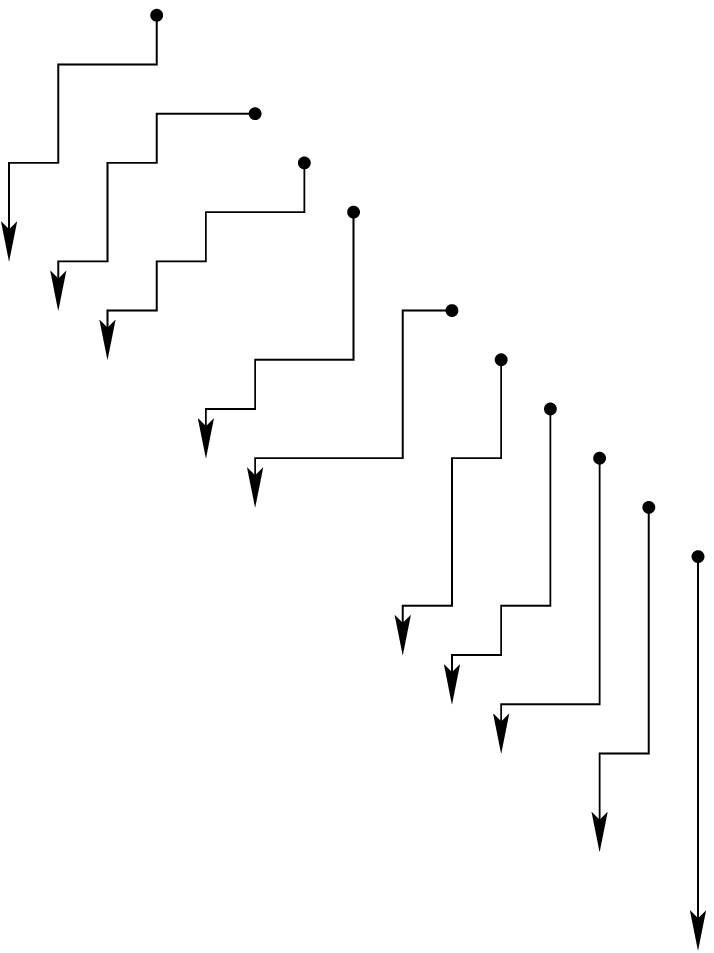}}}
$$
For families of paths that are not disjoint, cancellation can be realised as
before by choosing a pair of paths that meet, and exchanging their parts
beyond the first point of contact; for the matrices this means exchanging the
lower parts of the corresponding columns.

Now that we have two families of lattice paths related to any semistandard
tableau~$T$, a natural question is how the two are related. It turns out that
with a slight transformation the two families can be aligned so that there is
exactly one crossing of members of the two families for every square
of~$\set{\\/\mu}$. We show the two families shown above after transformation,
first individually and then in superposition.
\bigdisplay
\epsfxsize=0.27\hsize \vcenter{\hbox{\epsfbox{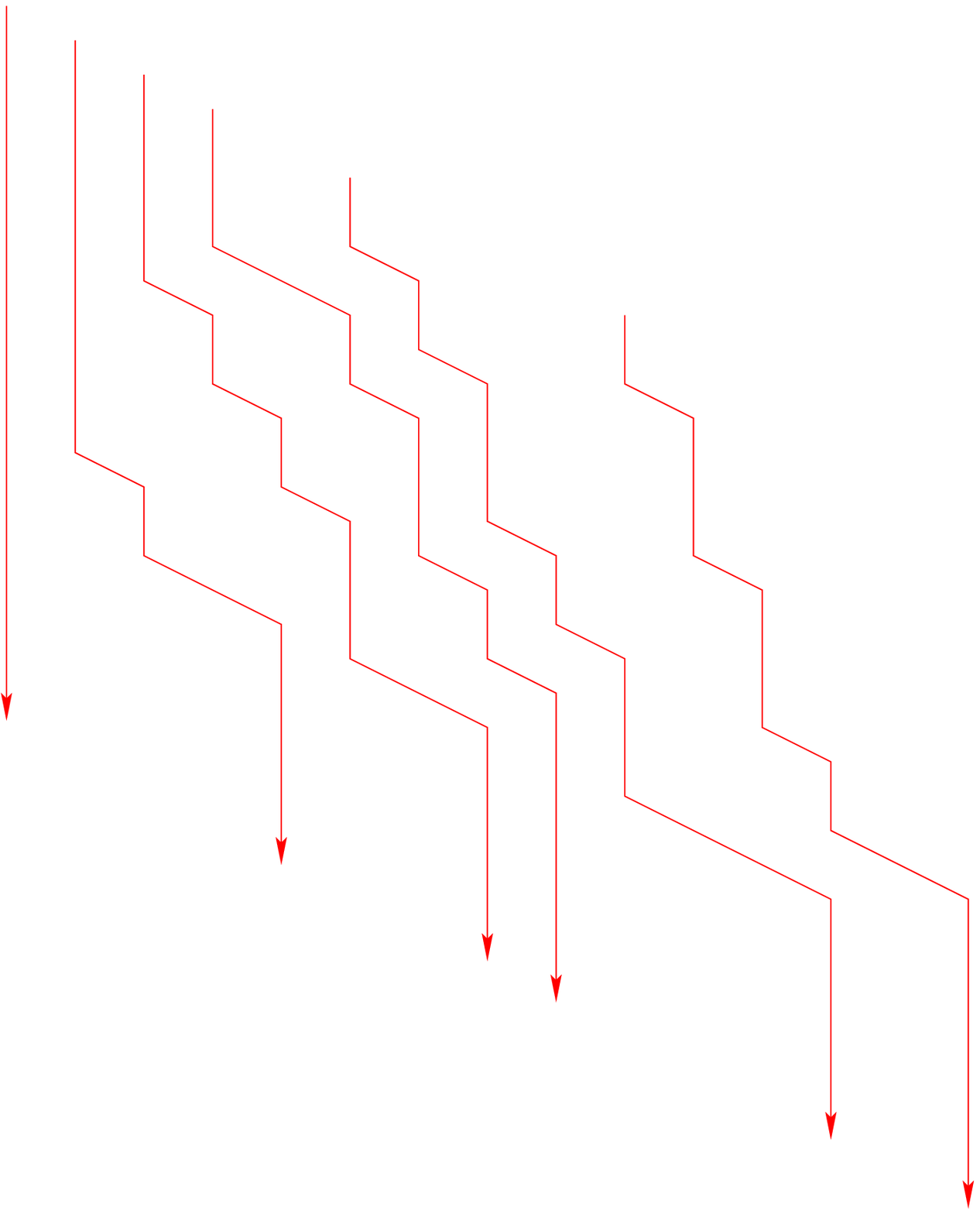}}}
\epsfxsize=0.27\hsize \vcenter{\hbox{\epsfbox{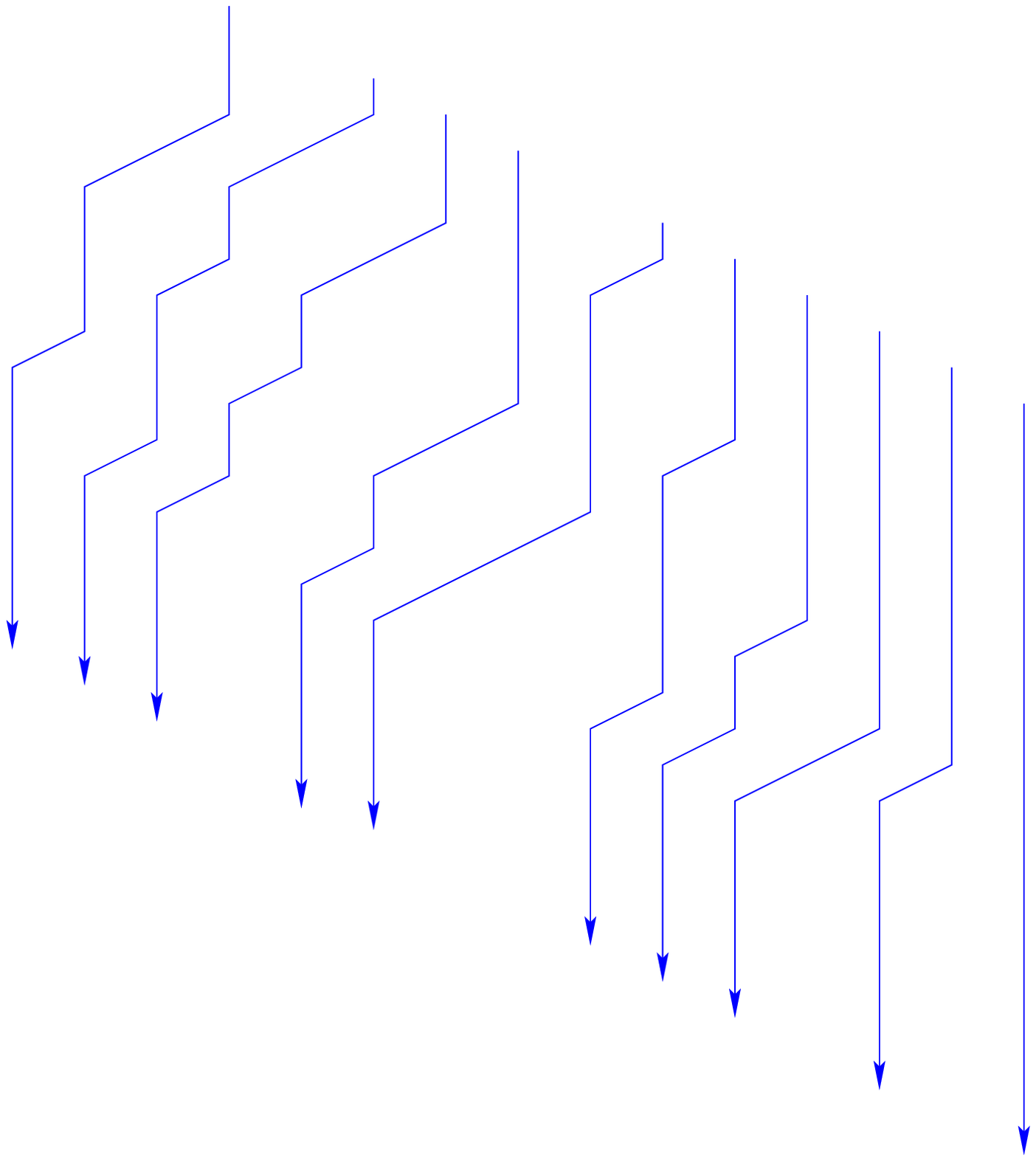}}}
\qquad
\epsfxsize=0.27\hsize \vcenter{\hbox{\epsfbox{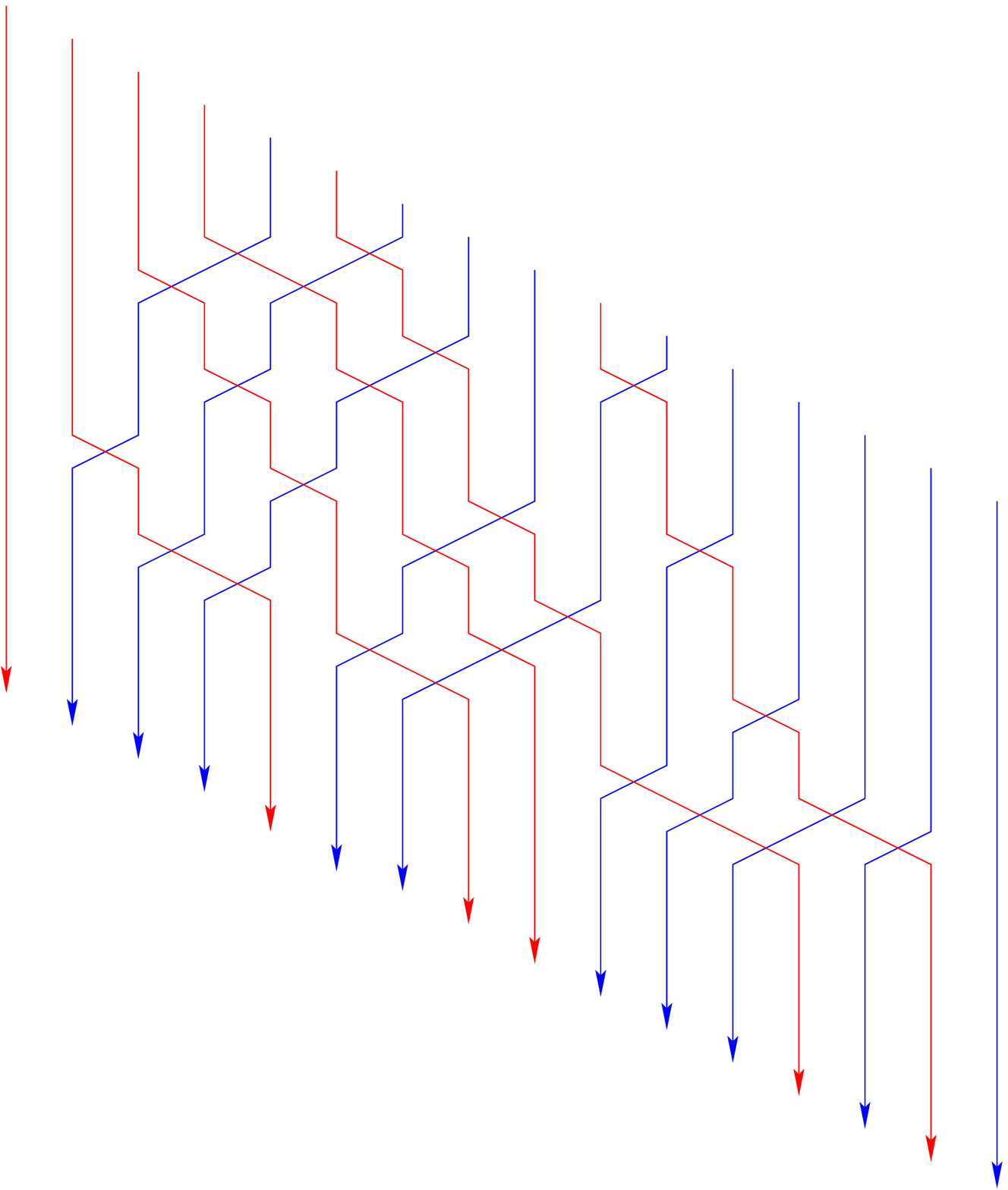}}}
$$
We have in fact, for aesthetic reasons, prefixed vertical segments of one
unit to all paths of the first family, and of half a unit to all paths of the
second family. The combinatorial explanation for this nice superposition is
that, if $T$ is given by the sequence of partitions $(\\^{(k)})_{k\in\N}$,
then the first family shows the evolution of the sets
$\setof\[\\^{(k)},i]:i\in\N\endset$ as~$k$ traverses~$\N$, and the second
family shows the evolution of the complementary sets
$\setof-1-\[(\\^{(k)})\tr,j]:j\in\N\endset$ (the term~$-1$ needed for
complementarity was left out of the coordinates of the second family for
simplicity, but is taken into account in the superposition). Indeed one finds
for matrices that are integral and binary encodings of a semistandard tableau
$T=(\\^{(k)})_{k\in\N}$ of shape~$\\/\mu$, that in our descriptions above
$\alpha^{(j)}=\\^{(j)}$, respectively that $\beta^{(i)}=(\\^{(i)})\tr$.

For the first family of paths, the values $\[\\^{(k)},i]$ are most clearly
read off at the middle of the vertical segment from $(k-1,\[\\^{(k)},i])$ to
$(k,\[\\^{(k)},i])$, while for the second family of paths, the values
of~$-1-\[(\\^{(k)})\tr,j]$ can be read off on a diagonal passing
\emph{through} the vertices; this explains the necessity of half a unit of
vertical displacement between the triangular lattices that are the images of
the lattice $\Z^2$ under the two transformations, in order to superimpose the
lines that correspond to the same sequence~$\[\\^{(k)},]$.

This section would not be complete without a \emph{visual} explanation of the
superposition of families of paths. The form of path~$i$ of the first family
describes the evolution of row~$i$, with $i$ fixed, of the Young diagrams
$\set{\\^{(k)}}$ for~$k\in\N$, so it has the same shape as the path along the
boundary of a single ``row'' (i.e., subset with first coordinate fixed) of the
$3$-dimensional ``block diagram'' corresponding to the semistandard tableau,
namely $\setof(i,j,k)\in\N^3:(i,j)\in\set{\\^{(k)}}\endset$ in which each
element is represented by a cube. To match the way we drew the lattice paths,
the coordinate~$k$ must be made to increase downwards, and the coordinate~$j$
to the bottom right. With the coordinate~$i$ increasing to the bottom left, we
obtain, for the tableau~$T$ of our example, the display below. The first
picture shows just the block diagram, with the top of the ``pillars'' for the
squares of the diagram $\set\mu=\set{(4,1)}$ omitted, so that the visible tops
of pillars correspond to the squares of $\set{\\/\mu}$. In the second picture
we have added the paths along the ``row boundaries'', which trace the boundary
of the block diagram at fixed value of~$i$.
\bigdisplay
\epsfxsize=0.30\hsize \vcenter{\hbox{\epsfbox{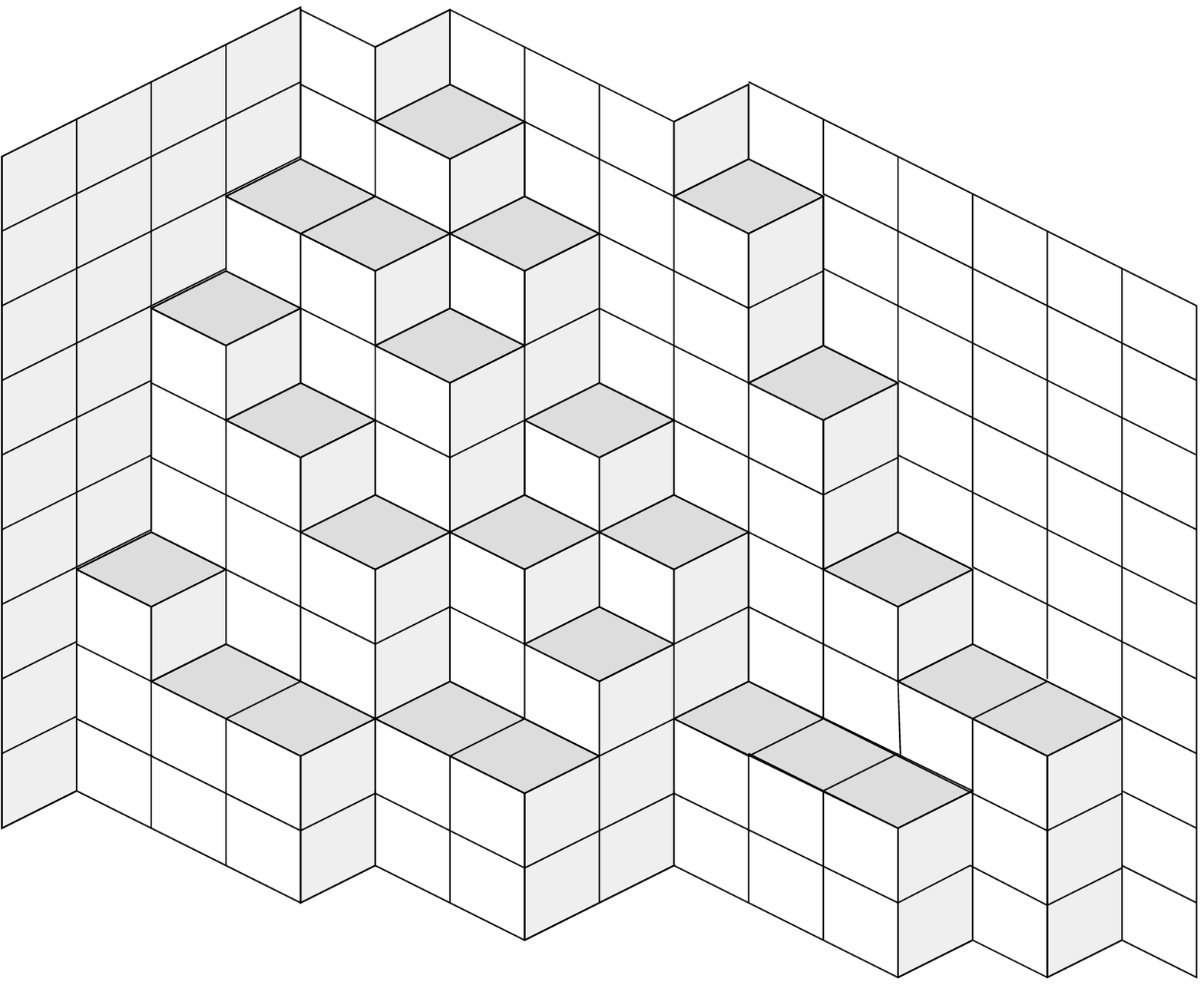}}}
\qquad\qquad
\epsfxsize=0.30\hsize \vcenter{\hbox{\epsfbox{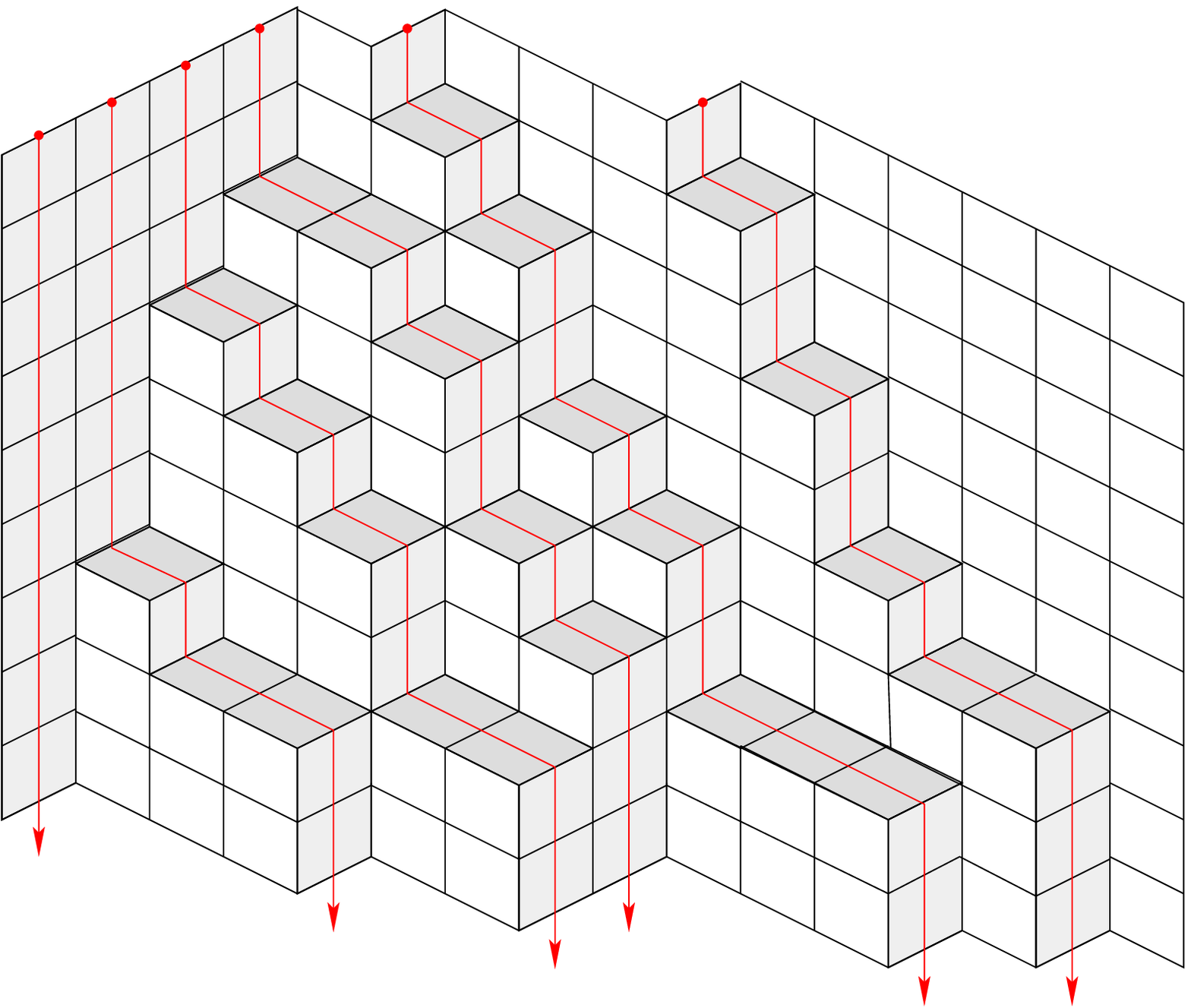}}}
$$
This display has as defects that the relative position of the paths is not
correct (so cases where paths should have a common point would still appear
disjoint), and the paths along ``column boundaries'' (with $j$~fixed, so
running from top right to bottom left) do not resemble the second family. Both
points can be remedied simultaneously, by raising (in the sense of
decreasing~$k$) each ``row''~$i$ (the plane of cubes with that value of~$i$)
by $i$ units. This amounts to using, instead of the block diagram of the
semistandard tableau $(\\^{(k)})_{k\in\N}$, the set
$\setof(i,j,k-i):(i,j)\in\set{\\^{(k)}}\endset\subset\N^3$. We show the
resulting diagram, and then add the two superimposed families of paths.
\bigdisplay
\epsfxsize=0.37\hsize \vcenter{\hbox{\epsfbox{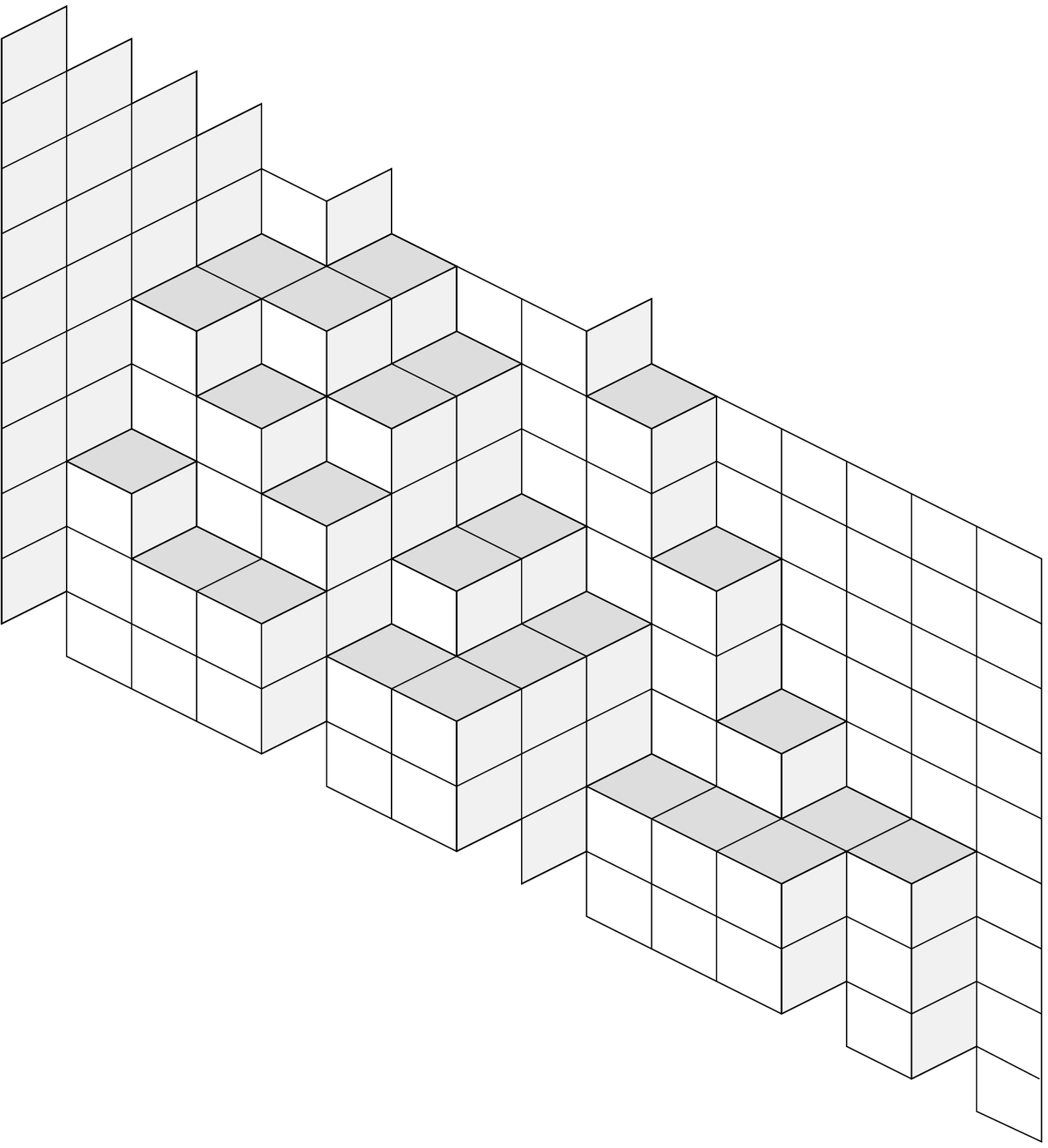}}}
\qquad
\epsfxsize=0.37\hsize \vcenter{\hbox{\epsfbox{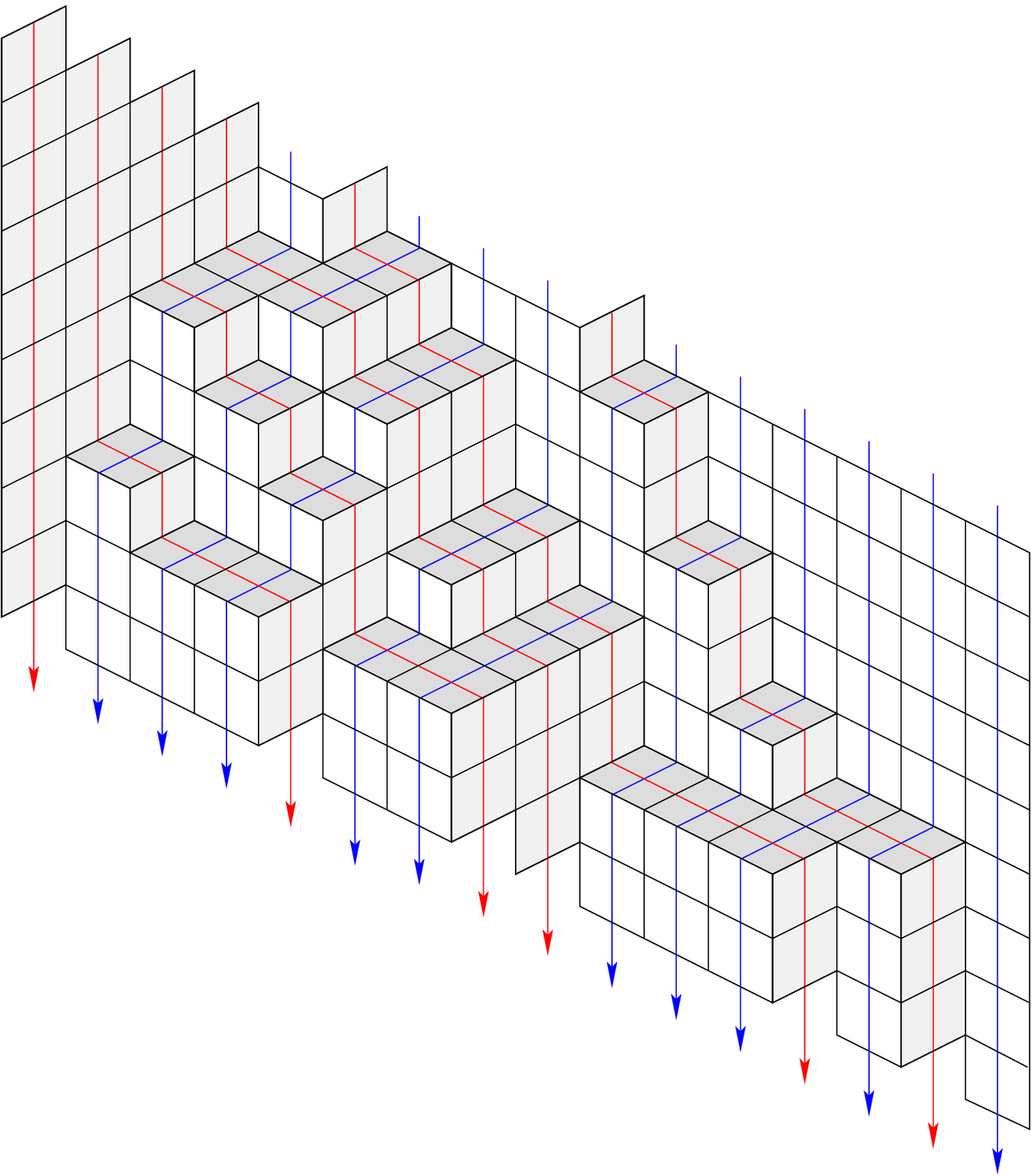}}}
$$

% Local IspellDict: default
% -*- mode: tex; tex-main-file: "root.tex"; -*-

\newsection The Littlewood-Richardson rule.
\labelsec\LRsec

We shall finally apply lemma~\Schurmul\ for $f=s_\\$ and even
for~$f=s_{\\/\kappa}$, so as to obtain an expression for multiplication by
(skew) Schur functions in the basis of Schur functions. As usual, the sport
will be to find a cancellation of the resulting alternating sum that leaves
only positive contributions (from representation theoretic interpretations of
Schur functions it is clear that the decomposition of the product in the basis
of Schur functions cannot involve negative coefficients). We shall thus find a
combinatorial description of the decomposition of a product $s_\\s_\mu$ into
Schur functions that is equivalent to the Littlewood-Richardson rule. In fact
we shall obtain a generalisation that was given in~\ref{Zelevinsky pictures},
which describes the scalar product of an arbitrary pair of skew Schur
functions. Our approach naturally leads to two different formulations of this
result, one in terms of binary matrices and another in terms of matrices with
coefficients in~$\N$. Both are equivalent to a more traditional formulation in
terms of Littlewood-Richardson tableaux, but the relation with cancellations
in alternating sums gives a perspective somewhat different from for instance
the one we took in~\ref{L-R intro}, where the Littlewood-Richardson rule was
approached via Robinson's correspondence and jeu de taquin. We would like to
note however that coplactic operations provide a very elegant means to unify
the two points of view, as we shall indicate in a sequel to this paper.

As an introduction we shall discuss the question that was evoked after
equation~(\Schurgenf), of proving bijectively that Schur functions are the
generating series of semistandard Young tableaux of given shape by weight. The
suggested proof was to apply lemma~\Schurmul\ for
$f[X_\N]=\sum_{T\in\SSTx\\}X^{\weight{T}}$ and $\beta=(0)$, and prove that its
right hand side reduces by cancellations to~$s_\\$. If we pretend not to know
already that this generating series $f[X_\N]$ is a Schur function, then there
is in fact another obligation, namely to prove that $f[X_\N]$ is a symmetric
function, so that lemma~\Schurmul\ can be applied in the first place. Proving
this symmetry can, surprisingly, be done in almost the same way as proving the
mentioned cancellation
$$
  \sum_{T\in\SSTx\\}s_{\weight{T}} = s_\\
.\label(\slambdacanc)
$$

The symmetry of~$f[X_\N]$ can be established by giving for each~$m$ an
involution on~$\SSTx\\$ that interchanges components $m$ and~$m+1$ of the
weight, which means, in terms of the display of the tableaux that the
frequencies of the entries $m$~and~$m+1$ must be interchanged, while
preserving the shape. Since there is no reason to change other entries than
$m$~and~$m+1$, this involution can be viewed as operating only on the
component shape $\\^{(m+1)}$ of tableaux $(\\^{(i)})_{i\in\N}$ in~$\SSTx\\$,
which shape is constrained by the relation
$\\^{(m)}\leh\\^{(m+1)}\leh\\^{(m+2)}$ in which the shapes
$\\^{(m)}$~and~$\\^{(m+2)}$ are fixed.

As for the cancellation~(\slambdacanc), the one tableau~$T$ in the summation
whose term survives the cancellation will be the one whose display has only
entries~$i$ in any row~$i$, and which therefore satisfies $\weight{T}=\\$. For
any other $T\in\SSTx\\$, let $i$ be minimal such that row~$i$ of the display
of~$T$ has some entry other than~$i$ (by column strictness, such entries must
be greater than~$i$), and let the final entry of that row, in
column~$j=\\_i-1$, have value~$m+1$ (so $m\geq{i}$). We shall consider this
entry as the one that triggers the cancellation of the term for~$T$, and in
constructing $T'\in\SSTx\\$ whose terms cancels that of~$T$, we shall leave
that entry intact, as well as all rows $i'<i$, so that the same entry~$m+1$
also triggers the cancellation for~$T'$. In fact, we shall only change entries
$m$~and~$m+1$ in the display of~$T$, if any, with the goal of interchanging
the values of $\[\weight{T},m]$ and of $\[\weight{T},m+1]$. The entries that
may be changed are in rows $i'\geq{i}$ and in columns $j'<j$ of the display
of~$T$, and one can treat the subtableau of entries $m,m+1$ in that region in
isolation: whatever subtableau of entries $m,m+1$ should replace them, it will
be compatible with the entry~$m+1$ at position~$(i,j)$. Moreover, since that
single entry causes the contribution to $\[\weight{T},]$ of the part outside
the specified region to have equal values at indices $m$~and~$m+1$, the goal
must be to interchange the frequencies of entries $m$~and~$m+1$ in the
mentioned subtableau. Thus, the problem is basically the same one as for
proving symmetry of the generating series~$f[X_\N]$.

The simplest involutions establishing that symmetry are those that were
introduced in \ref{Bender Knuth}, and which we shall call the Bender-Knuth
involutions. They were in fact introduced in the combinatorial proof of a
lemma that essentially states our equation~(\slambdacanc) (lemma~2 of the
cited paper), although an application to the symmetry of~$f$ immediately
follows it. The idea is straightforward: one ignores any columns containing
both entries $m$~and~$m+1$ (such entries cannot be changed), after which all
rows can be treated independently, and in each of them one separately
interchanges the frequencies of $m$~and~$m+1$; we refer to the proof
of~\ref{Stanley EC2, theorem 7.10.2} for details.

A simple example will illustrate these involutions and their application to
establishing~(\slambdacanc). Take the semistandard Young tableau displayed as
follows
\bigdisplay
   \Young(0,0,0,0,0,0,0|1,2,2,2,4,4,4|3,3,3,4,5,6,6|4,5,5,6,6|6,7,7)
$$
The Bender-Knuth involution that interchanges the components $5$~and~$6$ of
the weight, operates separately on the singleton sequence $[6]$ in row~$4$, on
the sequence $[5,5,6]$ in columns $1,2,3$ of row~$3$, and on the sequence
$[6,6]$ at the end of row~$2$ (the entries $5$~and~$6$ in column~$4$ are
frozen); the sequences are respectively turned into $[5]$, $[5,6,6]$
and~$[5,5]$. As for the cancellation~(\slambdacanc), it is triggered by the
final entry~$4$ of row~$1$, and therefore operates on entries with values
$3$~and~$4$; that final entry does not participate in the exchange, nor do the
frozen entries in column~$0$, so $[3,3,4]$ in row~$2$ and the initial $[4,4]$
in row~$1$ are respectively changed into $[3,4,4]$ and $[3,3]$. The tableaux
resulting for these two operations are
\bigdisplay
   \Young(0,0,0,0,0,0,0|1,2,2,2,4,4,4|3,3,3,4,5,5,5|4,5,6,6,6|5,7,7)
\ttext{and}
   \Young(0,0,0,0,0,0,0|1,2,2,2,3,3,4|3,3,4,4,5,6,6|4,5,5,6,6|6,7,7)
.
$$
One sees that the weights of these tableaux have the required relation to the
weight $\alpha=(7,1,3,3,5,3,5,2)$ of the initial tableau, notably the
weight $\alpha'=(7,1,3,4,4,3,5,2)$ of the final tableau satisfies
$s_\alpha+s_{\alpha'}=0$.

Now let us turn to more useful applications of lemma~\Schurmul, for which we
do not know the result beforehand. With respect to the application in the
previous example one could replace the choice $\beta=(0)$ by one taking
for~$\beta$ an arbitrary partition~$\mu$, which results in a formula for
products $s_\\s_\mu$ of Schur functions, or one could replace the choice
$f=s_\\$ by the more general $f=s_{\\/\kappa}$, resulting in a formula
for decomposing skew Schur functions into ordinary ones. In view of the
relation $\scalar<s_\\s_\mu,s_\nu>=\scalar<s_\\,s_{\nu/\mu}>$ the two problems
are equivalent, but it will be interesting to perform both generalisations at
once, which in fact produces no additional complications. So we apply
lemma~\Schurmul\ with $f=s_{\\/\kappa}$ (and hence
$f[X_\N]=\sum_{T\in\SSTx{\\/\kappa}}X^{\weight{T}}$) and $\beta=\mu\in\Part$,
which gives
$$
  s_{\\/\kappa}s_\mu=\sum_{T\in\SSTx{\\/\kappa}}s_{\mu+\weight{T}}
.\label(\LRgenbruteeq)
$$
We can obtain a somewhat more symmetric formulation by taking coefficients of
some Schur function~$s_\nu$ on both sides, which amounts to taking the scalar
product with it; using $\scalar<s_{\\/\kappa}s_\mu,s_\nu>
=\scalar<s_{\\/\kappa},s_{\nu/\mu}>$ this gives
$$
  \scalar<s_{\\/\kappa},s_{\nu/\mu}>
 =\sum_{T\in\SSTx{\\/\kappa}}\edgprm{\mu+\weight{T}}\nu
.\label(\LRgensymeq)
$$

We shall now describe a cancellation that can be applied to the right hand
side of~(\LRgenbruteeq), which was proposed in~\ref{Stembridge
Littlewood-Richardson}; in fact a rather similar cancellation argument, but
using a somewhat different terminology, can already be found
in~\ref{Berenstein Zelevinsky, Theorem~4}. Contrary to what we saw
for~(\slambdacanc) there will in general be more than one term that survives
the cancellation, but for those tableaux~$T$ whose term does cancel, an
entry~$m+1$ in its display will be indicated that triggers the cancellation,
and a tableau~$T'$ whose terms cancels that of~$T$ will be constructed by
possibly modifying entries $m$~and~$m+1$, in much the same way as before. In
particular all changes will be strictly to the left and weakly below the entry
triggering the cancellation.

The condition for cancellation is stated in terms of the weights of individual
columns of~$T$. Define the weight $\colweight{j}T$ of column~$j$ of~$T$ to be
the binary composition~$\alpha$ whose part~$\alpha_i\in\sing{0,1}$ tells
whether or not there is an entry~$i$ in column~$j$ of the display of~$T$, in
other words, $\colweight{j}T$ is the column $(M')\tr_j$ of the binary
encoding~$M'$ of~$T$. For later use also define the weight $\rowweight{i}T$ of
row~$i$ of~$T$ as the composition~$\beta$ whose part~$\beta_j$ gives the
number of entries~$j$ in row~$i$ of the display of~$T$; thus $\rowweight{i}T$
equals row $M_i$ of the integral encoding~$M$ of~$T$. One obviously has
$\weight{T} =\sum_{j\in\N}\colweight{j}T =\sum_{i\in\N}\rowweight{i}T$.

The terms in the right hand side of~(\LRgenbruteeq) that will survive the
cancellation are those for the tableaux~$T$ with the property that for all
$l\in\N$ one has $\mu+\sum_{j\geq{l}}\colweight{j}T\in\Part$; in particular
$\mu+\weight{T}$ is a partition, so that the term for such~$T$ is already
normalised. So once the cancellation is established, it will prove
$$
   s_{\\/\kappa}s_\mu
   =\sum_{T\in\SSTx{\\/\kappa}}
    \Kr{\forall{l\in\N}\:\mu+\sum_{j\geq{l}}\colweight{j}T\in\Part}
    s_{\mu+\weight{T}}
,\label(\LRgeneqcol)
$$
which is a generalisation of the Littlewood-Richardson rule. The tableaux that
satisfy the condition are those that we termed $\mu$-dominant in~\ref{L-R
intro, \Sec1.4} (although it requires a bit of work to see this), and for
$\mu=(0)$ these are the traditional Littlewood-Richardson tableaux; as one can
see, these describe the decomposition of $s_{\\/\kappa}$ in the basis of Schur
functions. By taking instead $\kappa=(0)$, one obtains a description of the
decomposition of $s_\\s_\mu$ in that basis, in terms of $\mu$-dominant
tableaux of shape~$\\$.

For any $T$ that fails the condition in~(\LRgeneqcol), there is a maximal~$l$
such that $\mu+\sum_{j\geq{l}}\colweight{j}T\notin\Part$; we fix this~$l$ and
put $\alpha=\mu+\sum_{j\geq{l}}\colweight{j}T$. Since $\alpha\notin\Part$
while $\mu+\sum_{j>l}\colweight{j}T=\alpha-\colweight{l}T\in\Part$, and
$\colweight{l}T$ is a binary composition, it must be that there exist some~$m$
with $\alpha_{m+1}=\alpha_m+1$, which can also be written as
$\[\alpha,m]=\[\alpha,m+1]$. We choose the minimal such~$m$, which then
necessarily also has $\colweight{l}T_m=0$ and $\colweight{l}T_{m+1}=1$, so
column~$l$ of the display of~$T$ contains no entry~$m$, but it does contain an
entry~$m+1$. That entry will be the one that is considered to trigger the
cancellation of the term for~$T$.

Since the entries in columns $j\geq{l}$ of the display of~$T$ are involved in
triggering the cancellation and determining the values of $l$~and~$m$, they
will not be altered during the construction of a tableau~$T'$ whose term
cancels that of~$T$; this will ensure that the cancellation is triggered
identically for~$T'$ as it was for~$T$. As before the subtableau of entries
$m,m+1$ in columns~$j<l$ can be treated in isolation, since the frozen
entry~$m+1$ in column~$l$ imposes no constraint. One has
$\mu+\weight{T}=\alpha+\sum_{j<l}\colweight{j}T$ and for~$T'$ one will
similarly have $\mu+\weight{T'}=\alpha+\sum_{j<l}\colweight{j}{T'}$; since
$\[\alpha,m]=\[\alpha,m+1]$, the desired cancellation
$s_{\mu+\weight{T}}+s_{\mu+\weight{T'}}=0$ can be assured by having
$\sum_{j<l}\colweight{j}{T'}$ differ from $\sum_{j<l}\colweight{j}T$ by the
interchange of the parts indexed by $m$~and~$m+1$. Therefore $T'$ can be
obtained from~$T$ by applying the Bender-Knuth involution for $m,m+1$ to the
subtableau in columns~$j<l$.

As an example consider the tableau~$T$ of~(\Tdef) with $\mu=(7,6,5,4,2)$; one
finds $l=5$, for which $\alpha=\mu+(0,0,1,0,1,2,3)=(7,6,6,4,3,2,3)$ and $m=5$,
and the entry $m+1=6$ at position~$(1,5)$ triggers the cancellation. The
tableau~$T'$ is found by applying the Bender-Knuth involution for
interchanging $5$~and~$6$ to the subtableau formed by the $5$~leftmost columns
of~$T$, which transforms
\bigdisplay
T=\Skew(4:0,2,4,5,5|1:0,1,3,4,6,6,6|0:1,1,2,4,5|0:2,3,4,6,6|0:5,6,6)
\ttext{into}
T'=\Skew(4:0,2,4,5,5|1:0,1,3,4,6,6,6|0:1,1,2,4,5|0:2,3,4,5,6|0:5,5,6)
.\label(\BeKnex)
$$
Indeed one has $s_{\mu+\weight{T}}+s_{\mu+\weight{T'}} =\s9986647 +\s9986665
=0$.

The condition for $\mu$-dominance that appears in~(\LRgeneqcol) may appear
somewhat strange at first sight, but it has a remarkable resemblance to the
condition that characterises binary encodings of tableaux. To bring out this
resemblance, we shall replace the tableau~$T$ in the summation by its binary
encoding~$M$, so that $M\tr_j$ will replace $\colweight{j}T$, and $\rsum{M}$
will replace $\weight{T}$. Also we shall apply our cancellations to the right
hand side of~(\LRgensymeq) rather than of~(\LRgenbruteeq). One has
$\csum{M}=\\\tr-\kappa\tr$ for any $M\in\BE{\\/\kappa}$, and those~$M$ giving
nonzero terms after cancellation satisfy $\rsum{M}=\nu-\mu$. Spelling out
explicitly the condition that $M\in\BE{\\/\kappa}$, we therefore obtain the
following statement, equivalent to~(\LRgeneqcol).

\proclaim Theorem. \LRthmone
(generalised Littlewood-Richardson rule, binary formulation)
For all $\kappa,\\,\mu,\nu\in\Part$:
$$
  \scalar<s_{\\/\kappa},s_{\nu/\mu}>
 =\sum_{M\in\bmat{\nu-\mu}{\\\tr-\kappa\tr}}
  \Kr{\forall k\: \kappa\tr+\sum_{i<k}M_i\in\Part}
  \Kr{\forall l\: \mu+\sum_{j\geq{l}}M\tr_j\in\Part}
.\eqno\copy\QEDbox
$$

There is another formulation of the condition of $\mu$-dominance than the one
that appears in~(\LRgeneqcol), which uses row weights rather than column
weights, and which is actually closer to the way the notion was defined
in~\ref{L-R intro}. One can determine independently for each value of~$m$
whether any entry~$m+1$ in the display of~$T$ is a candidate for triggering
cancellation. Although in case of multiple candidates, only the topmost one in
the rightmost column containing any candidates will be selected, that rule has
no influence on the set of tableaux that will survive cancellation, namely
those for which no candidates exist at all. In order for an entry~$m+1$ in
column~$l$ to be candidate for triggering cancellation, one must have
$\[\alpha,m]=\[\alpha,m+1]$ for $\alpha=\mu+\sum_{j\geq{l}}\colweight{j}T$,
and $l$ must be the largest value for which this condition holds. Then there
will not be any entry~$m$ above that entry~$m+1$ in the same column~$l$
of~$T$, and this means that if the entry~$m+1$ is in row~$k$, one also has
$\[\alpha',m]=\[\alpha',m+1]$ for $\alpha'
=\mu+\sum_{i<k}\rowweight{i}T+\beta$, where $\beta$ is the weight of the part
of row~$k$ of~$T$ consisting of the candidate entry and all entries to its
right.

Conversely any entry~$m+1$ for which this condition in terms of~$\alpha'$
holds, and for which $k$ is minimal, cannot have an entry~$m$ in the same
column, and is therefore a candidate for triggering cancellation. Then,
assuming there are no candidate entries~$m+1\binoppenalty=10000$ in
rows~$i<k$, one may test if there is \emph{any} such entry~$m+1$ in row~$k$,
by considering $\alpha''=\mu+\sum_{i<k}\rowweight{i}T+\beta'$ where $\beta'$
the weight of the part of row~$k$ consisting of entries exceeding~$m$ (in
other words $\beta'$ is obtained from $\rowweight{k}M$ by clearing its parts
at indices~$\leq{m}$): there will be such a candidate if and only if
$\alpha''_{m+1}>\alpha''_m$. This use of weights of partial rows can be
avoided in the formulation, by comparing parts of two different compositions:
if one defines $\alpha^{(i_0)}=\mu+\sum_{i<i_0}\rowweight{i}T$ for~$i_0\in\N$,
then the final condition $\alpha''_{m+1}>\alpha''_m$ is equivalent to
$\ifblownup\else\expandafter\smash\fi{\alpha^{(k)}_m}
<\ifblownup\else\expandafter\smash\fi{\alpha^{(k+1)}_{m+1}}$. For instance in
our example above, the fact that some entry~$6$ in row~$1$ is a candidate for
triggering cancellation, can be deduced from $\alpha^{(1)}_5<\alpha^{(2)}_6$
where $\alpha^{(1)} =\mu+(1,0,1,0,1,2)=(8,6,6,4,3,\underline2)$ and
$\alpha^{(2)} =\alpha^{(1)}+(1,1,0,1,1,0,3) =(9,7,6,5,4,2,\underline3)$ (only
the underlined entries are relevant).

Now let us return to formulating $\mu$-dominance, the condition stating that
the term for~$T$ is not cancelled at all. This means that the condition just
stated fails for all~$k$ and for all~$m$. The failure for given~$k$ and for
all~$m$ simply becomes $\alpha^{(k)}\leh\alpha^{(k+1)}$, which means that we
can reformulate~(\LRgeneqcol) as follows:
$$
   s_{\\/\kappa}s_\mu
   =\sum_{T\in\SSTx{\\/\kappa}}
    \Kr{\forall{k\in\N}\:\mu+\sum_{i<k}\rowweight{i}T\leh
                         \mu+\sum_{i\leq{k}}\rowweight{i}T}
    s_{\mu+\weight{T}}
.\label(\LRgenerow)
$$
Note that, somewhat surprisingly, the sequence $(\alpha^{(i)})_{i\in\N}$ of
compositions actually is a semistandard tableau when $T$ survives the
cancellation; its shape is $\mu+\weight{T}/\mu$, and in the terminology
of~\ref{L-R intro} it is the companion tableau of~$T$ witnessing the
$\mu$-dominance of~$T$. If this time we replace $T$ in the summation by its
integral encoding, we obtain a symmetry in the formulation that is even more
striking than the one in theorem~\LRthmone. Again we take (\LRgensymeq) as a
starting point, and obtain the following statement, which is equivalent
to~(\LRgenerow), and therefore to~(\LRgeneqcol), like the previous theorem.
This time the condition that selects the matrices whose term survives the
cancellation is exactly like the condition for being an integral encoding in
the first place, except for the interchange of rows and columns; therefore we
add a succinct reformulation of the conditions, using the notation introduced
in definition~\encodingdef.

\proclaim Theorem. \LRthmtwo
(generalised Littlewood-Richardson rule, integral formulation)
For all $\kappa,\\,\mu,\nu\in\Part$
$$
\def\breakline{\kern-11mm\cr}% must be a macro to be used conditionally
\eqalignno
{ \scalar<s_{\\/\kappa},s_{\nu/\mu}>
\ifblownup\breakline\fi
&=\sum_{M\in\mat{\\-\kappa}{\nu-\mu}}
  \Kr{\forall l\:\kappa+\sum_{j<l}M\tr_j\leh\kappa+\sum_{j\leq{l}}M\tr_j}
  \Kr{\forall k\:   \mu+\sum_{i<k}M   _i\leh   \mu+\sum_{i\leq{k}}M   _i}
\cr
&=\sum_{M\in\Mat}
  \Kr{M\in\IE{\\/\kappa}}\Kr{M\tr\in\IE{\nu/\mu}}
. & \copy\QEDbox
\cr}
$$

We have completed the description of our final application of lemma~\Schurmul,
and of the result of applying cancellations to the expression it yields. There
is however one more variation that we would like to present. In theorems
\LRthmone~and~\LRthmtwo\ we have given expressions for
$\scalar<s_{\\/\kappa},s_{\nu/\mu}>$ in terms of matrices that treat rows and
columns similarly. But our starting point for these expressions, equation
(\LRgensymeq), does not exhibit such a symmetry when we replace the
semistandard tableaux in the summation by their binary or integral encoding
matrices. We would like to give expressions involving alternating sums with
the same kind of symmetry as the mentioned theorems, and from which those
theorems can be obtained by cancellations. This will require that we proceed
backwards, introducing a second layer of signs by a kind of inverse
cancellation process. In fact we shall use equations
(\Kostkaalteq)~and~(\Kostkatralteq) to translate the Kostka numbers, which are
implicit in the summation over semistandard tableaux, into alternating sums
over matrices.

Concretely, we write the right hand side of~(\LRgensymeq) as
$\sum_{\beta\in\comp}K_{\\/\kappa,\beta}\edgprm{\mu+\beta}\nu$, and then
we can apply (\Kostkaalteq), which gives
$$
 \scalar<s_{\\/\kappa},s_{\nu/\mu}>
 =\sum_{M\in\Mat} \edgprm{\kappa+\rsum{M}}\\\edgprm{\mu+\csum{M}}\nu
,\label(\LRgeninteq)
$$
or we can apply (\Kostkatralteq) to
$K_{\\/\kappa,\beta}=K'_{\\\tr/\kappa\tr,\beta}$, which gives (transposing the
binary matrices~$M$)
$$
 \scalar<s_{\\/\kappa},s_{\nu/\mu}>
 =\sum_{M\in\bMat} \edgprm{\kappa\tr+\csum{M}}{\\\tr}\edgprm{\mu+\rsum{M}}\nu
.\label(\LRgenbineq)
$$

It is clear from the way we have obtained these equations, that one can go
from them to theorems \LRthmtwo~and~\LRthmone, respectively, by applying two
successive phases of cancellation. The first phase is a Gessel-Viennot
cancellation of terms with opposite values of $\edgprm{\kappa+\rsum{M}}\\$
respectively of $\edgprm{\kappa\tr+\csum{M}}{\\\tr}$, as in~(\GVbasiseq)
and~(\GVtreq), which gives us back~(\LRgensymeq) with tableaux replaced by
their integral or binary encodings. The second cancellation phase cancels
within the remaining terms pairs with opposite values of the remaining factor
$\edgprm{\mu+\csum{M}}\nu$, using the cancellation based on the Bender-Knuth
involutions that was given following~(\LRgeneqcol), translated in terms of
integral or binary encodings.

It is interesting to compare the way in which these two cancellation phases
proceed. The Gessel-Viennot cancellations search for a failure of the
condition $M\in\IE{\\/\kappa}$ respectively $M\in\BE{\\\tr/\kappa\tr}$ by
successively computing the compositions~$\alpha^{(i)}$ that would define the
shapes (or their transposes in the binary case) forming the semistandard
tableau encoded by~$M$ if one exists; as soon as one finds
$\alpha^{(i)}\nleh\alpha^{(i+1)}$ respectively
$\alpha^{(i)}\nlev\alpha^{(i+1)}$, the cancellation of the term for~$M$ is
detected, and the construction of a cancelling matrix~$M'$ starts. Up to this
point the second cancellation phase proceeds in quite the same way, differing
only in whether weights of rows or columns of~$M$ are successively added to
the compositions~$\alpha^{(i)}$, and in the fact that for the binary case
these columns are taken in the opposite order.

But after detecting cancellation and selecting a pair of indices responsible
for it, the Gessel-Viennot cancellations construct~$M'$ from~$M$ by simply
interchanging entries for those indices in the rows or columns that had not
yet been considered; in the integral case there is also a transfer at those
indices in the column in which cancellation was detected itself, which can be
interpreted as interchange of those parts of the values of the entries
concerned that remain after the cancellation has been triggered (the parts of
their paths after the first point of contact). Such a simple rule would not
work in the Bender-Knuth case: in the binary case weak increase along rows
would be violated, while in the integral case strict increase in columns would
be violated (in either case this is the only one of the two conditions that
can be violated, as the other one is assured by the way tableaux are encoded
by matrices).

What the Bender-Knuth involutions actually do in terms of integral and binary
encodings is not easily described directly, although one can say in general
that the transfers between entries are more limited than what they would be
for a Gessel-Viennot cancellation. This can be seen by considering the
encodings of the tableaux $T,T'$ in~(\BeKnex). The changes only involve
entries $5$~and~$6$, so we just show the changes to columns $5$~and~$6$ of the
integral encoding and to rows $5$~and~$6$ of the binary encoding:
\bigdisplay
  \pmatrix{2&0\cr0&\underline3\cr1&0\cr0&2\cr1&2\cr}
  \leftrightarrow
  \pmatrix{2&0\cr0&\underline3\cr1&0\cr1&1\cr2&1\cr}
\text{and}
\matrix
{ \pmatrix{1&0&0&0&1&0&0&1&1\cr0&1&1&1&1&\underline1&1&1&0\cr}
 \cr
  \strut\updownarrow
 \cr
  \pmatrix{1&1&0&1&1&0&0&1&1\cr0&0&1&0&1&\underline1&1&1&0\cr}
 \cr
}
.\nn
$$
The entries where cancellation is triggered have been underlined (in the
integral case it is the final unit of the entry~$3$, i.e., a value~$2$ at the
same position would not trigger cancellation). One sees that the parts of the
columns below that entry (in the integral case) or of the rows to the left of
that entry (in the binary case) have not been completely interchanged. The
precise changes depend on the shape $\\^{(5)}=(7,5,4,3)$ in the tableau~$T$ to
which the entries~$5$ in its display are added, and in terms of the matrices
this shape depends both on~$\mu$ and on the columns to the left of those that
were shown, respectively the rows above those that were shown.

One may wonder, both for (\LRgeninteq) and for~(\LRgenbineq), if it is
possible to find a pair of cancellations that reduces them to theorem
\LRthmtwo\ respectively to theorem~\LRthmone, in a manner that better respects
the symmetry present both in the initial and in the final expressions. In
particular one would like to be able to apply the two cancellations in either
order, with the cancellation being applied last operating only on the terms
that survived the first cancellation. This requires that either cancellation
respects the condition expressing survival under the other cancellation: it
will never cancel a survivor of the other cancellation against a non-survivor.
Yet while the involutions defining such cancellations should allow restriction
to the survivors of the other cancellation, their definition should not be
limited to those survivors (as is the case for the Bender-Knuth involutions,
which are only defined for semistandard tableaux): they should be as general
as the Gessel-Viennot cancellations. We shall show in a forthcoming sequel to
this paper that both for (\LRgeninteq) and for~(\LRgenbineq) such a pair of
cancellations does indeed exist. Moreover the operations used to define these
cancellations are quite interesting in their own right, and they give rise to
ramifications that relate to Knuth correspondences and to jeu de taquin.

% Local IspellDict: default

\Supplement References.

\input \jobname.ref

\iftoc 
  \closeout\tocfile \vfil\eject 
  \def\tocitem#1=#2\onpage#3.
    {\par\line{\hbox to 1.5pc{\bf #1\hss}\quad#2\dotfill#3}\ignorespaces}
  \catcode`@=11
  \input \jobname.toc
\fi

\bye